\documentclass[%
 amsmath,amssymb,
 aps,
 prx,
]{revtex4-2}
\usepackage{graphicx}
\usepackage[usenames, dvipsnames]{xcolor}
\usepackage{amssymb, latexsym}
\usepackage{amsmath}
\usepackage{mathrsfs}
\usepackage{soul}
\usepackage{natbib}
\usepackage{amsbsy}
\usepackage{psfrag}
\usepackage{mathtools}
\usepackage{color}
\usepackage{subfig}
\usepackage{makecell}
\usepackage[makeroom]{cancel}
\definecolor{myBlue}{rgb}{0, 0.4470, 0.7410}
\definecolor{myRed}{rgb}{0.8500, 0.3250, 0.0980}

\usepackage{tikz}
\usetikzlibrary{decorations.markings,arrows.meta,calc,shapes,intersections,backgrounds,patterns.meta,bending,angles}
\usepackage{xifthen}

\usepackage{hyperref}
\hypersetup{
        breaklinks,
        colorlinks=true,
        linkcolor=blue,
        urlcolor=Purple,
        citecolor=Purple,
 }

\usepackage{gacaps}
\usepackage[normalem]{ulem}

\def\beq{\begin{equation}}
\def\eeq{\end{equation}}

\def\beq{\begin{equation}}
\def\eeq{\end{equation}}

\def\dis{\varepsilon}

\definecolor{olivegreen}{rgb}{0,0.6,0}

\def\drawline#1#2{\raise 2.5pt\vbox{\hrule width #1pt height #2pt}}

\def\trian{\raise 1.25pt\hbox{$\scriptstyle\triangle$}\nobreak}

\def\dtrian{\raise 1.25pt\hbox%
{$\scriptscriptstyle\bigtriangledown$}\nobreak}

\def\squar{\raise 1.25pt\hbox{$\scriptstyle\Box$}\nobreak}

\def\diamon{\raise 1.25pt\hbox{$\scriptstyle\diamond$}\nobreak}

\def\beq{\begin{equation}}
\def\eeq{\end{equation}}

\def\citalajim03{Del \'Alamo \& Jim\'enez (2003)}

\definecolor{C0}{HTML}{1F77B4} 
\definecolor{C1}{HTML}{FF7F0E}
\definecolor{C2}{HTML}{2CA02C}
\definecolor{C3}{HTML}{D62728}
\definecolor{C4}{HTML}{9467BD}
\definecolor{C5}{HTML}{8C564B}
\definecolor{C6}{HTML}{E377C2}
\definecolor{C7}{HTML}{7F7F7F}
\definecolor{C8}{HTML}{BCBD22}
\definecolor{C9}{HTML}{17BECF}

\let\bs\boldsymbol
\let\provc\providecommand
\newcommand{\mybQ}{\bs{Q}}        
\newcommand{\mybq}{\bs{q}}        
\newcommand{\mybY}{\bs{Y}}        
\newcommand{\mybql}{\bs{q}}        
\newcommand{\mybyl}{\bs{y}}        

\newcommand{\myindexvar}[3]{%
  \ifthenelse{\isempty{#2}\and\isempty{#3}}%
  {#1}{#1_{#2}^{#3}}} 

\newcommand{\Q}[2]{\myindexvar{Q}{#1}{#2}}       
\newcommand{\bQ}[2]{\myindexvar{\mybQ}{#1}{#2}}  

\newcommand{\bQtilde}[2]{\myindexvar{\widetilde{\mybQ}{\mkern 0mu}}{#1}{#2}}
\newcommand{\bQprime}[2]{\myindexvar{{\mybQ'}{\mkern 0mu}}{#1}{#2}} \newcommand{\bqtilde}[2]{\myindexvar{\widetilde{\mybql}}{#1}{#2}}

\newcommand{\bQhat}[2]{\myindexvar{\widehat{\mybQ}{\mkern 0mu}}{#1}{#2}}
\newcommand{\bYhat}[2]{\myindexvar{\widehat{\mybY}{\mkern 0mu}}{#1}{#2}}
\newcommand{\bYtilde}[2]{\myindexvar{\widetilde{\mybY}{\mkern 0mu}}{#1}{#2}}

\newcommand{\bqhat}[2]{\myindexvar{\widehat{\mybql}}{#1}{#2}}
\newcommand{\byhat}[2]{\myindexvar{\widehat{\mybyl}}{#1}{#2}}
\newcommand{\bytilde}[2]{\myindexvar{\widetilde{\mybyl}}{#1}{#2}}

\newcommand{\Y}[2]{\myindexvar{Y}{#1}{#2}}       
\newcommand{\bY}[2]{\myindexvar{\mybY}{#1}{#2}}  

\provc{\sq}[2]{\myindexvar{q}{#1}{#2}}       
\provc{\bsq}[2]{\myindexvar{\mybq}{#1}{#2}}  

\newcommand{\mybS}{\bs{S}}        
\newcommand{\mybA}{\bs{A}}        
\newcommand{\mybth}{\bs{\theta}}  
\newcommand{\mybJ}{\bs{J}}        

\provc{\bS}[2]{\myindexvar{\mybS}{#1}{#2}}     
\provc{\bA}[2]{\myindexvar{\mybA}{#1}{#2}}     
\provc{\bth}[2]{\myindexvar{\mybth}{#1}{#2}}   

\provc{\bWs}[2]{\myindexvar{\bs{W}}{#1}{#2}}               
\provc{\bWa}[2]{\myindexvar{\bs{V}}{#1}{#2}}   

\provc{\bJ}[2]{\myindexvar{\mybJ}{#1}{#2}}                  
\provc{\bJtilde}[2]{\myindexvar{\widetilde{\mybJ}}{#1}{#2}} 
\provc{\bJhat}[2]{\myindexvar{\widehat{\mybJ}}{#1}{#2}}     

\provc{\J}[2]{\myindexvar{J}{#1}{#2}}                  
\provc{\Jtilde}[2]{\myindexvar{\widetilde{J}}{#1}{#2}} 
\provc{\Jhat}[2]{\myindexvar{\widehat{J}}{#1}{#2}}     

\provc{\mun}{\mu}
\provc{\sgn}{\Xi}
\provc{\bmun}{\bs{\mu}}
\provc{\bsgn}{\bs{\Xi}}
\provc{\muopt}{\mu^*}
\provc{\sgopt}{\Xi^*}
\provc{\bmuopt}{\bs{\mu}^*}
\provc{\bsgopt}{\bs{\Xi}^*}
\provc{\mutar}{\hat{\mu}}
\provc{\sgtar}{\widehat{\Xi}}
\provc{\bmutar}{\hat{\bs{\mu}}}
\provc{\bsgtar}{\widehat{\bs{\Xi}}}

\provc{\relf}{\alpha}
\provc{\relfmu}{\relf_\mu}
\provc{\relfsg}{\relf_\xi}

\provc{\bths}{\mybth_s} 
\provc{\bthpa}{\mybth_p}
\provc{\bthaa}{\mybth_a}

\provc{\ths}{\theta_s} 
\provc{\thpa}{\theta_p}
\provc{\thaa}{\theta_a}

\provc{\shiftM}{\bs{a}}
\provc{\shiftmi}{a}

\definecolor{myc1}{HTML}{003049}
\definecolor{myc2}{HTML}{d62828}
\definecolor{myc3}{HTML}{f77f00}
\definecolor{myc4}{HTML}{6ca13b}

\newcommand{\myls}[2]{\lineSymbol[#1]{none}{#2}{.8ex}{.8pt}{#2}{#2}}

\let\ig\includegraphics
\let\tw\textwidth

\newcommand{\bse}{\begin {subequations}}
\newcommand{\ese}{\end {subequations}}
\newcommand {\ba} {\begin {array}}
\newcommand {\ea} {\end {array}}

\usepackage{pifont}
\begin{document}

\title{Information-theoretic formulation of dynamical systems: causality, modeling, and control}

%
\author{Adri\'an Lozano-Dur\'an}
\author{Gonzalo Arranz}%
\affiliation{%
Department of Aeronautics and Astronautics, Massachusetts Institute of Technology, Cambridge, MA 02139, USA
}%

\date{\today}

\begin{abstract}
    The problems of causality, modeling, and control for chaotic,
high-dimensional dynamical systems are formulated in the language of
information theory.  The central quantity of interest is the Shannon
entropy, which measures the amount of information in the states of the
system. Within this framework, causality is quantified by the
information flux among the variables of interest in the dynamical
system. Reduced-order modeling is posed as a problem related to the
conservation of information in which models aim at preserving the
maximum amount of relevant information from the original
system. Similarly, control theory is cast in information-theoretic
terms by envisioning the tandem sensor-actuator as a device reducing
the unknown information of the state to be controlled. The new
formulation is used to address three problems about the causality,
modeling, and control of turbulence, which stands as a primary example
of a chaotic, high-dimensional dynamical system.  The applications
include the causality of the energy transfer in the turbulent cascade,
subgrid-scale modeling for large-eddy simulation, and flow control for
drag reduction in wall-bounded turbulence.

\end{abstract}

\maketitle

\newcommand{\corr}[1]{\textcolor{black}{#1}}
\section{Introduction}
 
Information theory is the science about the laws governing information
or, in other words, the mathematics of message communication. A
message can be thought of as the bits (ones and zeros) of an image
transfer via the Internet, but also as the cascade of energy in a
turbulent flow or the drag reduction in an airfoil when applying a
particular control strategy. Information theory is one of the few
scientific fields fortunate enough to have an identifiable beginning:
\citet{shannon1948}, who ushered us into the Information Age with a
quantitative theory of communication. Since then, a field that started
as a branch of mathematics dealing with messages, ultimately matured
into a much broader discipline applicable to engineering, biology,
medical science, sociology,
psychology\ldots~\citep[e.g.,][]{cover2006}.  The success of
information theory relies on the notion of information as a
fundamental property of physical systems, closely tied to the
restrictions and possibilities of the laws of
physics~\citep{landauer1996}. The fundamental nature of information
provides the foundations for the principles of conservation of
information and maximum entropy highly regarded within the physics
community~\citep{landauer1991, brillouin2013,
  susskind2013}. Interestingly, despite the accomplishments of
information theory in many scientific disciplines, applications to
some branches of physics are remarkably limited. The goal of the
present work is to advance our physical understanding of chaotic,
high-dimensional dynamical systems by looking at the problems of
causality, reduced-order modeling, and control through the lens of
information theory.

The grounds for causality as information are rooted in the intimate
connection between information flux and the arrow of time: the laws of
physics are time-symmetric at the microscopic level, and it is only
from the macroscopic viewpoint that time-asymmetries arise in the
system~\citep{eddington1929}. Such asymmetries determine the direction
of time, that can be leveraged to measure the causality of events
using information-theoretic metrics based on the Shannon entropy.
Modeling can be posed as a problem of conservation of information:
reduced-order models contain a smaller number of degrees of freedom
than the original system, which in turn entails a loss of
information. Thus, the goal of a model is to guarantee as much as
possible the conservation of relevant information from the original
system.
Similarly, control theory for dynamical systems can be cast in
information-theoretic terms if we envision the tandem sensor-actuator
as a device aimed at reducing the unknown information associated with
the state of the system to be controlled. In all of the cases above,
the underlying idea advanced is that the evolution of the information
content in a chaotic system greatly aids the understanding and
manipulation of the quantities of interest.

\corr{In the present work, \emph{i)} a new definition of causality is
  proposed based on the information required to attain total knowledge
  of a variable in the future,
  \emph{ii)} the conditions for maximum information-preserving models are
  derived and leveraged to prove that accurate models maximize the
  information shared between the model state and the true state, and
  \emph{iii)} new definitions of open/closed-loop control,
  observability, controllability, and optimal control are introduced
  in terms of the information shared among the variable to control
  and/or the sensors and actuators.}
%
The information-theoretic formulation of causality, reduced-order
modeling and control is introduced in \S\ref{sec:causality},
\S\ref{sec:modeling}, and \S\ref{sec:control}, respectively.  The
sections are self-contained and follow a consistent notation. Each
section provides a brief introduction of the topic and closes with the
application of the theory to tackle a problem in turbulent
flows. Nonetheless, the theory is broadly applicable to any chaotic
dynamical system. Given our emphasis on turbulent flows, we provide
next a summary of current approaches for turbulence research. The
reader interested in the formulation of the theory is directly
referred to \S\ref{sec:definitions}.

\subsection{Abridged summary of modeling, control, and causality in turbulence research}

Turbulence, i.e., the multiscale motion of fluids, stands as a primary
example of a chaotic, high-dimensional phenomenon. Broadly speaking,
efforts in turbulence research can be subdivided into physical
insight, modeling, and control.  The three branches are intimately
intertwined, yet they provide a conceptual partition of the field
which is useful in terms of goals and methods.  Even a brief survey of
the methods for modeling, control, and causality would entail a
monumental task that will not be attempted here. Instead, common
pathways to tackle these problems are discussed along with the some
pitfalls and limitations.

In the context of modeling, the field of fluid mechanics is in the
enviable position of owning a set of equations that describes the
motion of a fluid to near-perfect accuracy: the Navier-Stokes
equations. Thus, significant ongoing
efforts are devoted to capturing the essential flow physics in the
form of reduced-order models. Prominent techniques include proper
orthogonal decomposition and Galerkin projection~\citep{berkooz1993},
balanced truncation and dynamic mode decomposition~\citep{schmid2007},
or extensions by Koopman theory~\citep{rowley2017}. Machine learning
also provides a modular and agile framework that can be tailored to
address reduced-order modeling~\citep{duraisamy2019, brunton2019}.
Linear theories are still instrumental for devising reduced-order
models, while other approaches rely on phenomenological arguments. The
modeling application of this work is centered on large-eddy simulation
(LES), in which the large eddies in the flow are resolved and the effect of the
small scales is modeled through a subgrid-scale model (SGS).  Most
SGS models are derived from a combination of theory and physical
intuition.~\citep{meneveau2000, piomelli2002, bose2018}.  In addition,
Galilean invariance, along with the principles of mass, momentum, and
energy conservation, are invoked to constrain the admissible
models~\citep[e.g.,][]{speziale1991, silvis2017}.
However, although we do possess a crude practical
understanding of turbulence, flow predictions  from
the state-of-the-art models are still unable to comply with the
stringent accuracy requirements and computational efficiency demanded
by the industry~\citep{slotnick2014}.

Control, the ability to alter flows to achieve the desired outcome, is
a matter of tremendous consequence in engineering. In system control,
sensors measure the state of the flow, while actuators create the flow
disturbances to prevent or trigger a targeted condition (e.g., drag
reduction, mixing enhancement, etc.).  Recent decades have seen a
flourishing of activity in various techniques for control of turbulent
flows --active, passive, open-loop, closed-loop~\citep{gadelhak2000,
  bewley2001, gunzburger2002, kim2003, collis2004, brunton2019}. A
common family of methods originates from linear theories, which
constitutes the foundation of many control strategies~\citep{kim2007,
  schmid2012, mckeon2017, rowley2017, zare2020,
  jovanovic2020}. However, linear methods have sparked criticism as
turbulence is a highly nonlinear phenomenon, and universal control
strategies cannot be anticipated from a single set of linearized
equations. Nonlinear control strategies are less common, but they have
also been available for years~\citep[e.g.,][]{king2005,
  luchtenburg2010, aleksic2010}.  They are, nonetheless, accompanied
by a considerable penalty in the computational cost, which renders
nonlinear control impractical in many real-world applications.

Causality is the mechanism by which one event contributes to the
genesis of another~\citep{pearl2009}. Whereas control and modeling are
well-established cornerstones of turbulence research, the same cannot
be said about the elusive concept of causality, which has been
overlooked within the fluids community except for a handful of
works~\citep{tissot2014, liang2016, lozano2019b, lozano2021}.  In the
case of turbulence research, causal inference is usually simplified in
terms of the cross-time correlation between pairs of time signals
representing the events of interest (e.g., kinetic energy,
dissipation, etc.). The correlation method dates back to the work of
the mathematician A.-M. Legendre in the 1800s, and undoubtedly
constitutes an outdated legacy tool. Efforts to infer causality using
time-correlation include the investigations of the turbulent kinetic
energy~\citep{jimenez2018, cardesa2015} and the space-time signature
of spectral quantities~\citep[e.g.,][]{choi1990, wallace2014,
  wilczek2015, kat2015, he2017, wang2020}, to name a few. However, it
is universally accepted that correlation does not imply
causation~\citep{beebee2012}, as the former lacks the directionality
and asymmetry required to quantify causal interactions. Despite this
limitation, the correlation between time signals stands as the
state-of-the-art tool for (non-intrusive) causality quantification
in fluid mechanics.

The goal of the present work is to further advance the field of
turbulence research by introducing a new information-theoretic
formalism for causality, modeling, and control.  To date, the use of
information-theoretic tools in the fluid mechanics community is still
in its infancy. \citet{betchov1964} was one of first authors to
propose an information-theoretic metrics to quantify the intricacy of
turbulence.  \citet{cerbus2013} applied the concept of conditional
Shannon entropy to analyze the energy cascade in 2-D turbulence. The
work by \citet{cerbus2014} also contains additional discussions on the
use of established tools in information theory for fluid
dynamics. \corr{\citet{materassi2014} used normalized transfer entropy
  to study the cascading process in synthetic turbulence generated via
  the shell model. In a series of works, \citet{granero-belinchon2016,
    granero-belinchon2018, granero-belinchon_thesis,
    granero-belinchon2021} investigated the information content,
  intermittency, and stationarity characteristics of isotropic
  turbulence using information-theoretic tools applied to experimental
  velocity signals}. \citet{liang2016} and \citet{lozano2019b} applied
information-theoretic definitions of causality to unveil the dynamics
of energy-containing eddies in wall-bounded turbulence.  A similar
approach was followed by \citet{wang2021} to study cause-and-effect
interactions in turbulent flows over porous media.
\citet{lozano2019b} also discussed the use of information transfer
among variables to inform the design of reduced-order
models. \citet{shavit2020} used singular measures and information
capacity to study the turbulent cascade and explore the connection
between information and modeling.  More recently, \citet{lee2021}
leveraged the principle of maximum-entropy to analyze the turbulence
energy spectra. The aforementioned studies have offered a new
perspective of turbulence using established tools in information
theory. In the following, we further develop the theory of information
for dynamical systems with the aim of advancing the field of
turbulence research.

\section{Basics of information theory}
\label{sec:definitions}

Let us introduce the concepts of information theory required to
formulate the problems of causality, modeling, and control.  The first
question that must be addressed is the meaning of \emph{information},
as it departs from the intuitive definition used in our everyday
life. Let us consider the discrete random variable $X$ taking values
equal to $x$ with probability mass function $p(x)=\mathrm{Pr}\{X =
x\}$ over the finite set of outcomes of $X$.  
The information of observing the event $X=x$ is defined
\corr{as~\citep{shannon1948}}:
\begin{equation}\label{eq:info}
  \mathcal{I}(x) = -\log_2[p(x)].
\end{equation}
The units of $\mathcal{I}(x)$ are set by the base chosen, in this case
`bits' for base 2. The base of the logarithm is arbitrary and can be
changed using the identity $\log_a p = \log_a b \log_b p$. For
example, consider tossing a fair coin with
$X\in\{\mathrm{heads},\mathrm{tails}\}$ such that
$p(\mathrm{heads})=p(\mathrm{tails})=0.5$. The information of getting
heads after flipping the coin once is $\mathcal{I}(\mathrm{heads}) =
-\log_2(0.5) = 1$ bit, i.e., observing the outcome of flipping a fair
coin provides one bit of information. If the coin is completely biased
towards heads, $p(\mathrm{heads})=1$, then
$\mathcal{I}(\mathrm{heads}) = -\log_2(1) = 0$ bits (where $0\log 0 =
0$), i.e., no information is gained as the outcome was already known
before flipping the coin. This simple but revealing example
illustrates the meaning of information in the present work:
information is the statistical notion of how unlikely it is to observe
an event.  Low probability events provide more information than high
probability events. Thus, within the framework of information theory,
the statement ``tomorrow the sun will rise in the west'' contains more
information than ``tomorrow the sun will rise in the east'', simply
because the former is extremely unlikely. The information can also be
interpreted in terms of uncertainty: $\mathcal{I}(x)$ is the number of
bits required to unambiguously determine the state $x$. The latter
interpretation will be frequently evoked in this work.

The reader might ask why not choosing information to be directly
proportional to $p(x)$ rather than to $-\log_2[p(x)]$. However, the
logarithm function is the most natural choice for a measure of
information that is additive in the number of states of the system
considered. This might be illustrated by tossing a fair coin $n$
times. The information gathered for a particular sequence of events is
\begin{equation}\label{eq:info_n}
  \mathcal{I}(\mathrm{heads},\mathrm{heads},\mathrm{tails},...) = -\log_2(0.5^n) =  n \ \mathrm{bits}.
\end{equation}
In general, when two systems with $N$ different states are combined,
the resulting system contains $N^2$ states (i.e., Cartesian product of
the states of both systems), but the amount of information is
$2N$~\citep{mackay2002} as illustrated in the example above. Another
viewpoint of Eq. (\ref{eq:info_n}) is that the probability of
observing a sequence of events is a multiplicative process given by
$p(x_1)p(x_2)...p(x_N)$, whereas it would be preferable to work with
an additive process. The latter is attained by taking the logarithms
of the probabilities, $-\log_2[p(x_1)] -
\log_2[p(x_2)]-...-\log_2[p(x_N)]$, where the minus sign is introduced
for convenience to obtain an outcome that is equal or larger than
zero.

Equations (\ref{eq:info}) and (\ref{eq:info_n}) provide the
information gained observing one particular event or a sequence of
events, respectively. Usually, we are interested in the average
information in $X$ given by the expectation $\langle \cdot \rangle$
over all the possible outcomes
\begin{equation}\label{eq:entropy}
   H(X) = \langle \mathcal{I}(x) \rangle = \sum_x -p(x) \log_2[p(x)]\geq 0, 
\end{equation}
which is referred to as the Shannon entropy and represents the
generalization to arbitrary variables of the well-known thermodynamic
entropy~\citep{boltzmann1877,jaynes1957}.  Following the example
above, the entropy of the system ``flipping a fair coin $n$ times'' is
$H = -\sum 0.5^n \log_2(0.5^n) = n$ bits, where the sum is performed
over all the possible outcomes of flipping a coin $n$ times (namely,
$2^n$). As expected, flipping $n$ times a biased coin with
$p(\mathrm{heads})=1$ provides no information ($H=0$).
\citet{shannon1948} showed that Eq.~(\ref{eq:entropy}) corresponds to
the minimum average number of bits needed to encode a source of $n$
states with probability distribution $p$.  As measure of uncertainty,
$H$ is maximum when all the possible outcomes are equiprobable (large
uncertainty in the state of the system) and zero when the process is
completely deterministic (no uncertainty in the outcome).

The Shannon entropy can be generalized to $m$ random variables
$\boldsymbol{X} = [X_1,X_2,...,X_m]$ as
\begin{subequations}
  \begin{gather}
H(\boldsymbol{X}) = \langle \mathcal{I}(\boldsymbol{x}) \rangle = 
\sum_{x_1,...,x_m} -p(x_1,x_2,...,x_m) \log_2[p(x_1,x_2,...,x_m)]
\end{gather}
\end{subequations}
where $p(x_1,x_2,...,x_m)$ is the joint probability mass function
$\mathrm{Pr}\{X_1 = x_1, X_2 = x_2,...,X_m = x_m\}$. Similarly, given
the random variables $X$ and $Y$ and the conditional distribution
$p(x|y) = p(x,y)/p(y)$ with $p(y) = \sum_x p(x,y)$ as the marginal
probability distribution of $Y$, the entropy of $X$ conditioned on $Y$
is defined as~\citep{stone2013}
\begin{equation}
H( X | Y ) = \sum_{x,y} -p(x,y) \log_2[p(x|y)].
\end{equation}
It is useful to interpret $H(X|Y)$ as the uncertainty in the state $X$
after conducting the `measurement' of the state $Y$. This
interpretation is alluded to in the following sections. If $X$ and $Y$
are independent random variables, then $H(X|Y)=H(X)$, i.e., knowing
the state $Y$ does not reduce the uncertainty in $X$. Conversely,
$H(X|Y)=0$ if knowing $Y$ implies that $X$ is completed
determined. Finally, the mutual information between the random
variables $X$ and $Y$ is
\begin{equation}\label{eq:optimal}
   I(X;Y) =  H(X) - H(X|Y) = H(Y) - H(Y|X),
\end{equation}
which is a symmetric metric $I(X;Y)=I(Y;X)$ representing the
information shared among the state variables $X$ and $Y$.  The mutual
information between variables will be also central to the formalism
presented below. Figure \ref{fig:basics} depicts the relationship
between the entropy, conditional entropy, and mutual information.

The concepts above provide the foundations to the
information-theoretic formulation of causality, modeling, and control
detailed in the following sections.  First, we introduce the formalism
of information in dynamical systems.  Several studies have already
discussed this topic, mostly in the context of predictability and
chaos~\citep[e.g.][]{shaw1981, delsole2004, garbaczewski2006,
  liang2005, kleeman2011}.  The exposition here is extended and
tailored to our particular interests.
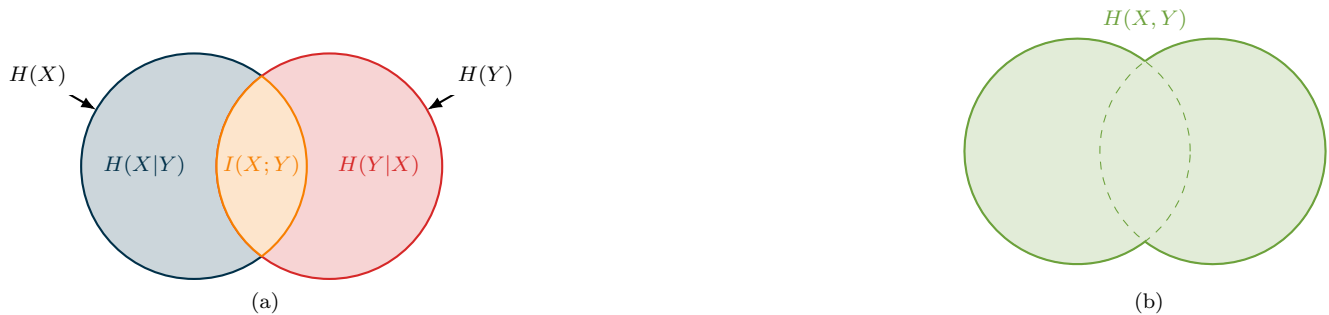
\begin{figure}
    \subfloat[]{\begin{tikzpicture}[scale=1.5,thick,>={Latex[length=.2cm]}]
    \colorlet{c1}{myc1}
    \colorlet{c2}{myc2}

    \colorlet{MI}{myc3}
    \pgfmathsetmacro{\dis}{.6}
    \begin{footnotesize}
        \draw[c1,fill=c1!20] (-\dis,0) circle (1);
        \draw[c2,fill=c2!20] (\dis,0) circle (1);
        
        \begin{scope}
            \clip (-\dis,0) circle (1);
            \draw[draw=MI,fill=MI!20] (\dis,0) circle (1);
        \end{scope}
        \begin{scope}
            \clip (\dis,0) circle (1);
            \draw[MI] (-\dis,0) circle (1);
        \end{scope}

        \node[anchor=east,c1] at (-\dis,0) {$H(X|Y)$};
        \node[anchor=west,c2] at (+\dis,0) {$H(Y|X)$};

        \node[anchor=center,MI] at (0,0) {$I(X;Y)$};

        \draw[<-] (+\dis,0) ++ (30:1)  --+ (30:.6)  node[fill=white,pos=1] {$H(Y)$};
        \draw[<-] (-\dis,0) ++ (150:1) --+ (150:.6) node[fill=white,pos=1] {$H(X)$};

    \end{footnotesize}
\end{tikzpicture}} \hfill
    \subfloat[]{\begin{tikzpicture}[scale=1.5,thick]
    \pgfmathsetmacro{\dis}{.6}
    \begin{footnotesize}

        \colorlet{c1}{myc4}
        
        \fill[c1!20!white] (-\dis,0) circle (1);
        \fill[c1!20!white] (\dis,0) circle (1);
        
        \draw[dashed,c1,thin] (-\dis,0) circle (1);
        \draw[dashed,c1,thin] (\dis,0) circle (1);

        \pgfmathsetmacro{\Angle}{acos(\dis)};
        \pgfmathsetmacro{\inters}{sqrt(1 - \dis^2)};

        \draw[c1] (0,\inters) arc (\Angle:360-\Angle:1) --
                     (0,-\inters) arc (-180+\Angle:180-\Angle:1);


        \node[anchor=center,c1,anchor=south] at (0,1) {$H(X,Y)$};

    \end{footnotesize}
\end{tikzpicture}}
    \caption{Venn diagrams of (a) the conditional entropy and mutual
      information between two random variables $X$ and $Y$, and (b)
      total entropy in two random variables. \label{fig:basics}}
\end{figure}

\section{Information in dynamical systems}
\label{sec:information_dynamical}

Let us consider the continuous \corr{deterministic} dynamical system
with state variable $\bsq{}{}=\bsq{}{}(\boldsymbol{x},t)$, where
$\boldsymbol{x}$ is the vector of spatial coordinates, and $t$ is
time. The dynamics of $\bsq{}{}$ are governed by the partial
differential equation
\begin{equation}
  \label{eq:cau:q_cont}
\frac{\partial\bsq{}{}}{\partial t} =
\boldsymbol{F}(\bsq{}{}),
\end{equation}  
that might represent, for example, the equations of conservation of
mass, momentum, and energy.  Equation (\ref{eq:cau:q_cont}) can be
integrated in time from $t_n$ to $t_{n+1}$ to yield
\begin{equation}
  \label{eq:cau:qn_cont}
\bsq{}{}(\bs{x},t_{n+1})= \bsq{}{}(\bs{x},t_n) + \int_{t_n}^{t_{n+1}}\boldsymbol{F}(\bsq{}{}) \mathrm{d}t,
\end{equation}
where $t_{n+1}-t_n$ is an arbitrary time span.
%
We consider a spatially coarse-grained version of $\bsq{}{}$ at time
$t_n$ denoted by $\bsq{}{n} = [\sq{1}{n},....,\sq{N}{n}]$, where $N$
is the number of degrees of freedom of the system. We assume that the
dimensionality of $\bsq{}{n}$ is large enough to capture all the
relevant dynamics of Eq. (\ref{eq:cau:qn_cont}). In the context of
fluid dynamics $\bsq{}{n}$ might represent, for example, the three
velocities components and pressure at discrete spatial locations, the
Fourier coefficients of the velocity, the coefficients from the
Karhunen-Lo{\`e}ve decomposition of the flow~\citep{berkooz1993}, or
in general, any spatially-finite representation of the continuous
system.

We treat $\bsq{}{n}$ as a random variable, indicated by
$\bQ{}{n}=[\Q{1}{n},\dots,\Q{N}{n}]$, and consider a finite partition
of the phase space $D =\{D_1, D_2, \dots, D_{N_q}\}$, \corr{where
  $N_q$ is the number of partitions,} such that $ D = \cup_{i=1}^{N_q}
D_i$ and $D_i \cap D_j = \emptyset$ for all $i \neq j$. The system is
said to be in the state $D_i$ if $\corr{\bsq{}{n}} \in D_i$. The
probability of finding the system at state $D_i$ at time $t_n$ is
$p^q_i = \mathrm{Pr}\{\bQ{}{n} \in D_i\}$.  For simplicity, we refer
to the latter probability simply as $p(\bsq{}{n})$.  The dynamics of
$\bQ{}{n}$, are determined by
\begin{equation}
  \label{eq:cau:Q_total}
  \bQ{}{n+1} = \boldsymbol{f}(\bQ{}{n}),
\end{equation}
where the map $\boldsymbol{f}$ is derived from
Eq. (\ref{eq:cau:qn_cont}). \corr{Note that the system considered in
Eq. (\ref{eq:cau:Q_total}) is closed in the sense that no external
stochastic forcing is applied.}

The information contained in the system at time $t_n$ is given by the
entropy of $\bQ{}{n}$, namely, $H(\bQ{}{n})$. As the system evolves in
time according to Eq. (\ref{eq:cau:Q_total}), its information content
is bounded by
\begin{equation}
  \label{eq:cau:H_reduction}
  H(\bQ{}{n+1}) = H(\boldsymbol{f}(\bQ{}{n})) \leq H(\bQ{}{n}),
\end{equation}
which is the result of the entropy of transformed random
variables~\citep{cover2006}.  A consequence of
Eq. (\ref{eq:cau:H_reduction}) is that the dynamical system in
Eq. (\ref{eq:cau:Q_total}) either conserves or destroys information,
but never creates information. Another interpretation of
Eq. (\ref{eq:cau:H_reduction}) is that, for deterministic systems, the
information of future states is completely determined by the initial
condition, whereas the converse is not always true.  For example,
dissipative systems cannot be integrated backwards in time to
univocally recover its initial state.

The entropy of the system at $t_{n}$ can be related to the entropy at
$t_{n+1}$ through the Perron-Frobenious operator
$\mathbb{P}[\cdot]$~\citep{beck1995} which advances the probability
distribution of the system
\begin{equation}
  \label{eq:cau:Perron}
  p(\bsq{}{n+1}) = \mathbb{P}[ p(\bsq{}{n}) ].
\end{equation}
By construction of the system in Eq. (\ref{eq:cau:Q_total}), we can
derive the zero conditional-entropy condition
\begin{subequations}
\label{eq:caus:H_total}
  \begin{gather}
  H(\bQ{}{n+1} | \bQ{}{n} ) =  
  \sum -p( \bs{q}^{n+1}, \bs{q}^n ) \log[p( \bs{q}^{n+1} | \bs{q}^n)] = \\
 = \sum -\mathbb{P}[p( \bsq{}{n}| \bsq{}{n} )] p(\bsq{}{n} ) \log\{ \mathbb{P}[p( \bsq{}{n}| \bsq{}{n} )]\} =0,
\end{gather}
\end{subequations}
which shows that there is no uncertainty in the future state
$\bQ{}{n+1}$ given the past state $\bQ{}{n}$. Equation
(\ref{eq:caus:H_total}) merely echoes the deterministic nature of the
governing equations, and will be instrumental in the formulation of
the principles for causality, modeling, and control.  Additionally, if
the map $\boldsymbol{f}$ is reversible, namely, $\bQ{}{n} =
\boldsymbol{f}^{-1}(\bQ{}{n+1})$, then we obtain the conservation of
information for dynamical systems
\begin{equation}
  \label{eq:cau:conservation_info}
 H(\bQ{}{n+1}) = H(\bQ{}{n}),
\end{equation}
which can be regarded as a fundamental principle underlying the rest of conservation laws.

The condition in Eq. (\ref{eq:caus:H_total}) may be generalized by
adding the noise, $\bs{W}^n$, which accounts for uncertainties in the
system state $\bQ{}{n}$, numerical errors, unknown physics in the map
$\bs{f}$, etc. The new governing equation is then $\bQ{}{n+1} =
\boldsymbol{f}(\bQ{}{n},\bs{W}^n)$ which implies that
\corr{$H(\bQ{}{n+1} | \bQ{}{n} ) \geq 0$ (information can be created)}
unless the effect of noise is taken into consideration,
i.e. $H(\bQ{}{n+1} | \bQ{}{n}, \boldsymbol{W}^n ) = 0$. \corr{A
  consequence of the latter is that for long integration times in
  chaotic systems, a small amount of noise will also result in
  $H(\bQ{}{n+1} | \bQ{}{n} ) \geq 0$}. Hereafter, we center our
attention on fully deterministic systems and assume that the impact of
the noise on $\bQ{}{n+1}$ is negligible
($\boldsymbol{W}^n=\boldsymbol{0}$) for the problems of causality and
modeling. The effect of the noise will be introduced in the
formulation of control.

\section{Information flux as causality}
\label{sec:causality}

\corr{One of the most intuitive definitions of causality relies on the
  concept of interventions: manipulation of the causing variable leads
  to changes in the effect~\citep{pearl2009,
    eichler2013}. Interventions provide a pathway to evaluate the
  causal effect that a process $A$ exerts on another process $B$ by
  setting $A$ to a modified value $\widetilde{A}$ and observing the
  post-intervention consequences on $B$. Despite the intuitiveness of
  interventions as a measure of causality, the approach is not free of
  shortcomings. Causality with interventions is intrusive (i.e., it
  requires modification of the system) and costly (the simulations
  need to be recomputed if numerical experiments are used). When the
  data are collected from physical experiments, the causality with
  interventions might be even more challenging or directly impossible
  to practice (for instance, we cannot use interventions to assess the
  causality of the prices in the stock market in 2008). Causality with
  interventions also poses the question of what type of intervention
  must be introduced in $A$ and whether that would affect the outcome
  of the exercise as a consequence of forcing the system out of its
  natural attractor. The framework of information theory provides an
  alternative, non-intrusive definition of causality as the
  information transferred from the variable $A$ to the variable $B$.}
\corr{The idea can be traced back to the work of \citet{wiener1956}
  and was first quantified by \citet{granger1969} using signal
  forecasting via linear autoregressive models. In the context of
  information theory, the definition was formalized by
  \citet{massey1990} and \citet{kramer1998} through the use of
  conditional entropies with the so-called directed information.}
\citet{schreiber2000} introduced an heuristic definition of causality
inspired by the direction of information in Markov chains.
\citet{liang2006} and later \citet{sihna2016} proposed to infer
causality by measuring the information changes in the disturbed
dynamical system. \corr{The new formulation of causality proposed here
  is motivated by the information required to attain total knowledge
  of a future state. Similar to previous works, the definition relies
  on conditional entropies. However, our information-theoretic
  quantification of causality is directly grounded on the zero
  conditional-entropy condition for deterministic systems (i.e.,
  Eq. \ref{eq:caus:H_total}) and generalizes previous definitions of
  causality to multivariate systems. We also introduce the concept of
  information leak as the amount of information unaccounted for by the
  observable variables}.

\subsection{Formulation}
\label{sec:cau:info_flux}


\subsubsection{Information flux}
\label{sec:information_flux}
 
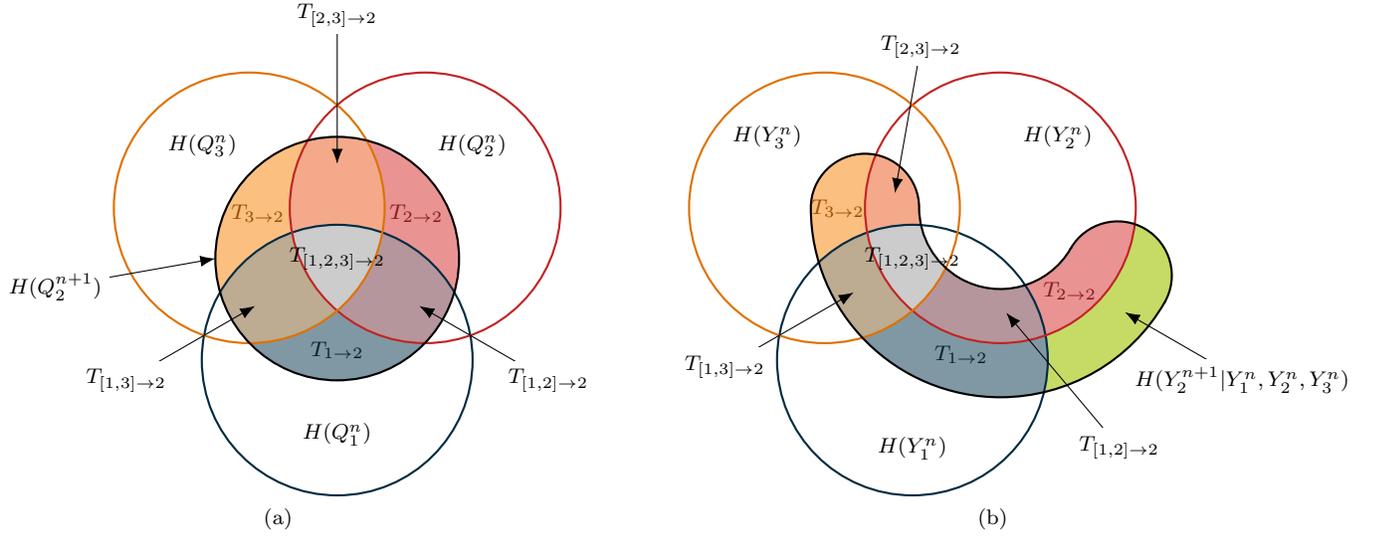
\begin{figure}
    \subfloat[\label{fig:cau:TQ}]{

\colorlet{c1}{myc1}
\colorlet{c2}{myc2}
\colorlet{c3}{myc3}

\colorlet{c1s}{c1!50}
\colorlet{c2s}{c2!50}
\colorlet{c3s}{c3!50}

\begin{tikzpicture}[scale=1.8,thick,>={Latex[length=.2cm]}
]
    
    \begin{footnotesize}
    \pgfmathsetmacro{\rdi}{1}
    \pgfmathsetmacro{\rdiO}{.9*\rdi}
    \pgfmathsetmacro{\dis}{.75}

    \coordinate (O1) at (270:\dis);
    \coordinate (O2) at (30:\dis);
    \coordinate (O3) at (150:\dis);
   
    \fill [c1s]  (210:\rdiO) arc (210:330:\rdiO) -- (0,0) -- cycle;
    \fill [c2s] (90:\rdiO) arc (90:-30:\rdiO) -- (0,0) -- cycle;
    \fill [c3s]    (90:\rdiO) arc (90:210:\rdiO) -- (0,0) -- cycle;

    \begin{scope} 
        \clip (O1) circle(\rdi);
    	\clip (O2) circle(\rdi);
        \path [name path=int12,draw=none, fill=c1s!50!c2s] (0:0) circle (\rdiO);
    \end{scope} 

    \begin{scope} 
    	\clip (O2) circle(\rdi);
        \clip (O3) circle(\rdi);
        \path [name path=int23,draw=none, fill=c2s!50!c3s] (0:0) circle (\rdiO);
    \end{scope} 
    
    \begin{scope} 
    	\clip (O1) circle(\rdi);
        \clip (O3) circle(\rdi);
        \path [name path=int13,draw=none, fill=c1s!50!c3s] (0:0) circle (\rdiO);
    \end{scope} 

    \begin{scope} 
    	\clip (O1) circle(\rdi);
        \clip (O3) circle(\rdi);
        \clip (O2) circle(\rdi);
        \path [name path=int0,draw=none, fill=white!80!black] (0:0) circle (\rdiO);
    \end{scope} 

    \draw [name path=C1,draw=c1!90!black] (O1) circle (\rdi);
    \draw [name path=C2,draw=c2!90!black] (O2) circle (\rdi);
    \draw [name path=C3,draw=c3!90!black] (O3) circle (\rdi);
    \draw [name path=C0,draw=black] (0:0) circle (\rdiO);

    \path (0,0) -- (O1) node[c1!60!black,pos=.9] {$T_{1\rightarrow2}$};
    \path (0,0) -- (O2) node[c2!60!black,pos=.9] {$T_{2\rightarrow2}$};
    \path (0,0) -- (O3) node[c3!60!black,pos=.9] {$T_{3\rightarrow2}$};

    \draw[thin,<-] (330:.7*\rdi) --++ (330:1.1) node[fill=white]{$T_{[1,2]\rightarrow2}$};
    \draw[thin,<-] (210:.7*\rdi) --++ (210:1.1) node[fill=white]{$T_{[1,3]\rightarrow2}$};
    \draw[thin,<-] (90:.7*\rdi) --++ (90:1.1)   node[fill=white]{$T_{[2,3]\rightarrow2}$};

    \node at( 0,0) {$T_{[1,2,3]\rightarrow2}$};
   
    \node at (270:1.3*\rdi) {$H(\Q{1}{n})$};
    \node at (40:1.3*\rdi)  {$H(\Q{2}{n})$};
    \node at (140:1.3*\rdi) {$H(\Q{3}{n})$};

    \draw[thin,<-] (180:\rdiO) --++ (190:1.2*\rdi) node[fill=white] {$H(\Q{2}{n+1})$};

    \end{footnotesize}
\end{tikzpicture}}~\hfill~
    \subfloat[\label{fig:cau:TY}]{

\colorlet{c1}{myc1}
\colorlet{c2}{myc2}
\colorlet{c3}{myc3}

\colorlet{c1s}{c1!50}
\colorlet{c2s}{c2!50}
\colorlet{c3s}{c3!50}

\colorlet{cHY}{SpringGreen}

\begin{tikzpicture}[scale=1.8,thick,>={Latex[length=.2cm]}]
    
    \begin{footnotesize}
    \pgfmathsetmacro{\rdi}{1}
    \pgfmathsetmacro{\dis}{.75}

    \pgfmathsetmacro{\rdn}{.4}
    \pgfmathsetmacro{\angS}{180}
    \pgfmathsetmacro{\angE}{330}

    \coordinate (O1) at (270:\dis);
    \coordinate (O2) at (30:\dis);
    \coordinate (O3) at (150:\dis);

    \fill[cHY] ($(O2)+(\angS:\rdi-\rdn)$) arc (\angS-180:\angS:\rdn)
          --+ (\angS:0) arc (\angS:\angE:\rdi+\rdn)
          --+ (\angE:0) arc (\angE:\angE+180:\rdn)
          --+ (0:0) arc (\angE:\angS:\rdi-\rdn);

    \foreach \i/\ci in {1/c1s,2/c2s,3/c3s} { 
    \begin{scope} 
        \clip ($(O2)+(\angS:\rdi-\rdn)$) arc (\angS-180:\angS:\rdn)
             --+ (\angS:0) arc (\angS:\angE:\rdi+\rdn)
             --+ (\angE:0) arc (\angE:\angE+180:\rdn)
             --+ (0:0) arc (\angE:\angS:\rdi-\rdn);
        \clip (O\i) circle(\rdi);
        \path [name path=int\i,draw=none, fill=\ci] (O\i) circle (\rdi);
    \end{scope} }
    
    \foreach \i/\j/\ci/\cj in {1/2/c1s/c2s,1/3/c1s/c3s,2/3/c2s/c3s} {
    \begin{scope} 
        \clip ($(O2)+(\angS:\rdi-\rdn)$) arc (\angS-180:\angS:\rdn)
             --+ (\angS:0) arc (\angS:\angE:\rdi+\rdn)
             --+ (\angE:0) arc (\angE:\angE+180:\rdn)
             --+ (0:0) arc (\angE:\angS:\rdi-\rdn);
    	\clip (O\i) circle(\rdi);
        \clip (O\j) circle(\rdi);
        \path [name path=int\i\j,draw=none, fill=\ci!50!\cj] (O\i) circle (\rdi);
    \end{scope}
    }
    
    \begin{scope} 
        \clip ($(O2)+(\angS:\rdi-\rdn)$) arc (\angS-180:\angS:\rdn)
             --+ (\angS:0) arc (\angS:\angE:\rdi+\rdn)
             --+ (\angE:0) arc (\angE:\angE+180:\rdn)
             --+ (0:0) arc (\angE:\angS:\rdi-\rdn);
    	\clip (O2) circle(\rdi);
    	\clip (O3) circle(\rdi);
        \clip (O1) circle(\rdi);
        \path [name path=int123,draw=none, fill=white!80!black] (O1) circle (\rdi);
    \end{scope} 

    \draw [name path=C1,draw=c1!90!black] (O1) circle (\rdi);
    \draw [name path=C2,draw=c2!90!black] (O2) circle (\rdi);
    \draw [name path=C3,draw=c3!90!black] (O3) circle (\rdi);
    
    \draw[thick] ($(O2)+(\angS:\rdi-\rdn)$) arc (\angS-180:\angS:\rdn)
          --+ (\angS:0) arc (\angS:\angE:\rdi+\rdn)
          --+ (\angE:0) arc (\angE:\angE+180:\rdn)
          --+ (0:0) arc (\angE:\angS:\rdi-\rdn);

    \pgfmathparse{.5*(\angE+\angS)}
    \path (O2) --+ (\pgfmathresult:\rdi+\rdn)   
        node[c1!60!black,pos=.8] {$T_{1\rightarrow2}$};

    \pgfmathparse{.5*(290+\angE)}
    \path (O2) --+ (\pgfmathresult:\rdi-\rdn/2) 
        node[c2!60!black,pos=1] {$T_{2\rightarrow2}$};
    \path (O2) --+ (180:\rdi+\rdn/2) node[c3!60!black,pos=1] {$T_{3\rightarrow2}$};

    \draw[thin,<-] (330:.8*\rdi) --++ (310:1.3) node[fill=white]{$T_{[1,2]\rightarrow2}$};
    \draw[thin,<-] (210:.5*\rdi) --++ (210:1.1) node[fill=white]{$T_{[1,3]\rightarrow2}$};
    \draw[thin,<-] (105:.5*\rdi) --++ (80:1.1)   node[fill=white]{$T_{[2,3]\rightarrow2}$};

    \node at (0,0) {$T_{[1,2,3]\rightarrow2}$};
   
    \node at (270:1.4*\rdi) {$H(\Y{1}{n})$};
    \node at (40:1.4*\rdi)  {$H(\Y{2}{n})$};
    \node at (140:1.4*\rdi) {$H(\Y{3}{n})$};

    \draw[thin,<-] ($(O2)+(320:\rdi+\rdn/2)$)--++(330:1*\rdi) node [fill=white]
        {$H(\Y{2}{n+1}|\Y{1}{n},\Y{2}{n},\Y{3}{n})$};

    \end{footnotesize}
\end{tikzpicture}}
    \caption{(a) Schematic of the entropies at time $t_n$ for a system
      with three variables $[Q^n_1,Q^n_2,Q^n_3]$ and their relation to
      the entropy of the future variable $Q^{n+1}_2$. Note that the
      entropy of $Q^{n+1}_2$ must be contained within the entropy of
      $H(Q^n_1,Q^n_2,Q^n_3)$ for consistency with
      Eq.~\eqref{eq:cau:HQ}.
      \corr{(b) Similar to (a), but for a system of observables states
        (Eq.~\eqref{eq:cau:Y}).  In this case, the entropy of
        $Y_{2}^{n+1}$ in not contained within $H(Y^n_1,Y^n_2,Y^n_3)$,
        leading to $T^Y_{\mathrm{leak},j} = H(Y_{2}^{n+1} | Y_{1}^{n},
        Y_{2}^{n}, Y_{3}^{n}) > 0$. }}
\end{figure}

The goal of this section is to leverage the information flux from
present states of the system to future states as a proxy for causal
inference.  Without loss of generality, let us derive the information
transferred from $\bQ{}{n}$ to $\Q{j}{n+1}$.  The dynamics of
$\Q{j}{n+1}$ is governed by the $j$-th component of
Eq. (\ref{eq:cau:Q_total}),
\begin{equation}
  \label{eq:cau:Qj}
  \Q{j}{n+1} = f_j(\bQ{}{n}).
\end{equation}
From Eq. (\ref{eq:cau:Qj}) and the propagation of information in
deterministic systems (Eq. (\ref{eq:caus:H_total})), it follows that
\begin{equation}
\label{eq:cau:HQ}
H(\Q{j}{n+1} | \bQ{}{n} ) =  0, 
\end{equation}
which shows that all the information contained in $\Q{j}{n+1}$
originates from $\bQ{}{n}$. Let us define the subset of variables
$\bQ{\bs{\bar\imath}}{n} =
[\bQ{\bar\imath_1}{n},... ,\bQ{\bar\imath_M}{n}]$, where
$\bs{\bar\imath}=[\bar\imath_1,...,\bar\imath_M]$ is a vector of
indices with $M\leq N$, and the vector of remaining variables
$\bQ{\cancel{\bs{\bar\imath}}}{n}$, such that $\bQ{}{n} =
[\bQ{\cancel{\bs{\bar\imath}}}{n},\bQ{\bs{\bar\imath}}{n}]$.  If only
the information from $\bQ{\cancel{\bs{\bar\imath}}}{n}$ is accessible,
then the uncertainty in the future state $\Q{j}{n+1}$ can be non-zero,
\begin{equation}
\label{eq:cau:H_part}
H(Q^{n+1}_j | \boldsymbol{Q}_{\cancel{\bs{\bar\imath}}}^n ) \geq  0. 
\end{equation}
%
%
Equation~\eqref{eq:cau:H_part} quantifies the average number of bits
required to completely determine\corr{\sout{d}} the state of $Q_j^{n+1}$ when
$\bQ{\bs{\bar\imath}}{n}$ is unknown, or in other words, the
information in $\bQ{\bs{\bar\imath}}{n}$ contributing to the dynamics
of $\Q{j}{n+1}$.  The interpretation of Eq. (\ref{eq:cau:H_part}) as
the information flux from $\bQ{\bs{\bar\imath}}{n}$ to $Q_j^{n+1}$
motivates our definition of causality. The information-theoretic
causality from $\bQ{\bs{\bar\imath}}{n}$ to $\Q{j}{n+1}$, denoted by
$T_{\bs{\bar\imath} \rightarrow j }$, is defined as the information
flux from $\bQ{\bs{\bar\imath}}{n}$ to $\Q{j}{n+1}$,
\begin{equation}
\label{eq:cau:T}
T_{\bs{\bar\imath} \rightarrow j } = 
\sum_{k=0}^{M-1}  \sum_{ \bs{\bar\imath}(k) \in \mathcal{C}_k}
(-1)^k H(\Q{j}{n+1} | \bQ{\cancel{\bs{\bar\imath}}(k)}{n} ),
\end{equation}
where $\bs{\bar\imath}(k)$ is equal to $\bs{\bar\imath}$ removing
$k$-components and $\mathcal{C}_k$ is the group of all the
combinations of $\bs{\bar\imath}(k)$. Equation~\eqref{eq:cau:T}
represents how much the past information in $\bQ{\bs{\bar\imath}}{n}$
improves our knowledge of the future state $\Q{j}{n+1}$, which is
consistent with the intuition of causality \corr{\citep{wiener1956}}.
Note that the information flux from $T_{\bs{\bar\imath} \rightarrow j
}$ does not overlap with the information flux from
$T_{\bs{\bar\imath}' \rightarrow j }$
for $\bs{\bar\imath}$ different from $\bs{\bar\imath}'$ even if
$\bs{\bar\imath} \cap \bs{\bar\imath}' \neq \emptyset$. For example,
the information flux from $T_{[1,2]\rightarrow j}$ does not overlap
with $T_{1\rightarrow j}$.
This implies that $T_{\bs{\bar\imath}
  \rightarrow j}$ only accounts for the information flux exclusively
due to the joint effect of all the variables in
$\bQ{\bs{\bar\imath}}{n}$.  Figure \ref{fig:cau:TQ} illustrates the
Venn diagram of entropies and information fluxes for a system with
three variables.  The information flux can be cast in compact form
using the generalized conditional mutual information,
\begin{equation}
    T_{\bs{\bar\imath} \rightarrow j } = I( \Q{j}{n+1} ; 
    \Q{i_1}{n} ; \Q{i_2}{n};...;\Q{i_M}{n} | \bQ{\cancel{\bs{\bar\imath}}}{n} ),
\end{equation}
where $I( \cdot ; 
\cdot ; \cdot ;... | \cdot )$
is the conditional co-information~\citep{yeung1991,bell2003}
recursively defined by
\begin{subequations}
  \label{eq:cau:Tco}
\begin{gather}
I(\Q{j}{n+1} ; 
\Q{i_1}{n} ; ...; \Q{i_{M-1}}{n} ; \Q{i_M}{n} | \bQ{\cancel{\bs{\bar\imath}}}{n} ) 
= \\
I(\Q{j}{n+1} ; 
\Q{i_1}{n} ; ...; \Q{i_{M-1}}{n} |  \bQ{\cancel{\bs{\bar\imath}}}{n} ) 
-
I(\Q{j}{n+1} ; 
\Q{i_1}{n} ; ...;\Q{i_{M-1}}{n} | [ \Q{i_M}{n} \bQ{\cancel{\bs{\bar\imath}}}{n}] ).
\end{gather}
\end{subequations}
The recursion in Eq. (\ref{eq:cau:Tco}) is repeated until obtaining
the pairwise definition of conditional mutual information $I(X ; Y |
Z) = H(X|Z) - H(X| [Y,Z])$.

By construction of Eq.~\eqref{eq:cau:T}, it is satisfied that the
amount of information in the state $\Q{j}{n+1}$ is equal to the sum of
all the information fluxes from $\bQ{\bs{\bar\imath}}{n}$ \corr{and
$Q^{n}_{\cancel{\bs{\bar{\imath}}}}$} to $\Q{j}{n+1}$,
\begin{equation}
  \label{eq:cau:P1_pre}
  H(\Q{j}{n+1}) =  \sum_{ \bs{\bar\imath}' \in \mathcal{C}} T_{\bs{\bar\imath}' \rightarrow j },
\end{equation}
where $\mathcal{C}$ is the group of all combinations of vectors
\corr{$\bs{\bar\imath}'$ of length 1 to $N$ with components taken from
  $\bs{\bar\imath} \cup \cancel{\bs{\bar\imath}}$}.
%
Another important property of the information flux is that
$T_{\boldsymbol{\bar\imath} \rightarrow j }=0$ if the dynamics of
$\Q{j}{n+1} $ does not depend explicitly on the states
$\bQ{\boldsymbol{\bar\imath}}{n}$, namely,
\begin{equation}
  \label{eq:cau:P2}
   \Q{j}{n+1} =
f_j(\bQ{\cancel{\bs{\bar\imath}}}{n}) \Rightarrow T_{\boldsymbol{\bar\imath} \rightarrow j } = 0.
\end{equation}
The zero-information-flux condition above is again consistent with the
intuition that no direct causality should emerge from
$\bQ{\bs{\bar\imath}}{n}$ to $\Q{j}{n+1}$ unless the latter depends on
the former.
\corr{Additionally, the information flux is based on probability
  distributions and, as such, is invariant under shifting, rescaling
  and, in general, nonlinear $C^1$-diffeomorphism transformations of
  the signals~\citep{Kaiser2002}. One more attractive feature of the
  information flux is that $T_{\boldsymbol{\bar\imath} \rightarrow j
  }$ accounts for direct causality excluding intermediate
  variables. For example, if the causality flow is $\Q{i}{}
  \rightarrow \Q{j}{} \rightarrow \Q{k}{}$, then there is no causality
  from $\Q{i}{}$ to $\Q{k}{}$ (i.e., $T_{i \rightarrow k} = 0$)
  provided that the three components $\Q{i}{}$, $\Q{j}{}$, and
  $\Q{k}{}$ are contained in $\bQ{}{}$.}

The definition from Eq. (\ref{eq:cau:T}) is trivially generalized to
quantify the information flux from $\bQ{\bs{\bar\imath}}{n}$ to a set
of variables with indices
$\bs{\bar\jmath}=[\bar\jmath_1,\bar\jmath_2,...]$ by replacing
$\Q{j}{n+1}$ by $\bQ{\bs{\bar\jmath}}{n+1}$.
It can be shown that given two arbitrary sets of variables with index
vectors $\bs{\bar\imath}$ and $\bs{\bar\jmath}$, in general, it is
satisfied that $T_{\bs{\bar\imath} \rightarrow \bs{\bar\jmath}} \neq
T_{\bs{\bar\jmath} \rightarrow \bs{\bar\imath}}$, and thus the
information flux is asymmetric.

\subsubsection{Information flux of observable states}

In many occasions, we are interested in, or only have accessed to, an
observable state
\begin{equation}
   \label{eq:cau:Y}
   \bY{}{n} = \bs{h}(\bQ{}{n}),
 \end{equation}
where $\bY{}{n} = [\Y{1}{n},...,\Y{N_Y}{n}]$ with $N_Y \leq N$.
Equation (\ref{eq:cau:Y}) generally entails a loss of information
 \begin{equation}
   \label{eq:cau:Y_loss}
   H(\bY{}{n}) = H(\bs{h}(\bQ{}{n})) \leq H(\bQ{}{n}),
 \end{equation}
such that complete knowledge of $\bY{}{n}$ does not necessarily imply
that the future state $\bY{}{n+1}$ is known. This is revealed by the
inequality $H(\bY{}{n+1} | \bY{}{n}) \geq H( \bY{}{n+1} | \bQ{}{n}) =
0$, that particularized for the $j$-component of $\bY{}{n}$ results in
\begin{equation}
  \label{eq:cau:HY}    
  H(\Y{j}{n+1} | \bY{}{n}) \geq 0.
\end{equation} 

In light of Eq. (\ref{eq:cau:HY}), the definition of information flux
from Eq. (\ref{eq:cau:T}) should be modified to account for the lack
of knowledge from unobserved states. The information flux from an
observable state $\bY{\bs{\bar\imath}}{n}$ to a future observable
state $\Y{j}{n+1}$ is defined as
\begin{equation}
  \label{eq:cau:TY_pre}
  T^Y_{\bs{\bar\imath} \rightarrow j } = 
\left[ \sum_{k=0}^{M-1} \sum_{ \bs{\bar\imath}(k) \in \mathcal{P}_k}
(-1)^{k} H(\Y{j}{n+1} | \bY{\cancel{\bs{\bar\imath}}(k)}{n} )\right] + (-1)^{M} H(\Y{j}{n+1} | \bY{}{n}),
\end{equation}
where ${\bs{\bar\imath}}$ is again a vector of indices,
$\bs{\bar\imath}=[\bar\imath_1,...,\bar\imath_M]$,
now with $M\leq N_Y$.
The term $H(\Y{j}{n+1} | \bY{}{n})$ in Eq. (\ref{eq:cau:TY_pre})
quantifies the information loss from unobserved states and is
naturally absorbed into the summation as
\begin{equation}
  \label{eq:cau:TY}
  T^Y_{\bs{\bar\imath} \rightarrow j} =  \sum_{k=0}^{M} \sum_{ \bs{\bar\imath}(k) \in \mathcal{P}_k}
(-1)^{k} H(\Y{j}{n+1} | \bY{\cancel{\bs{\bar\imath}}(k)}{n} ).
\end{equation}
The definition above can be written in compact form using again the
conditional co-information as
\begin{equation}
    T^Y_{\boldsymbol{\bar\imath} \rightarrow j } = I( \Y{j}{n+1} ; 
    \Y{i_1}{n} ; \Y{i_2}{n};...;\Y{i_{M}}{n} | \bY{\cancel{\bs{\bar\imath}}}{n} ).
\end{equation}
When $\bY{}{n}=\bQ{}{n}$, then $H(\bY{j}{n+1} | \bY{}{n})=0$ and
Eq. (\ref{eq:cau:T}) is recovered. Moreover, for
$\bs{\bar\imath}=[i]$, the information flux for observable states is
\begin{equation}
  \label{eq:cau:TY_1}
  T^Y_{ i \rightarrow j } = 
 H(\Y{j}{n+1} | \bY{\cancel{i}}{n} ) - H(\Y{j}{n+1} | \bY{}{n}),
\end{equation}
that is the multivariate generalization of the transfer entropy
proposed by \citet{schreiber2000}. \corr{Another difference from
  \citet{schreiber2000} is that the new definition of causality
  accounts for the information flux due to the joint effect of
  variables when $\bs{\bar\imath}$ has more than one component.  We
  will show in the example below that the joint information flux might
  be of the same order as the information flux from individual
  variables.}

We have shown in Eq. (\ref{eq:cau:P1_pre}) that the total information
in a future state is determined by the sum of all information fluxes.
However, limiting the knowledge of the system to a reduced set of
observable variables $\bY{}{n}$ entails a leak of information towards
future state $\Y{j}{n+1}$ such that now
\begin{equation}
  \label{eq:cau:P1Y_pre}
  H(\Y{j}{n+1}) =  \sum_{ \bs{\bar\imath} \in \mathcal{P}} T^Y_{\bs{\bar\imath} \rightarrow j }
  + T^Y_{\mathrm{leak},j},
\end{equation}
where $T^Y_{\mathrm{leak},j}$ is the information leak or amount of
information in $\Y{j}{n+1}$ that cannot be explained by $\bY{j}{n}$,
and it is given by
\begin{equation}
\label{eq:cae:leak}
      T^Y_{\mathrm{leak},j} = H(\Y{j}{n+1}| \bY{}{n}).
\end{equation}

The information flux in Eq.~\eqref{eq:cau:T} and Eq.~\eqref{eq:cau:TY}
has units of bits. It is then natural to introduce the normalized
information flux as
\begin{subequations}
  \begin{gather}
      \mathit{TN}^Y_{\boldsymbol{\bar\imath} \rightarrow j} =  \frac{T^Y_{\boldsymbol{\bar\imath} \rightarrow j }}{H(Y^j_{n+1})}, \\
      \mathit{TN}^Y_{\mathrm{leak}, j} =  \frac{T^Y_{\mathrm{leak},j }}{H(Y^j_{n+1})},
\end{gather}
\end{subequations}
which satisfies
\begin{equation}
  \sum_{
  \boldsymbol{\bar\imath} \in \mathcal{P}}
  \mathit{TN}^Y_{\boldsymbol{\bar\imath} \rightarrow j } + \mathit{TN}^Y_{\mathrm{leak}, j} = 1.
\end{equation}
\corr{A similar normalization was proposed by \citet{materassi2014} in
  the context of delayed mutual information. Other normalizations are
  discussed in \citet{Duan2013}.}

As an example, we elaborate on the information flux formulas for
$Y^{n+1}_2$ in a system with three observable variables
($N_Y=3$)\corr{, as depicted in Fig.~\ref{fig:cau:TY}}.  For
$\boldsymbol{\bar\imath} = [1]$ and $j=2$, we get
\begin{equation}
    T^Y_{1 \rightarrow 2 } = H(\Y{2}{n+1} | \Y{2}{n}, \Y{3}{n}) - 
    H(\Y{2}{n+1} | \Y{1}{n}, \Y{2}{n}, \Y{3}{n}),
\end{equation}
for  $\boldsymbol{\bar\imath} = [1,3]$ and $j=2$, 
%
  \begin{gather}    T^Y_{[1,3] \rightarrow 2 } = H(\Y{2}{n+1} | \Y{2}{n}) - H(\Y{2}{n+1} | \Y{1}{n}, \Y{2}{n})  
    -  H(\Y{2}{n+1} | \Y{2}{n}, \Y{3}{n}) +\nonumber\\ 
    +     H(\Y{2}{n+1} | \Y{1}{n}, \Y{2}{n}, \Y{3}{n}),
  \end{gather}
%
and for $\bs{\bar\imath} = [1,2,3]$ and $j=2$, we have
%
  \begin{gather}
      T^Y_{[1,2,3] \rightarrow 2 } = H(\Y{2}{n+1}) - \nonumber\\ 
  - H(\Y{2}{n+1} | \Y{1}{n}) - H(\Y{2}{n+1} | \Y{2}{n}) - H(\Y{2}{n+1} | \Y{3}{n}) \nonumber
  + 
  \\ + H(\Y{2}{n+1} | \Y{1}{n},\Y{2}{n}) + H(\Y{2}{n+1} | \Y{1}{n},\Y{3}{n}) + H(\Y{2}{n+1} | \Y{2}{n},\Y{3}{n}) - \nonumber\\ 
  -     H(\Y{2}{n+1} | \Y{1}{n}, \Y{2}{n}, \Y{3}{n}).
  \end{gather}
%
In the three examples above, $T^Y_{1 \rightarrow 2 }$, $T^Y_{[1,3]
  \rightarrow 2 }$ and $T^Y_{[1,2,3] \rightarrow 2 }$ contain
non-overlapping information similar to the sketch shown in figure
\ref{fig:cau:TQ}. The total information in $\Y{2}{n+1}$ is given by
\begin{gather}
    H(\Y{2}{n+1}) = T^Y_{1 \rightarrow 2 }  + T^Y_{2 \rightarrow 2 }  + T^Y_{3 \rightarrow 2 } +\nonumber \\
    T^Y_{[1,2] \rightarrow 2 } + T^Y_{[1,3] \rightarrow 2 } + T^Y_{[2,3] \rightarrow 2 } + T^Y_{[1,2,3] \rightarrow 2 } + T^Y_{\mathrm{leak},j},
\end{gather}
where $T^Y_{\mathrm{leak},j} = H(\Y{2}{n+1} | \Y{1}{n}, \Y{2}{n}, \Y{3}{n})$.

We close this section by noting that $T^Y_{\bs{\bar\imath}\rightarrow
  j}$ (similarly for $T_{\bs{\bar\imath}\rightarrow j}$) is not
constrained to be larger or equal to zero \corr{when the number of
  variables considered is odd}. This might not be obvious from
figure~\ref{fig:cau:TY}, as conditional entropies do not obey the
conservation of areas depicted by the Venn diagram shown in the
plot. However, the possibility of negative values of the
co-information is a known property often discussed in the
literature~\citep[see, for instance,][]{bell2003}.  In general,
negative information flux will occur when there is backpropagation of
information from the future to the past, i.e, the knowledge of an
event in the future would provide information about an event in the
past, but not the other way around. The reader is referred to
\citet{james2016} for a deeper discussion on the topic.



\subsubsection{Optimal observable states and phase-space partition for information flux}

The analysis of information fluxes is considerably simplified when the
mutual information between pairs of components in $\bY{}{n}$ is zero,
\begin{equation}
I(\Y{i}{n} ; \Y{j}{n}) = 0, \ i\neq j.
\end{equation}  
In that case, we can focus on the information flux from one single
variable to the future state \corr{as shown in
  Eq. \eqref{eq:cau:TY_1}.}
%
%
Given an observable state $\bY{}{n}$, we define the optimal observable
representation for information flux as $\bY{}{n*} =
\bs{w}^*(\bY{}{n})$, where $\bs{w}^*$ is the reversible transformation
satisfying
\begin{equation}
  \boldsymbol{w}^* = \arg\min_{\bs{w}(\bY{}{n})} \left( \sum_{i,j, \ i\neq j} I(\Y{i}{n} ; \Y{j}{n}) \right), \
  \mathrm{s.t.} \ H(\bY{}{n*})=H(\bY{}{n}).
\end{equation}  
The new observable state $\bY{}{n*}$ has the advantage of minimizing
the causal links due to the joint effect of two or more variables
acting together, which might ease the identification of key physical
processes and facilitate the reduced-order modeling of the system.  A
similar argument can be applied to the phase-space partition $D
=\{D_1, D_2,..., D_{N_q}\}$ to define the optimal $D^*$ for causal
inference as
\begin{equation}
  D^* = \arg\min_{D} \left( \sum_{i,j, \ i\neq j} I(\Q{i}{n} ; \Q{j}{n}) \right) \
  \mathrm{s.t.} \ D_i \cap D_j = \emptyset \ \forall i \neq j.
\end{equation}
%

\subsection{Application: Causality of the energy cascade in isotropic turbulence}
\label{sec:cau:example}

\begin{figure}
  \begin{center}
   \includegraphics[width=0.7\textwidth]{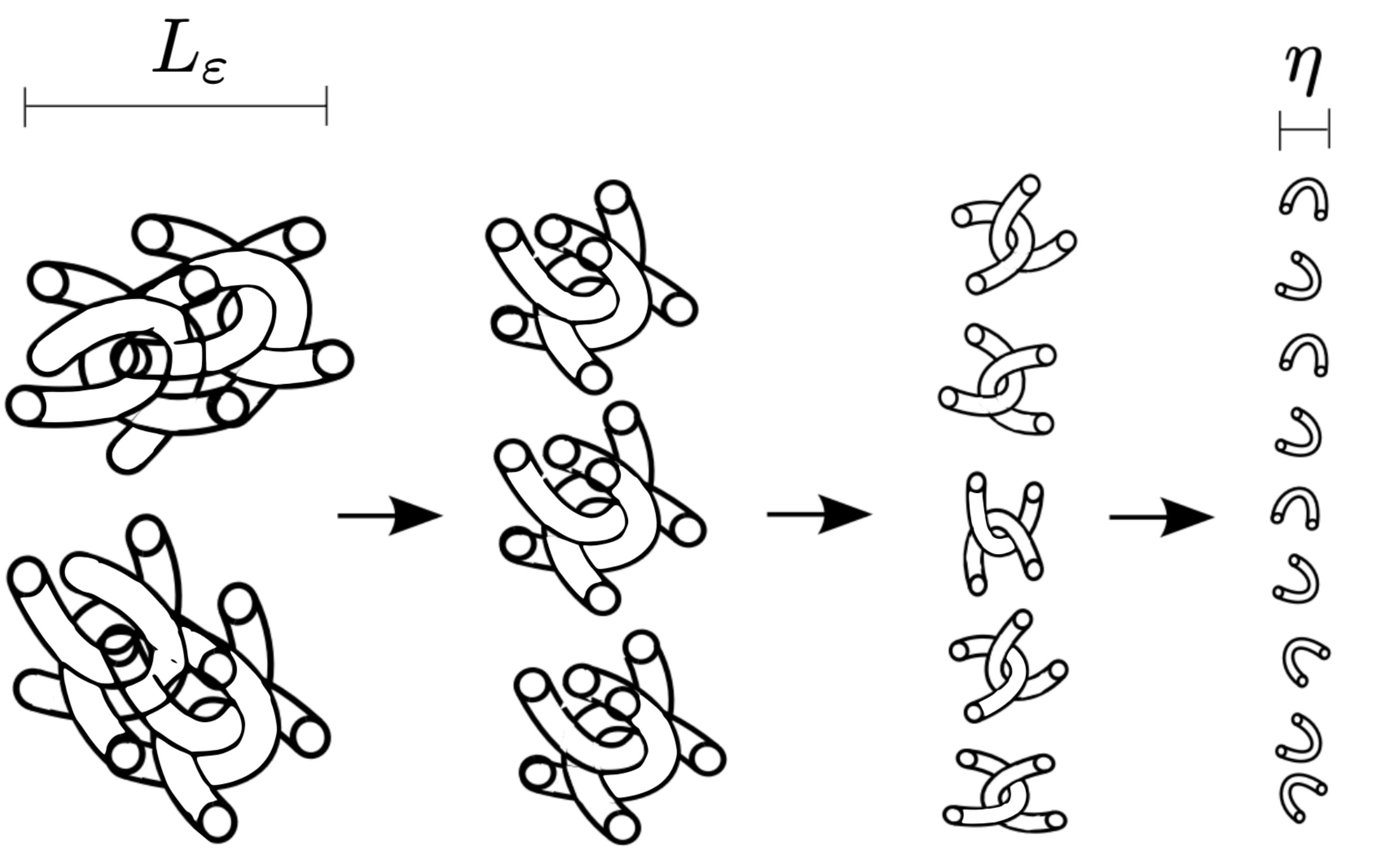}
 \end{center}
\caption{ \corr{Schematic of the Richardson's turbulent energy
    cascade~\citep{richardson1922} in which energy is transferred
    sequentially between eddies of decreasing size. The kinetic energy
    flows from the largest flow motions, characterized by the integral
    length-scale $L_\varepsilon$, to the Kolmogorov length-scale
    $\eta$, where it is finally dissipated.}
  \label{fig:cau:cascade_sketch}}
\end{figure}

The cascade of energy in turbulent flows, i.e., the transfer of
kinetic energy from large to small flow scales or vice versa (backward
cascade), is the cornerstone of most theories and models of turbulence
since the 1940s~\citep{richardson1922, obukhov1941, kolmogorov1941,
  kolmogorov1962, aoyama2005, falkovich2009, cardesa2017}. Yet,
understanding the dynamics of the kinetic energy transfer across
scales remains an outstanding challenge in fluid mechanics. Given the
ubiquity of turbulence, a deeper understanding of the energy transfer
among the flow scales would enable significant progress across various
fields ranging from combustion~\citep{veynante2002},
meteorology~\citep{bodenschatz2015}, and
astrophysics~\citep{young2017} to engineering applications of
aero/hydro-dynamics~\citep{sirovich1997, hof2010, marusic2010,
  kuhnen2018, ballouz2018}. In spite of the substantial advances in
the last decades, the causal interactions of energy among scales in
the turbulent cascade remain uncharted.  Here, we use the formalism
introduced in \S \ref{sec:cau:info_flux} to investigate the
information flux of the turbulent kinetic energy across different
scales. \corr{Our goal is to assess the local-in-scale cascade in
  which the kinetic energy is transferred sequentially from one scale
  to the next smaller scale as sketched in figure
  \ref{fig:cau:cascade_sketch}.}
%

The case selected to study the energy cascade is isotropic turbulence
in a triply periodic box with side $L$.  The data were obtained from
the DNS of \citet{cardesa2015}, which is publicly available in
\citet{torroja2021}. The conservation of mass and momentum equations
of an incompressible fluid are given by
\begin{eqnarray}\label{eq:cau:NS}
\frac{\partial  u_i}{\partial t} + \frac{\partial  u_i  u_j  }{\partial x_j}
 = - \frac{1}{\rho}\frac{\partial  \Pi}{\partial x_i} +
\nu \frac{\partial^2  u_i}{\partial x_j\partial x_j} + f_i,
\quad\frac{\partial  u_i}{\partial x_i}  = 0,
\end{eqnarray}
where repeated indices imply summation, $\boldsymbol{x}=[x_1, x_2,
  x_3]$ are the spatial coordinates, $u_i$ for $i=1,2,3$ are the
velocities components, $\Pi$ is the pressure, $\rho$ is the flow
density, $\nu$ is the kinematic viscosity, and $f_i$ is a linear
forcing sustaining the turbulent flow~\citep{rosales2005}.  The flow
setup is characterized by one nondimensional parameter, the Reynolds
number, \corr{which quantifies the separation of length-scales in the
  flow}. The Reynolds number based on the Taylor
microscale~\citep{pope2000} is $Re_\lambda\approx 380$. The simulation
was conducted by solving Eq. (\ref{eq:cau:NS}) with $1024^3$ spatial
Fourier modes, which is enough to accurately resolve all the relevant
length-scales of the flow.  The system is multiscale and highly
chaotic, with roughly $10^9$ degrees of freedom.

In the following, we summarize the main parameters of the
simulation. The reader is referred to \citet{cardesa2015} for more
details about the flow set-up.  The spatial- and time- average of the
turbulent kinetic energy ($K=u_iu_i/2$) and dissipation
($\varepsilon=2\nu S_{ij}S_{ij}$) are denoted by $K_\mathrm{avg}$ and
$\varepsilon_\mathrm{avg}$, respectively, where $S_{ij} = (\partial
u_i/\partial x_j + \partial u_j/\partial x_i)/2$ is the rate-of-strain
tensor. The ratio between the largest and smallest length-scales of
the problem can be quantified by $L_\varepsilon/\eta = 1800$, where
$L_\varepsilon = K_\mathrm{avg}^{3/2}/ \varepsilon_\mathrm{avg}$ is
the integral length-scale, and $\eta =
(\nu^3/\varepsilon_\mathrm{avg})^{1/4}$ is the Kolmogorov
length-scale.
The data generated is also time-resolved, with flow fields stored
every $\Delta t = 0.0076 T_\varepsilon$, where $T_\varepsilon =
K_\mathrm{avg}/\varepsilon_\mathrm{avg}$, and was purposely run for
long times to enable the reliable computation of conditional
entropies.  The total time simulated after transients was equal to
$165 T_\varepsilon$.

The next step is to quantify the kinetic energy carried by the
energy-containing eddies at various length-scales as a function of
time. To that end, the $i$-th component of the instantaneous flow
velocity $u_i(\boldsymbol{x},t)$ is decomposed into large- and small-
components according to $u_i(\boldsymbol{x},t) =
\bar{u}_i(\boldsymbol{x},t) + u'_i(\boldsymbol{x},t)$, where
$\bar{(\cdot)}$ denotes the low-pass Gaussian filter operator,
\begin{equation}
    \bar{u}_i(\boldsymbol{x},t) = 
    \int_V \frac{\sqrt{\pi}}{\bar{\Delta}} \exp\left[-\pi^2( \boldsymbol{x} -
      \boldsymbol{x}')^2/\bar{\Delta}^2\right] u_i(\boldsymbol{x}') \mathrm{d}\boldsymbol{x}',
\end{equation}
and $\bar{\Delta}$ is the filter width.  The kinetic energy of the
large-scale field evolves as
\begin{equation}
    \left( \frac{\partial}{\partial t} + \bar{u}_j \frac{\partial }{\partial x_j}\right)\frac{1}{2}\bar{u}_i\bar{u}_i 
    =-\frac{\partial }{\partial x_j}\left( \bar{u}_j \bar{\Pi} + \bar{u}_i\tau_{ij}^{\mathrm{SGS}} -  2\nu \bar{u}_i \bar{S}_{ij} \right)
    + \Sigma - 2\nu \bar{S}_{ij}\bar{S}_{ij} + \bar{u}_i \bar{f}_i, 
\end{equation}
where $\tau_{ij}^{\mathrm{SGS}} = (\overline{u_i u_j} - \bar{u}_i
\bar{u}_j)$ is the subgrid-scale stress tensor, \corr{which represents
  the effect of the (filtered) small-scale eddies on the (resolved)
  large-scale eddies}. \corr{The interscale energy transfer between
  the filtered and unfiltered scales is given by $\Sigma =
  \tau_{ij}^{\mathrm{SGS}} \bar{S}_{ij}$, which is the present
  quantity of interest.}

\begin{figure}
  \begin{center}
   \includegraphics[width=0.9\textwidth]{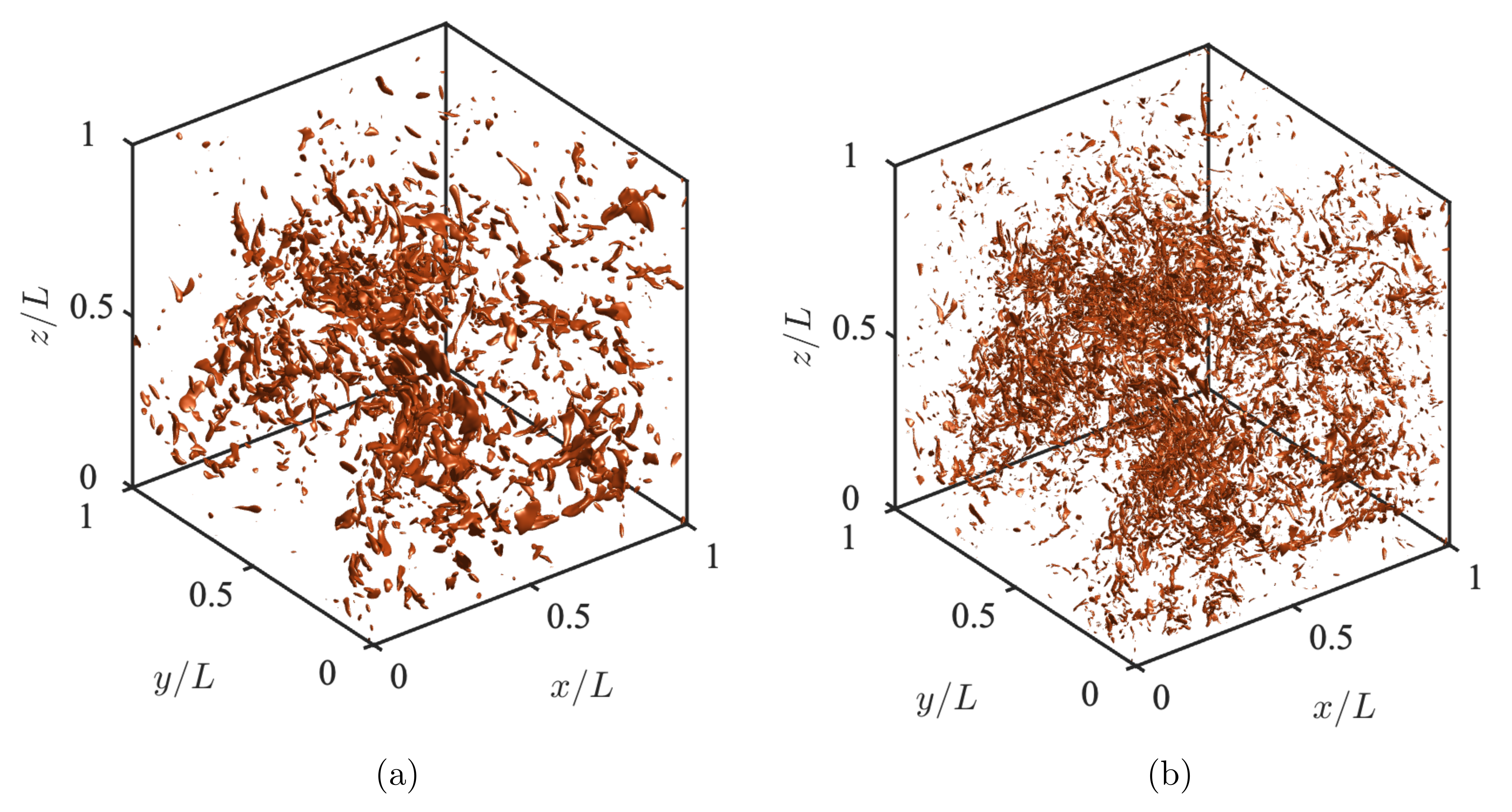}
 \end{center}
\caption{ Isosurfaces of the instantaneous kinetic energy transfer
  $\Sigma_i$ for filter sizes \corr{(a) $\bar\Delta = \bar\Delta_1=
    163 \eta$ (denoted by $\Sigma_1$) and (b) $\bar\Delta =
    \bar\Delta_4 =21 \eta$ (denoted by $\Sigma_4$)} at the same time.
  \label{fig:cau:TKE_iso}}
\end{figure}

\begin{figure}
  \begin{center}
   \includegraphics[width=0.8\textwidth]{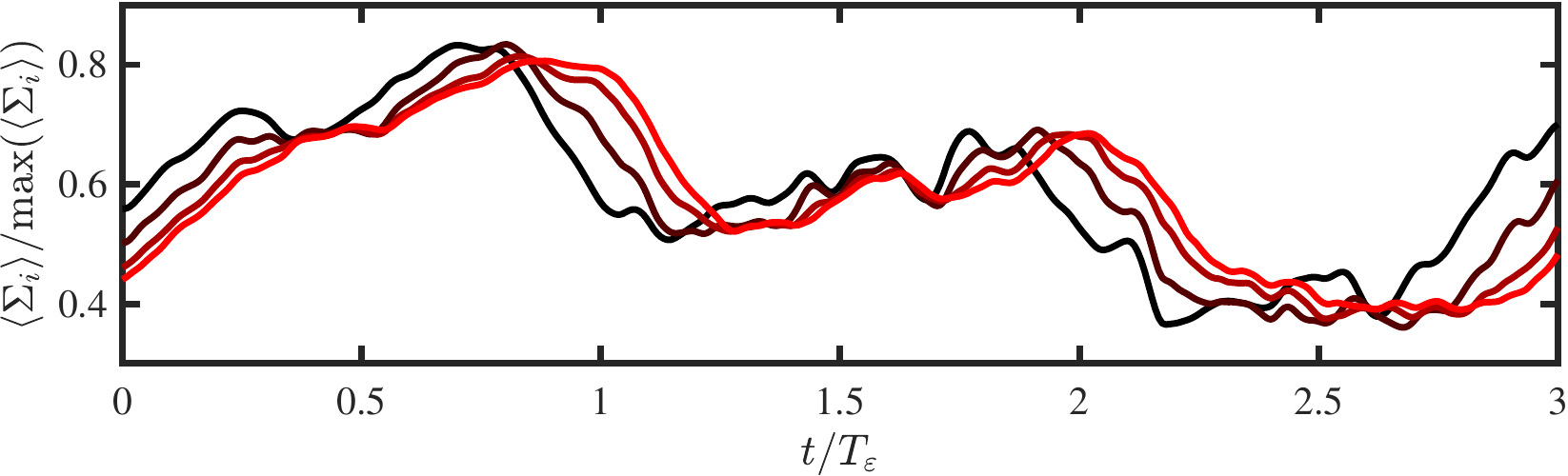}
 \end{center}
\caption{An extract of the time-history of $\langle \Sigma_1 \rangle$, $\langle
  \Sigma_2 \rangle$, $\langle \Sigma_3 \rangle$, and $\langle \Sigma_4
  \rangle$ (from black to red). Although not shown, the whole time-span of the signals is
  $165T_\varepsilon$.
  \label{fig:cau:signals}}
\end{figure}

The velocity field is low-pass filtered at four filter widths:
$\bar{\Delta}_1=163 \eta$, $\bar{\Delta}_2=81\eta$,
$\bar{\Delta}_3=42\eta$, and $\bar{\Delta}_4=21\eta$.  \corr{The
  filter widths selected lay in the inertial range of the simulation
  within the integral and Kolmogorov length-scales: $L_\varepsilon >
  \bar{\Delta}_i > \eta$, for $i=1,2,3$ and 4}. The resulting velocity
fields are used to compute the interscale energy transfer at scale
$\bar{\Delta}_i$, which is denoted by $\Sigma_i(\bs{x},t)$.  Examples
of three-dimensional isosurfaces of $\Sigma_1$ and $\Sigma_4$ are
featured in figure \ref{fig:cau:TKE_iso} to provide a visual reference
of the spatial organization of the interscale energy transfer.  We use
the volume-averaged value of $\Sigma_i$ over the whole domain, denoted
by $\langle \Sigma_i \rangle$, as a marker for the time-evolution of
the interscale energy transfer. Note that $\langle \Sigma_i \rangle$
is only a function of time. Figure \ref{fig:cau:signals} contains a
fragment of the time-history of $\langle \Sigma_i \rangle$ for
$i=1,2,3$ and $4$.
%
%

We can now relate the current formulation with the notation introduced
in \S \ref{sec:cau:info_flux}.  The full system state $\bQ{}{n}$ is
given by the velocity components $u_i$ and pressure $\Pi$ for the
$1024^3$ Fourier modes. The map $\bQ{}{n+1}=\bs{f}(\bQ{}{n})$ is
obtained from the spatio-temporal discretization of the Navier--Stokes
equations in Eq. (\ref{eq:cau:NS}). The observable states $\bY{}{n}$
are represented by the interscale turbulent kinetic energy transfer
$\bY{}{n}=[\langle \Sigma_1\rangle,\langle \Sigma_2\rangle,\langle
  \Sigma_3\rangle,\langle \Sigma_4\rangle]$, and the mapping for the
observable states $\bY{}{n} = \bs{h}(\bQ{}{n})$ is derived from the
definition of $\Sigma_i$ in conjunction with the discrete version of
Eq. (\ref{eq:cau:NS}).

\begin{figure}
  \begin{center}
    \subfloat[]{\includegraphics[width=0.43\textwidth]{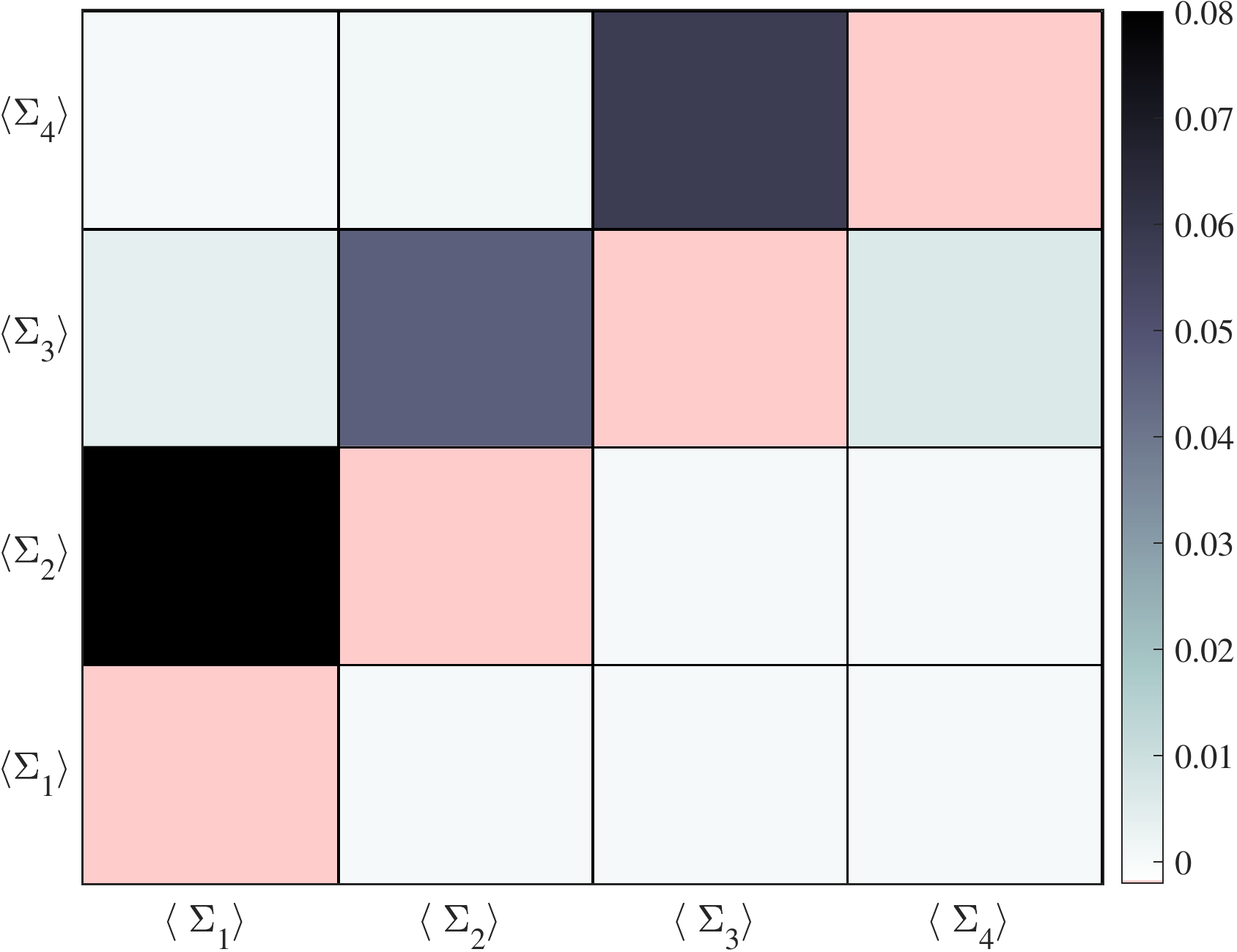}}
    \hspace{0.5cm}
   \subfloat[]{\includegraphics[width=0.43\textwidth]{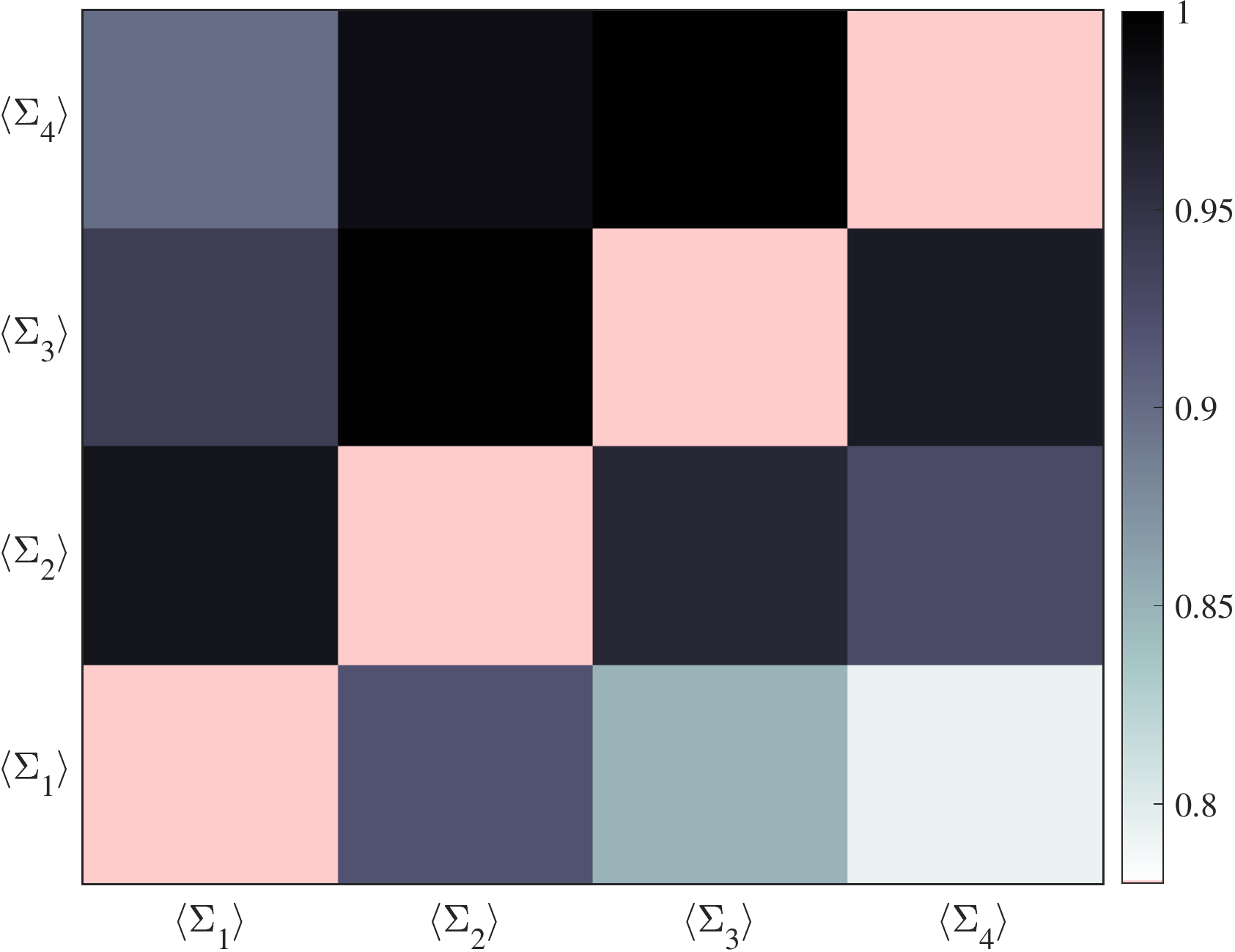}}
    \end{center}
\caption{ (a) Information flux $T^\Sigma_{i \rightarrow j}$ among
  interscale energy-transfer signals at different scales. (b)
  Correlation map $C_{i \rightarrow j}$ between interscale
  energy-transfer signals as defined by Eq. (\ref{eq:cau:Cij}). For
  simplicity, the labels in the axes are $\Sigma_i$, although they
  actually signify $\langle \Sigma_i \rangle$. The self-induced
  intrascale information fluxes $T^\Sigma_{i \rightarrow i}$ in are
  masked in light red.\label{fig:cau:maps1}}
\end{figure}

We examine the propagation of information among $\langle
\Sigma_i\rangle$ by evaluating the information flux defined in
Eq. (\ref{eq:cau:TY}). We focus first on the information flux from one
single energy transfer $\langle \Sigma_i\rangle$ at time $t$ to the
energy transfer $\langle \Sigma_j\rangle$ at time $t+\Delta t$, termed
as $T^\Sigma_{i\rightarrow j}$. The time-delay selected is $\Delta t =
0.046 T_\varepsilon$, which is consistent with the time-lag for energy
transfer reported in the literature~\citep{cardesa2015}. It was tested
that the conclusions drawn below are not affected when the value of
$\Delta t$ was halved and doubled. It was also assessed that $\sum
T^\Sigma_{\bs{\bar\imath}\rightarrow j} + T^\Sigma_{\mathrm{leak}, j}$
is equal to $H(\langle \Sigma_j \rangle)$ to within machine
precision. The information fluxes $T^\Sigma_{i\rightarrow j}$ are
normalized by $H(\langle \Sigma_j\rangle)$ and organized into the
causality map shown in figure \ref{fig:cau:maps1}(a). Our principal
interest is in the interscale propagation of information
($T^\Sigma_{i\rightarrow j}$ with $i\neq j$). Consequently, the
self-induced intrascale information fluxes ($T^\Sigma_{i\rightarrow
  i}$) in figure \ref{fig:cau:maps1}(a) are masked in light red, as
they tend to dominate \corr{(i.e., variables are mostly causal to
themselves)}.  The information fluxes in figure \ref{fig:cau:maps1}(a)
  vividly capture the forward energy cascade of information toward
  smaller scales, which is inferred from strongest information fluxes:
\begin{equation}
  \label{eq:cau:arrow_T}
  T^\Sigma_{1\rightarrow 2} \rightarrow T^\Sigma_{2\rightarrow 3}  \rightarrow  T^\Sigma_{3\rightarrow  4}.
\end{equation}
Backward transfer of information from smaller to larger scales is also
possible, but considerably feeble compared to the forward information
flux. Hence, the present analysis provides the first evidence of the
forward, sequential-in-scale turbulent energy cascade from the
information-theoretic viewpoint \corr{for the full Navier--Stokes
  equations}.  Our results are consistent with previous studies on the
forward the energy cascade using correlation-based
methods~\citep{cardesa2015,cardesa2017} \corr{and
  information-theoretic tools applied to the Gledzer–Ohkitana–Yamada
  shell model~\citep{materassi2014}}.
%

It is also revealing to compare the results in figure
\ref{fig:cau:maps1}(a) with an equivalent time-cross correlation, as
the latter is routinely employed for causal inference by the fluid
mechanics community. The time-cross-correlation `causality' from
$\langle \Sigma_i\rangle$ to $\langle \Sigma_j\rangle$ is defined as
\begin{equation}
\label{eq:cau:Cij}
C_{i\rightarrow j} = \frac{\sum_{n=1}^{N_t} \langle \Sigma_i\rangle(t_n) \langle \Sigma_j\rangle(t_n+\Delta t) }
{\left( \sum_{n=1}^{N_t}\langle \Sigma_i\rangle^2(t_n) \right)^{1/2}
\left( \sum_{n=1}^{N_t}\langle \Sigma_j\rangle^2(t_n) \right)^{1/2}},
\end{equation}
where $\langle \Sigma_i\rangle(t_n)$ signifies $\langle
\Sigma_i\rangle$ at time $t_n$, and $N_t$ is the total number of times
stored of the simulation.  The values of Eq. (\ref{eq:cau:Cij}) are
bounded between 0 and 1. The correlation map, $C_{i\rightarrow j}$, is
shown in figure \ref{fig:cau:maps1}(b). The process portrayed by
$C_{i\rightarrow j}$ is far more intertwined than its information-flux
counterpart offered in figure \ref{fig:cau:maps1}(a). Similarly to
$T^\Sigma_{i\rightarrow j}$, the correlation map also reveals the
prevailing nature of the forward energy cascade ($C_{i\rightarrow j}$
larger for $j>i$). However, $C_{i\rightarrow j}$ is always above 0.8,
implying that all the interscale energy transfers are tightly coupled.
This is inconsistent with the information flux in figure
\ref{fig:cau:maps1}(a) and is probably due to the inability of
$C_{i\rightarrow j}$ to compensate for the effect of intermediate
variables (e.g., a cascading process of the form $\Sigma_1 \rightarrow
\Sigma_2 \rightarrow \Sigma_3$ would result in non-zero correlation
between $\Sigma_1$ and $\Sigma_3$ via the intermediate variable
$\Sigma_2$). As a consequence, $C_{i\rightarrow j}$ also fails at
shedding light on whether the energy is cascading sequentially from
the large scales to the small scales (i.e.  $\langle \Sigma_1 \rangle
\rightarrow \langle \Sigma_2 \rangle \rightarrow \langle \Sigma_3
\rangle\rightarrow \Sigma_4$), or on the other hand, the energy is
transferred between non-contiguous scales (e.g.  $\langle \Sigma_1
\rangle\rightarrow \langle \Sigma_3 \rangle$ without passing through
$\langle \Sigma_2 \rangle$). We have seen that the information flux in
figure \ref{fig:cau:maps1}(a) supports the former: the energy is
predominantly transfer sequentially according to relation in Eq.
(\ref{eq:cau:arrow_T}). Overall, the inference of causality based on
the time-cross correlation is obscured by the often mild asymmetries
in $C_{i\rightarrow j}$ and the failure of $C_{i\rightarrow j}$ to
account for the effects of a third variable.  In contrast, the causal
map in figure \ref{fig:cau:maps1}(a) conveys a more intelligible
picture of the influence among energy transfers at different scales.

\begin{figure}
  \vspace{0.5cm}
   \begin{center}
    \subfloat[]{\includegraphics[width=0.47\textwidth]{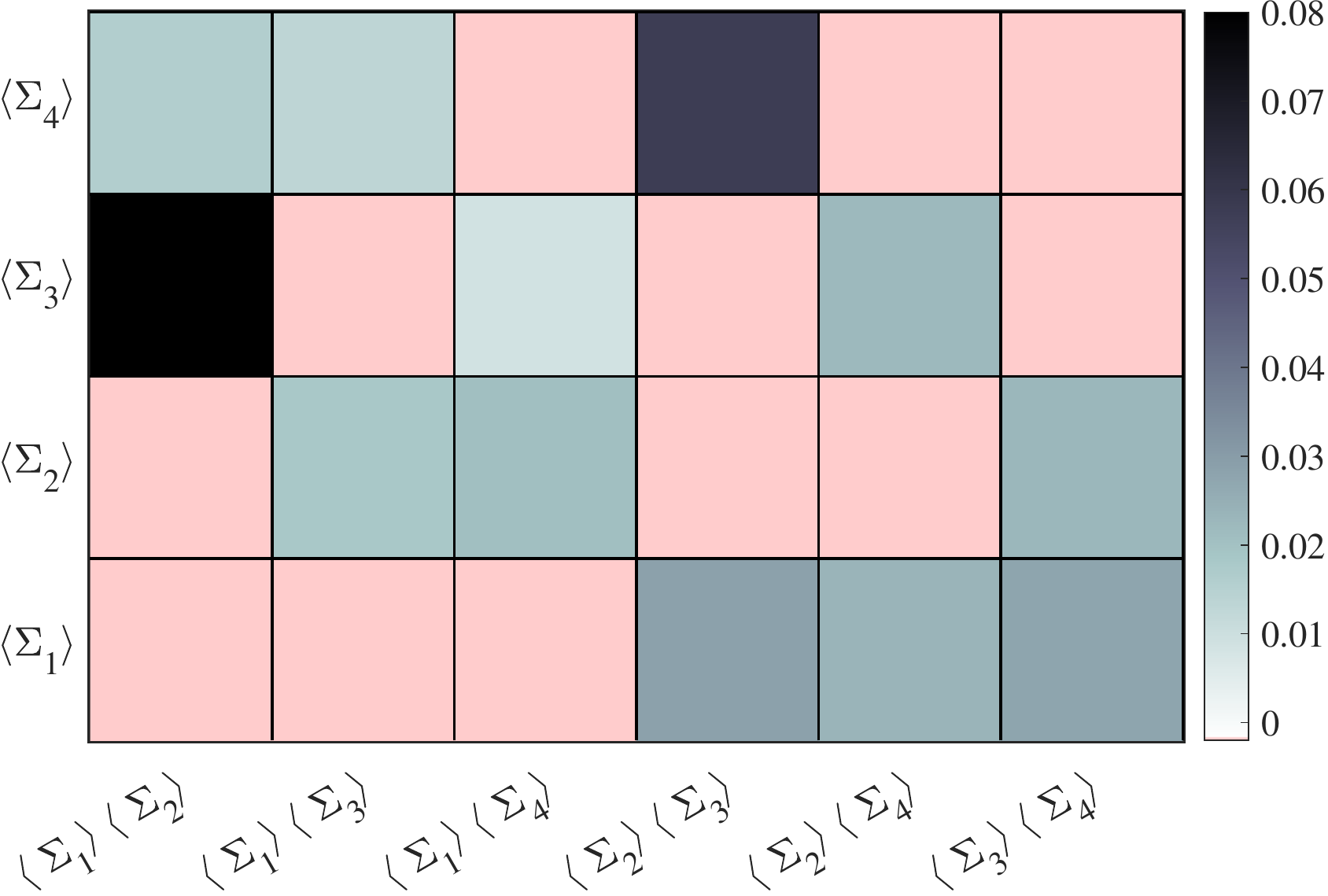}}
    \hspace{0.5cm}
    \subfloat[]{\includegraphics[width=0.38\textwidth]{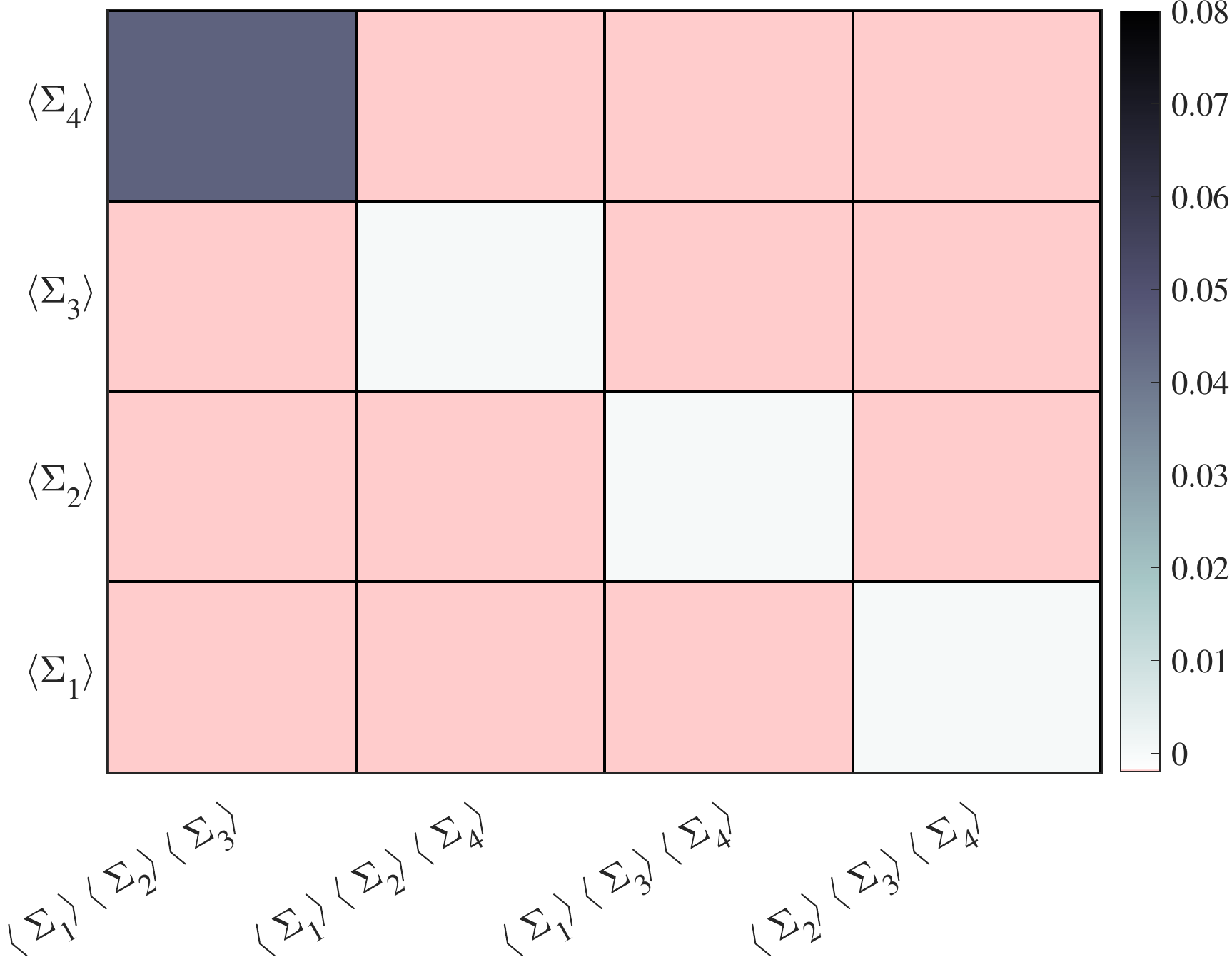}}
 \end{center}
\caption{ Information flux among energy-transfer signals at different
  scales for (a) $T^\Sigma_{[i,k] \rightarrow j}$ and (b)
  $T^\Sigma_{[i,k,l] \rightarrow j}$. The intrascale information
  fluxes $T^\Sigma_{[i,k] \rightarrow i}$ and $T^\Sigma_{[i,k,l]
    \rightarrow i}$ are masked with light red color. For simplicity,
  the labels in the axes are $\Sigma_i$, although they actually
  signify $\langle \Sigma_i \rangle$.
  \label{fig:cau:maps2}}
\end{figure}

For completeness, figure \ref{fig:cau:maps2} includes the information
flux due to the joint effect of two and three variables, where the
values of $T^\Sigma_{[i,j]\rightarrow j}$ and
$T^\Sigma_{[i,j,k]\rightarrow j}$ have also been masked for
clarity. The largest information fluxes are
\begin{equation}
T^\Sigma_{[1,2]\rightarrow 3}, \ T^\Sigma_{[2,3]\rightarrow 4}, \ \mathrm{and} \ T^\Sigma_{[1,2,3]\rightarrow 4},
\end{equation}
which are found to be of the same order of magnitude as $T^\Sigma_{i
  \rightarrow j}$.  The result is again consistent with the prevailing
downscale propagation of information of the energy cascade.

Finally, we calculate the information leak
($T^\Sigma_{\mathrm{leak},j}$) from Eq. (\ref{eq:cae:leak}) to
quantify the amount of information unaccounted for by the observable
variables. The ratios $T^\Sigma_{\mathrm{leak},j}/H(\langle \Sigma_j
\rangle)$ are found to be $0.35, 0.24, 0.18,$ and $0.13$ for $j=1, 2,
3$ and 4, respectively.  Therefore, the information from unobserved
states diminishes towards the smallest scales. The largest leak occurs
for $\langle \Sigma_1\rangle$, where $\sim$35\% of the information
comes from variables not considered within the set $[\langle
  \Sigma_1\rangle$, $\langle \Sigma_2\rangle$, $\langle
  \Sigma_3\rangle$, $\langle \Sigma_4\rangle]$.
%

\section{Modeling}
\label{sec:modeling}


A crucial step in reduced-order modeling of physical system consists
of the identification of transformations enabling parsimonious, yet
informative representations of the full system.  While in some cases
transformations can be carried out on the basis of intuition and
experience, straightforward discrimination of the most relevant
degrees of freedom is challenging for complex problems, most notably
for chaotic, high-dimensional physical systems.  In this context,
information theory has emerged as a valuable framework for model
selection and multimodel inference. Particularly noteworthy is the
work by Akaike~\citep{akaike1974, akaike1977, akaike1998}, where
models are selected on the basis of the relative amount of information
from observations they are capable of accounting for. The approach,
which shares similarities with Bayesian inference, offers an elegant
generalization of the maximum likelihood criterion via entropy
maximization~\citep{akaike1998}.  Akaike's ideas have also bridged the
use of the Kullback's information for parameter estimation and model
identification techniques~\citep{kullback1951}.  Other relevant
studies have leveraged information-theoretic tools for coarse-graining
of dynamical systems assisted by Monte-Carlo methods, renormalization
group, or mapping-entropy techniques~\citep[e.g.][]{baram1978,
  lenggenhager2020, giulini2020}.  A detailed survey of
information-theoretic model selection and inference can be found in
\citet{burnham2002} and \citet{anderson2008}.  In the last decade,
information theory has also become instrumental in machine learning,
mostly within the subfield of deep reinforcement
learning~\citep{sutton2018}. Examples of the latter are neural network
training via the cross-entropy cost functional and estimation of
confidence bounds for value functions and agent
policies~\citep[e.g.][]{still2009, ortega2013, russo2016,
  leibfried2017, koch2018, lu2019}.  \corr{In this section, we
  formulate the problem of reduced-order modeling for chaotic,
  high-dimensional dynamical systems within the framework of
  information theory.  We derived the equation that relates model
  accuracy with the amount of information preserved from the original
  system.  The conditions for maximum information-preserving models
  are also formulated in terms of the mutual information and
  Kullback-Leibler divergence of the quantities of interest}. The
theory is applied to devise a subgrid-scale model for large-eddy
simulation of isotropic turbulence.

\subsection{Formulation}
\label{sec:modeling:formulation}


Let us denote the state of the system to be modeled at time $t_n$ by
$\bQ{}{n}=[\Q{1}{n},...,\Q{N}{n}]$, where $N$ is the total number of
degrees of freedom.  The dynamics of the full system are completely
determined by
\begin{equation}
  \label{eq:model:Q_total}
  \bQ{}{n+1} = \bs{f}(\bQ{}{n}),
\end{equation}
where the map $\bs{f}$  advances the state of the system to an
arbitrary time in the future. It was shown in \S
\ref{sec:information_dynamical} that by construction of
Eq. (\ref{eq:model:Q_total}), it holds that
\begin{equation}
  \label{eq:model:H_total}
  H(\bQ{}{n+1} | \bQ{}{n} ) = 0,
\end{equation}
i.e., there is no uncertainty in the future state $\bQ{}{n+1}$ given
the past state $\bQ{}{n}$.

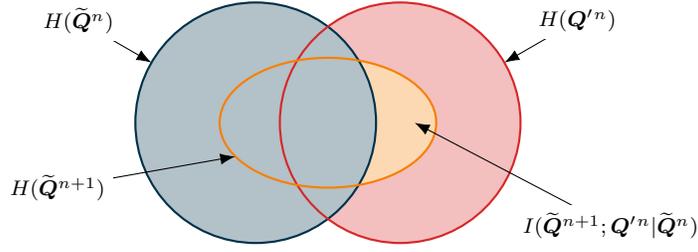
\begin{figure}
    \centering
    \colorlet{c1}{myc1}
\colorlet{c2}{myc2}
\colorlet{c3}{myc3}

\colorlet{c1s}{c1!30}
\colorlet{c2s}{c2!30}
\colorlet{c3s}{c3!30}

\begin{tikzpicture}[scale=1.6,thick,>={Latex[length=.2cm]}
]
    
    \begin{footnotesize}

    \pgfmathsetmacro{\rdi}{1}           
    \pgfmathsetmacro{\rdiO}{.9*\rdi}    
    \pgfmathsetmacro{\rdiOm}{.6*\rdiO}  
    \pgfmathsetmacro{\dis}{.6}          

    \coordinate (O1) at (-\dis,0);  
    \coordinate (O2) at (\dis,0);

    \fill [c2s] (O2) circle (\rdi);                           
    \fill [c3s] ellipse[x radius=\rdiO,y radius=\rdiOm];      

    \draw [c1,fill=c1s] (O1) circle (\rdi); 
    \draw [c2] (O2) circle (\rdi); 
    \draw [c3,name path= ellip] ellipse[x radius=\rdiO,y radius=\rdiOm];
   
    \draw[<-,thin] ($(-\dis,0)+(150:\rdi)$) --++ (150:.7)  
        node[anchor=center,fill=white]{$H(\bQtilde{}{n})$};
    \draw[<-,thin] ($(+\dis,0)+(30:\rdi)$) --++ (30:.7)  
        node[anchor=center,fill=white]{$H(\bQprime{}{n})$};

    \path [name path=line1] (0,0) -- (200:3*\rdi);
    \path [name intersections={of=ellip and line1,by=E}];
    \draw[<-,thin] (E) --++ (190:1.5)  
        node[anchor=center,fill=white]{$H(\bQtilde{}{n+1})$};

    \draw[<-,thin] (0:.8*\rdiO) -- (-20:2.5*\rdi)  
        node[anchor=center,fill=white]
        {$I(\bQtilde{}{n+1};\bQprime{}{n}|\bQtilde{}{n})$};

    \pgfmathsetmacro{\textd}{2.9}

%
%


    \end{footnotesize}

\end{tikzpicture}
    \caption{Schematic of the relationship among the entropies of the
      modeled (accessible) state $\widetilde{\bs{Q}}^{n}$, the future
      modeled state $\widetilde{\bs{Q}}^{n+1}$, and the inaccessible
      degrees of freedom $\bs{Q'}^{n}$.\label{fig:model:infoloss}}
\end{figure}

We aim at modeling a subset of the phase-space of the full system
denoted by $\bQtilde{}{n}=[\Q{1}{n},...,\Q{\tilde{N}}{n}]$ with
$\tilde{N}<N$, where $\tilde{N}$ are the degrees of freedom of the
model. Accordingly, the state vector of the full system is decomposed
as
\begin{equation}
    \bQ{}{n} = [\bQtilde{}{n} , \bQprime{}{n}],
\end{equation}
where $\bQtilde{}{n}$ is the state to be modeled (e.g., the
information accessible to the model) and $\bQprime{}{n}$ are the
inaccessible degrees of freedom that the model must account for.  The
exact dynamics of the modeled state is governed by
\begin{equation}
  \label{eq:model:model_exact}
  \bQtilde{}{n+1} = \tilde{\bs{f}}(\bQtilde{}{n},\bQprime{}{n}), 
\end{equation}
where $\tilde{\bs{f}}$ are the components of $\bs{f}$ 
corresponding to the states $\bQtilde{}{n}$.  It can be readily shown from
Eq. (\ref{eq:model:model_exact}) that disposing of $\bQprime{}{n}$ may
result in an increase of uncertainty in the future states quantified
by
\begin{equation}
  \label{eq:model:info_loss}
  H(\bQtilde{}{n+1} | \bQtilde{}{n} )
  = I(\bQtilde{}{n+1} ; \bQprime{}{n} | \bQtilde{}{n} ) 
 \geq  H(\bQtilde{}{n+1} | \bQtilde{}{n}, \bQprime{}{n}) = 0.
\end{equation}
Equation (\ref{eq:model:info_loss}) represents the fundamental loss of
information for truncated systems: given the initial truncated state
$\bQtilde{}{n}$, the uncertainty in the future state $\bQtilde{}{n+1}$
is equal to the information shared between $\bQprime{}{n}$ and
$\bQtilde{}{n+1}$ that cannot be accounted for by $\bQtilde{}{n}$.  
Figure \ref{fig:model:infoloss} provides a visual representation of 
Eq.~\eqref{eq:model:info_loss}.
If the system is reversible, then Eq. (\ref{eq:model:info_loss}) reduces
to
\begin{equation}
  \label{eq:model:info_loss_rev}
  H(\bQtilde{}{n+1} | \bQtilde{}{n} )
  = H( \bQprime{}{n} | \bQtilde{}{n}), 
\end{equation}
and the uncertainty in the future truncated state $\bQtilde{}{n+1}$ is
equal to the amount of information in $\bQprime{}{n}$ that cannot be
recovered from $\bQtilde{}{n}$.

\subsubsection{Information-theoretic bounds to model error}
\label{subsec:model:bounds}

Let us consider a model with access to the information contained in
$\bQtilde{}{n}$ (i.e., the exact initial condition for the truncated
state) but not to the information in $\bQprime{}{n}$ (i.e.,
inaccessible degrees of freedom).  The governing equation for the
model is denoted by
\begin{equation}
  \label{eq:model:model}
  \bQhat{}{n+1} = \widehat{\bs{f}}(\bQtilde{}{n}),
\end{equation}
where $\bQhat{}{n+1}$ is the model prediction, which does not need to
coincide with the exact solution $\bQtilde{}{n+1}$ obtained from
Eq. (\ref{eq:model:model_exact}) using the exact map
$\widetilde{\bs{f}}$.  We aim at finding a model map
$\widehat{\bs{f}}$ that predicts the future state to within the error
$\varepsilon$,
\begin{equation}
  \label{eq:model:error}
    ||\bQhat{}{n+1} - \bQtilde{}{n+1}|| \leq \varepsilon,
\end{equation}
where $||\cdot||$ is the L$_1$ norm. In particular, we are interested
in the bounds for the error expectation
\begin{equation}
  \label{eq:model:error_expected}
\mathbb{E}[||\bQhat{}{n+1} - \bQtilde{}{n+1}||].
\end{equation}
Given the model prediction from Eq. (\ref{eq:model:model}), the
uncertainty in the exact solution is quantified by
\begin{equation}
  \label{eq:model:info_loss_2}
    H(\bQtilde{}{n+1} | \bQhat{}{n+1} ),
\end{equation}
that we seek to relate to the model error $\varepsilon$.  In general,
nothing can be said about the relationship between $H(\bQtilde{}{n+1}
| \bQhat{}{n+1} )$ and the loss of information of the truncated system
$H(\bQtilde{}{n+1} | \bQtilde{}{n} )$. Thus, the latter might be
smaller, equal, or larger than $H(\bQtilde{}{n+1} | \bQhat{}{n+1} )$
contingent on $\widehat{\bs{f}}$.

Let us denote by $P_e$ the probability of obtaining a modeling error
above the prescribed tolerance $\varepsilon$,
\begin{equation}
P_e = \mathrm{Pr}(||\bQhat{}{n+1} - \bQtilde{}{n+1}|| > \varepsilon).
\end{equation}
The error from Eq. (\ref{eq:model:error}) can be related to the
uncertainty in Eq. (\ref{eq:model:info_loss_2}) via a generalized
Fano's inequality as
\begin{equation}
\label{eq:model:Pe}
P_e \geq \frac{H(\bQtilde{}{n+1} | \bQhat{}{n+1} ) -\log_2(\varepsilon/\Delta_Q) - 1}
{\log_2(\tilde{N})- \log_2(\varepsilon/\Delta_Q)},
\end{equation}
where $\Delta_Q$ is a measure of the size of the phase-space partition
$D_i$ introduced in \S \ref{sec:information_dynamical}. Equation
(\ref{eq:model:Pe}) reveals that, given a model $\widehat{\bs{f}}$,
the probability of incurring an error larger than $\varepsilon$ is
lower bounded by the information loss of the model. It is convenient
to rewrite Eq. (\ref{eq:model:Pe}) as
\begin{equation}
\label{eq:model:Pe_I}
P_e \geq \frac{H(\bQtilde{}{n+1}) -
  I( \bQtilde{}{n+1} ; \bQhat{}{n+1} ) -\log_2(\varepsilon/\Delta_Q) - 1}
{\log_2(\tilde{N}) - \log_2(\varepsilon/\Delta_Q)},
\end{equation}
where $I( \bQtilde{}{n+1} ; \bQhat{}{n+1} )$ is the
mutual information between the `true' state and the model
prediction. A lower bound for the expected error is found by applying
the Markov's inequality to Eq. (\ref{eq:model:Pe_I}),
\begin{equation}
  \label{eq:model:E_I}
 \mathbb{E}[|| \bQhat{}{n+1} - \bQtilde{}{n+1}||]
 \geq \varepsilon \frac{H(\bQtilde{}{n+1}) - I(\bQtilde{}{n+1} ;
 \bQhat{}{n+1} ) -\log_2(\varepsilon/\Delta_Q) - 1 }{\log_2(\tilde{N})- \log_2(\varepsilon/\Delta_Q)}.
\end{equation}
Note that $H(\bQtilde{}{n+1})$ in Eq. (\ref{eq:model:E_I}) is just the
information content of the true state, which is unaffected by the
model. Therefore, the potential predictive capabilities of a model are
attained by maximizing the mutual information between
$\bQtilde{}{n+1}$ and $\bQhat{}{n+1}$. Equation (\ref{eq:model:E_I})
echoes the intuition that the performance of a reduced-order model
improves with the amount of information preserved from the system to
be modeled.  Another advantage of formulating the modeling problem in
terms of mutual information $I(\bQtilde{}{n+1} ; \bQhat{}{n+1} )$ is
that the latter is a concave function of its arguments, which
facilitates the optimization of Eq. (\ref{eq:model:E_I}).

A second model condition can be derived by relaxing the error
constraint in Eq. (\ref{eq:model:error}) to
\begin{equation}
  \label{eq:model:error_P}
    ||p(\bqhat{}{n+1}) - p(\bqtilde{}{n+1})|| \leq \varepsilon'.
\end{equation}
where $p(\bqtilde{}{n+1})$ is the true probability distribution of the
system state and $p(\bqhat{}{n+1})$ is the probability distribution of
the model state.  The error constraint in Eq. (\ref{eq:model:error_P})
is weaker than the constraint in Eq. (\ref{eq:model:error}), as
evidenced by the fact that $||\bQhat{}{n+1} - \bQtilde{}{n+1}|| \geq
0$ even if $p(\bqhat{}{n+1}) = p(\bqtilde{}{n+1})$. Hence, a model can
flawlessly replicate the probability distribution (i.e., the
statistics) of the actual state, yet the sequential samples drawn from
the model (i.e., the dynamics) might not coincide with the ground
truth owing to the lack of mutual information between $\bQhat{}{n+1}$
and $\bQtilde{}{n+1}$.

The error defined by Eq. (\ref{eq:model:error_P}) allow us to estimate
an upper bound for the expectation of the modeling error of
probabilities. First, let us introduce the Kullback-Leibler (KL)
divergence between $p(\bqtilde{}{n+1})$ and $p(\bqhat{}{n+1})$,
\begin{equation}\label{eq:model:KL}
  \mathrm{KL}( \bQtilde{}{n+1} , \bQhat{}{n+1}) =  \sum p(\bqtilde{}{n+1})
  \log [ p(\bqtilde{}{n+1})/p(\bqhat{}{n+1}) ],
\end{equation}
which is a measure of the average number of bits required to recover
$p(\bqtilde{}{n+1})$ using the information in $p(\bqhat{}{n+1})$. From
a Bayesian inference viewpoint, $\mathrm{KL}( \bQtilde{}{n+1},
\bQhat{}{n+1})$ represents the information lost when
$p(\bqhat{}{n+1})$ is used to approximate $p(\bqtilde{}{n+1})$.
Equation (\ref{eq:model:KL}) is an extension of Shannon's concept of
information and is sometimes referred to as relative
entropy~\citep{hobson1973, soofi1994}.  It can be shown via the
Pinsker's inequality~\citep{weissman2003} that
\begin{equation}
  \label{eq:model:error_KL_bound}
  \mathrm{KL}( \bQtilde{}{n+1} , \bQhat{}{n+1}) \geq \frac{1}{2 \ln
    2} ||p(\bqhat{}{n+1}) - p(\bqtilde{}{n+1})||^{2},
\end{equation}
with KL$(\bQtilde{}{n+1},\bQhat{}{n+1})=0$ if and only if the model
predictions are statistically identical to those from the original
system. Equation (\ref{eq:model:error_KL_bound}) provides a connection
between information loss and probabilistic model performance.
Desirable maps $\widehat{\bs{f}}$ are those minimizing
Eq. (\ref{eq:model:KL}), which results in models containing the
coherent information in the data, while leaving out the incoherent
noise. Similarly to the mutual information, the KL divergence has the
advantage of being convex with respect to the input arguments, which
facilitates the search of the minimum.

\corr{Equation (\ref{eq:model:KL}) can be written as
\begin{equation}\label{eq:model:KL_expand}
  \mathrm{KL}( \bQtilde{}{n+1} , \bQhat{}{n+1}) =  
   \sum -p(\bqtilde{}{n+1}) \log [ p(\bqhat{}{n+1}) ] - H(\bQtilde{}{n+1}) 
\end{equation}
where the first term in the right-hand side of
Eq. (\ref{eq:model:KL_expand}) is referred to as the cross entropy
between $\bQtilde{}{n+1}$ and $\bQhat{}{n+1}$.  Taking into account
that $H(\bQtilde{}{n+1})$ is fixed and the cross entropy is equal or
larger than zero, minimizing $\mathrm{KL}( \bQtilde{}{n+1} ,
\bQhat{}{n+1})$ also implies minimizing the cross entropy, resulting
in the minimum cross-entropy principle~\citep{jaynes1957,
  abbas2017}. Additionally, if $p(\bqtilde{}{n+1})$ is taken to be the
uniform distribution, then
\begin{equation}\label{eq:model:KL_ML}
  \mathrm{KL}( \bQtilde{}{n+1} , \bQhat{}{n+1}) = \sum
  -\frac{1}{\tilde{N}} \log[ p(\bqhat{}{n+1}) ] - \log(\tilde{N}),
\end{equation}  
and minimizing Eq. (\ref{eq:model:KL_ML}) is equivalent to maximizing
$\sum \log[ p(\bqhat{}{n+1}) ]$. The latter is the well-known maximum
likelihood principle, which surfaces as a particular case of the
KL-divergence minimization proposed here.}

\subsubsection{Conditions for maximum information-preserving models}

The error bounds presented above provide the information-theoretic
foundations for model discovery. The discussion has been centered on
$\bQtilde{}{n+1}$; however, in most occasions, we are not interested
in the prediction of the full truncated state, but rather in some
quantity of interest
\begin{equation}
  \bYtilde{}{n+1} = \bs{h}(\bQtilde{}{n+1}),
\end{equation}
such that the dimensionality of $\bYtilde{}{n+1}$, denoted by $N_Y$, is much
smaller than the dimensionality of $\bQtilde{}{n+1}$, i.e., $N_Y \ll
\tilde{N}$.  One example is the aerodynamic modeling of an airfoil:
the modeled state $\bQtilde{}{n+1}$ is the flow around the airfoil,
which could contain millions of degrees of freedom, whereas
$\bYtilde{}{n+1}$ may be the surface forces, which contain a few
degrees of freedom.  Most of the discussion in \S
\ref{subsec:model:bounds} is applicable to the modeling of
$\bYtilde{}{n+1}$.  Given the quantity of interest predicted by the
model
\begin{equation}
  \bYhat{}{n+1} = \bs{h}(\bQhat{}{n+1}),
\end{equation}
the bounds for the modeling error are
\begin{subequations}
  \label{eq:model:summary}
  \begin{gather}
  \mathbb{E}[||\bYhat{}{n+1} -
   \bYtilde{}{n+1}||]  \geq \varepsilon_Y \frac{H(\bQ{}{n+1}) - I( \bYtilde{}{n+1} ;
    \bYhat{}{n+1} ) -\log_2 (\varepsilon_Y/\Delta_Y) - 1 }{\log_2(N_Y) - \log_2(\varepsilon_Y/\Delta_Y)}, \label{eq:model:summary_1} \\
   ||p(\byhat{}{n+1}) - p(\bytilde{}{n+1})|| \leq  \left( 2 \ln2 \ \mathrm{KL}( \bYtilde{}{n+1} , \bYhat{}{n+1}) \right)^{1/2}, \label{eq:model:summary_2}
\end{gather}
\end{subequations}
where $\Delta_Y$ is the partition size for the state
$\bYtilde{}{n+1}$. In summary, Eq. (\ref{eq:model:summary})
establishes that a faithful model must \emph{i)} maximize the mutual
information between the model state and the true state, and \emph{ii)}
minimize the KL divergence between their probabilities. The mutual
information assist the model to reproduce the dynamics of the original
system, while the KL divergence enables the accurate prediction of the
statistical quantities of interest. Whether we choose to optimize the
mutual information, the KL divergence, or a combination of both
depends on the scope of the model.  We close this section by remarking
that Eq. (\ref{eq:model:summary}) is a necessary condition for the
discovery of accurate models, but it is not
sufficient. Eq. (\ref{eq:model:summary}) does not provide the information of what
physical quantities should be preserved by the model, nor the modeling
assumptions to undertake. Those will rely on physical insight of the
system to model and clear understanding of the relevant characteristic scales 
(length, time, velocities,...)  involved in the problem.

\subsection{Application: Maximum information-preserving subgrid-scale model for LES}

\begin{figure}
  \begin{center}
   \includegraphics[width=0.6\textwidth]{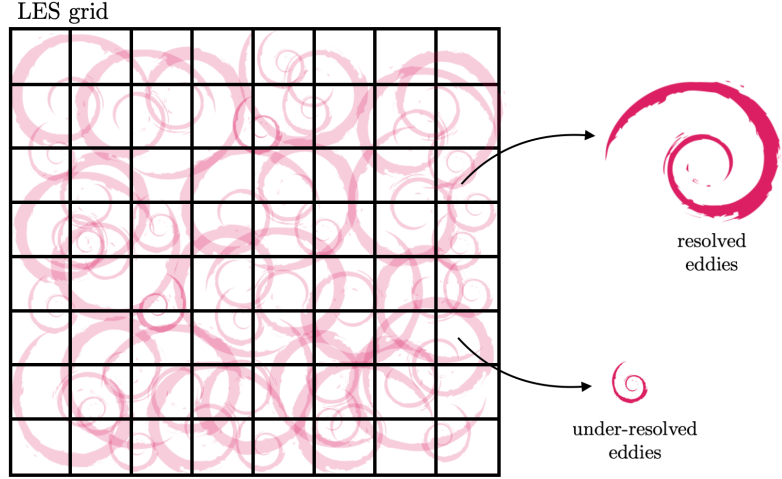}
 \end{center}
\caption{ \corr{Schematic of an LES grid and the turbulent eddies of
    different sizes. Only the large eddies are resolved by the grid,
    whereas the information from the small-scale eddies is lost.}
  \label{fig:model:LES_sketch}}
\end{figure}

Most turbulent flows of engineering significance cannot be simulated
by solving all the fluid motions of the Navier-Stokes equations
because the range scales involved is so large that the computational
cost becomes prohibitive. In LES, only the large eddies are resolved,
and the effect of the small scales on the larger eddies is modeled
through an SGS model~\citep{sagaut2006} as illustrated in figure
\ref{fig:model:LES_sketch}. The approach enables a reduction of the
computational cost by several orders of magnitude while still
capturing the statistical quantities of interest.  In the present
section, we demonstrate the principle of maximum conservation of
information discussed in \S \ref{sec:modeling} by devising an SGS
model for LES.
%

The governing equations for LES are formally derived by applying a
spatial filter to Eq. (\ref{eq:cau:NS}),
\begin{eqnarray}\label{eq:model:LES}
\frac{\partial \bar u_i}{\partial t} + \frac{\partial \bar u_i \bar u_j  }{\partial x_j}
+ \frac{\partial \tau^{\mathrm{SGS}}_{ij}  }{\partial x_j}  = - \frac{1}{\rho}\frac{\partial \bar \Pi}{\partial x_i} +
\nu \frac{\partial^2 \bar u_i}{\partial x_j\partial x_j},
\quad\frac{\partial \bar u_i}{\partial x_i}  = 0,
\end{eqnarray}
where $\bar{(\cdot)}$ denotes spatially filtered quantity, and
$\tau^{\mathrm{SGS}}_{ij}$ is the effect of the subgrid scales on the
resolved eddies, which has to be modeled. The filter operator on a
variable $\phi$ is defined as
\begin{equation}\label{eq:model:filter_in}
\bar \phi(\boldsymbol{x},t) \equiv \int_V
G(\boldsymbol{x}-\boldsymbol{x'} ; \bar{\Delta})
\phi(\boldsymbol{x'},t) \mathrm{d}\boldsymbol{x'},
\end{equation}
where $G$ is the filter kernel with filter size $\bar{\Delta}$, and
$V$ is the domain of integration. The system in
Eq. (\ref{eq:model:LES}) is assumed to be severely truncated in the
number of the degrees of freedom with respect to
Eq. (\ref{eq:cau:NS}).  The objective of LES is to model the SGS
tensor as function of known filtered quantities,
\begin{subequations}
  \begin{gather}
    \tau^{\mathrm{SGS}}_{ij} = \tau^{\mathrm{SGS}}_{ij}( \bar S_{ij},\bar \Omega_{ij}, \bar \Delta ; \boldsymbol{\theta}),
\end{gather}
\end{subequations}
where $\bar S_{ij} = (\partial \bar{u}_i/\partial x_j + \partial
\bar{u}_j/\partial x_i)/2$, and $\bar \Omega_{ij}=(\partial
\bar{u}_i/\partial x_j - \partial \bar{u}_j/\partial x_i)/2$ are the
filtered rate-of-strain and rate-of-rotation tensors, respectively,
and $\boldsymbol{\theta}$ are model parameters.

In the present formulation, the map $\boldsymbol{f}$ in
Eq. (\ref{eq:model:Q_total}) corresponds to a discrete version of
Eq. (\ref{eq:cau:NS}), in which all the space and time scales are
accurately resolved. The state vector $\bQ{}{n}$ is given by the
discretization of $u_i$ and $\Pi$ in a grid fine enough to capture all
the relevant scales of motion. The map for the model,
$\widehat{\bs{f}}$, is derived from the discretization of
Eq. (\ref{eq:cau:NS}), and the model state $\bQhat{}{n}$ corresponds
to the filtered velocities and pressure, $\bar{u}_i$ and $\bar{\Pi}$.

A common misconception in LES modeling is that the closure problem of
determining $\tau^{\mathrm{SGS}}_{ij}$ arises from introducing the
filter operator. Interestingly, this is not entirely
accurate. Instead, the formalism introduced in \S \ref{sec:modeling}
shows that the closure problem is a consequence of the loss of
information introduced by the filter rather than the action of
filtering itself. This is easy to demonstrate by noting that the
analytic form of $\tau^{\mathrm{SGS}}_{ij}$ is completely determined
when the filter is reversible, i.e., the information is
conserved~\citep{yeo1987, carati2001, bae2017, bae2018}. One example
of reversible filter is given by the differential
filter~\citep{germano1986}
\begin{equation}
  \label{eq:model:SGS_exact}
G(\boldsymbol{x}-\boldsymbol{x}';\bar{\Delta}) = \frac{1}{4\pi
{\bar{\Delta}}^2}\frac{\exp(-|\boldsymbol{x}-\boldsymbol{x}'|/\bar{\Delta})}{|\boldsymbol{x}-\boldsymbol{x}'|},
\end{equation}
such that the analytic form of the $\tau_{ij}^{\mathrm{SGS}}$ is
exactly given by
\begin{equation}
  \label{eq:model_SGS_exact_e}
\tau_{ij}^{\mathrm{SGS}}= 
  \overline{\bar{u}_i \bar{u}_j} - \bar{\Delta}^2 \overline{\bar{u}_j \frac{\partial\bar{u}_i}{\partial x_k\partial x_k}}
  - \bar{\Delta}^2 \overline{ \bar{u}_i \frac{\partial\bar{u}_j}{\partial x_k\partial x_k}}
+  \bar{\Delta}^4 \overline{\frac{\partial\bar{u}_i}{\partial x_k\partial x_k} \frac{\partial\bar{u}_j}{\partial x_k\partial x_k}}
  - \bar{u}_i \bar{u}_j.
\end{equation}
Equation (\ref{eq:model_SGS_exact_e}) is a function of the filtered
velocities $\bar{u}_i$, which are accessible to the model, and does
not pose any closure problem. The actual closure problem emerges from
the application of irreversible filters and/or the coarse
discretization of the governing equations, which entail a truncation
of the number of degrees of freedom in the system. In those
situations, the LES grid resolution is unable to represent the small
scales (i.e., subgrid scales), which in turn entails a loss of
information. Hence, the paradigm of conservation of information
discussed in \S \ref{sec:modeling:formulation} arises naturally as a
fundamental aspect of LES modeling.

\begin{figure}
  \begin{center}
     \includegraphics[width=1.07\textwidth]{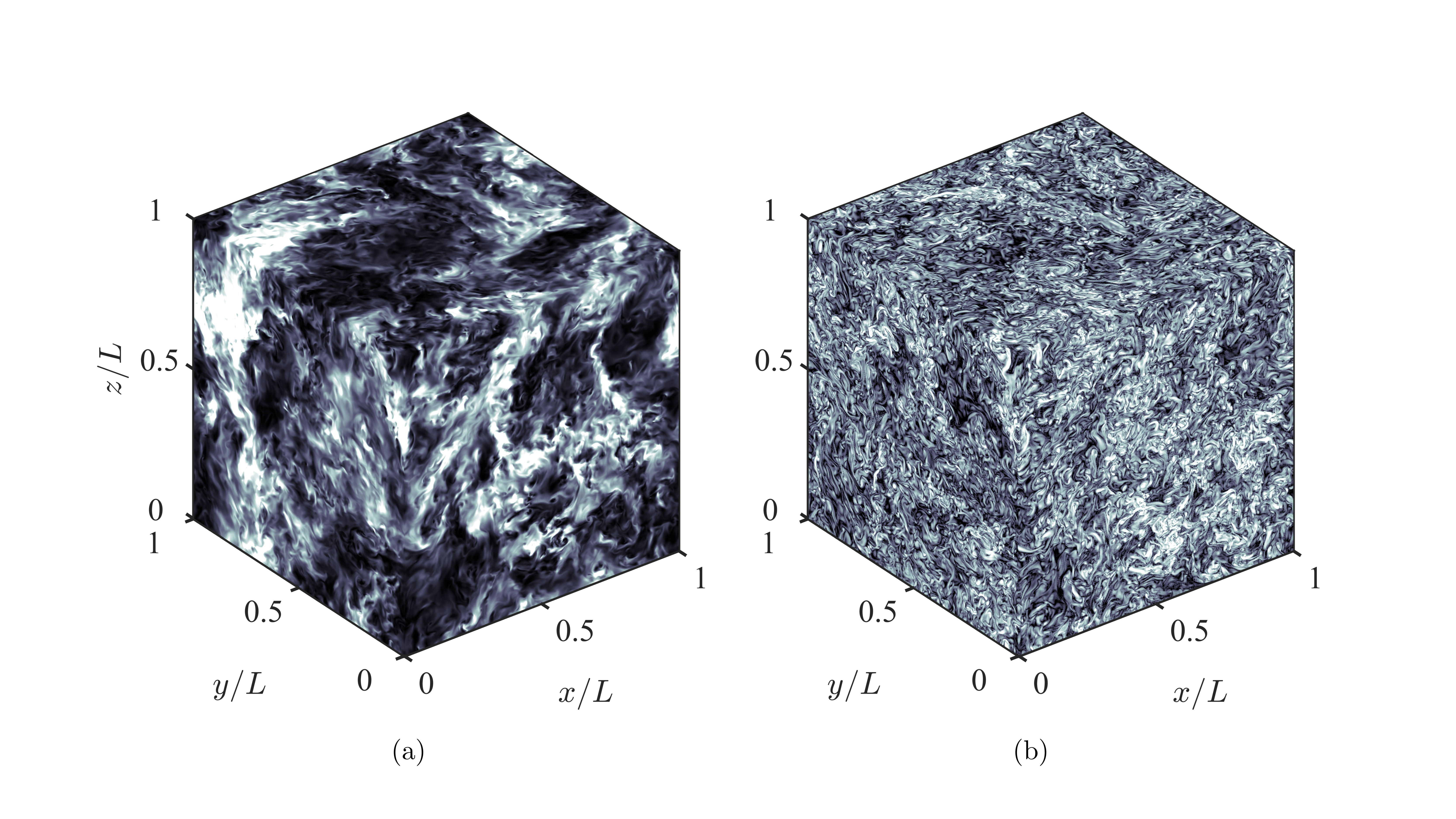} 
 \end{center}
\caption{ Visualization of (a) the instantaneous turbulent kinetic
  energy $u_i u_i/2$ and (b) enstrophy $\omega_i \omega_i$, where
  $\omega_i$ is the vorticity. The colormap ranges from 0 (dark) to
  0.8 of the maximum value (light) of turbulent kinetic energy and
  enstrophy, respectively.
  \label{fig:model:HIT_u}}
\end{figure}

We leverage Eq. (\ref{eq:model:summary_2}) to construct an SGS model
for turbulent flows. The flow considered is similar to that presented
in \S \ref{sec:cau:example}: forced isotropic
turbulence~\citep{lesieur2008} in a triply periodic, cubic domain with
size equal to $L$ as shown in figure \ref{fig:model:HIT_u}. The exact
solution is obtained from a DNS using $512^3$ dealiased Fourier modes
for the spatial discretization and a fourth-order Runge-Kutta time
stepping method. The turbulence is sustained by adding a linear
forcing to the right-hand side of Eq. (\ref{eq:cau:NS}) equal to $f_i
= A u_i$, where $A$ was adjusted to maintain on average $Re_\lambda
\approx 260$. The simulation was run for 50 integral times after
initial transients. Figure \ref{fig:model:HIT_u} features a
visualization of the energy field and dissipation field from the DNS,
highlighting the separation of scales in the system.
%
The functional form considered for the SGS stress tensor is
\begin{equation}
  \label{eq:model:SGS_model}
  \tau_{ij}^{\mathrm{SGS}} - \frac{1}{3} \tau_{kk}^{\mathrm{SGS}} \delta_{ij} =
  \theta_1 \bar{\Delta}^2 \bar{S}_{ij} \sqrt{\bar{S}_{nm}\bar{S}_{nm}}  +
  \theta_2 \bar{\Delta}^2 (\bar{S}_{ik}\bar{\Omega}_{kj} - \bar{\Omega}_{ik}\bar{S}_{kj}),
\end{equation}
where $\delta_{ij}$ is the Kronecker delta, and $\theta_1$ and
$\theta_2$ are modeling parameters to be determined.  Equation
(\ref{eq:model:SGS_model}) is derived by retaining the two leading
terms from the general expansion of the SGS tensor in terms of
$\bar{S}_{ij}$ and $\bar{\Omega}_{ij}$ proposed by
\citet{lund1992}.

\begin{figure}
  \begin{center}
    \subfloat[]{\includegraphics[width=0.47\textwidth]{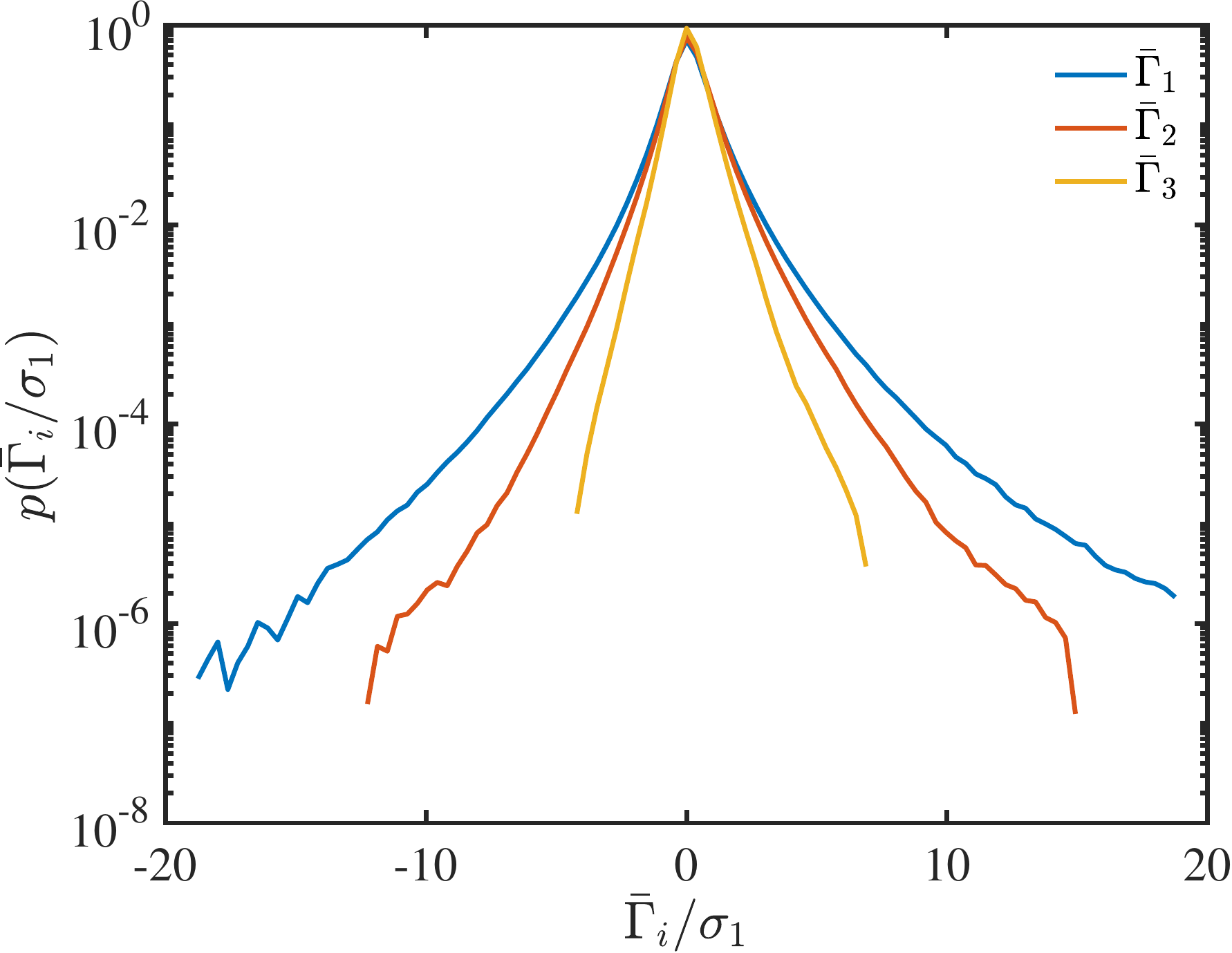}}
    \hspace{0.05cm}
    \subfloat[]{\includegraphics[width=0.47\textwidth]{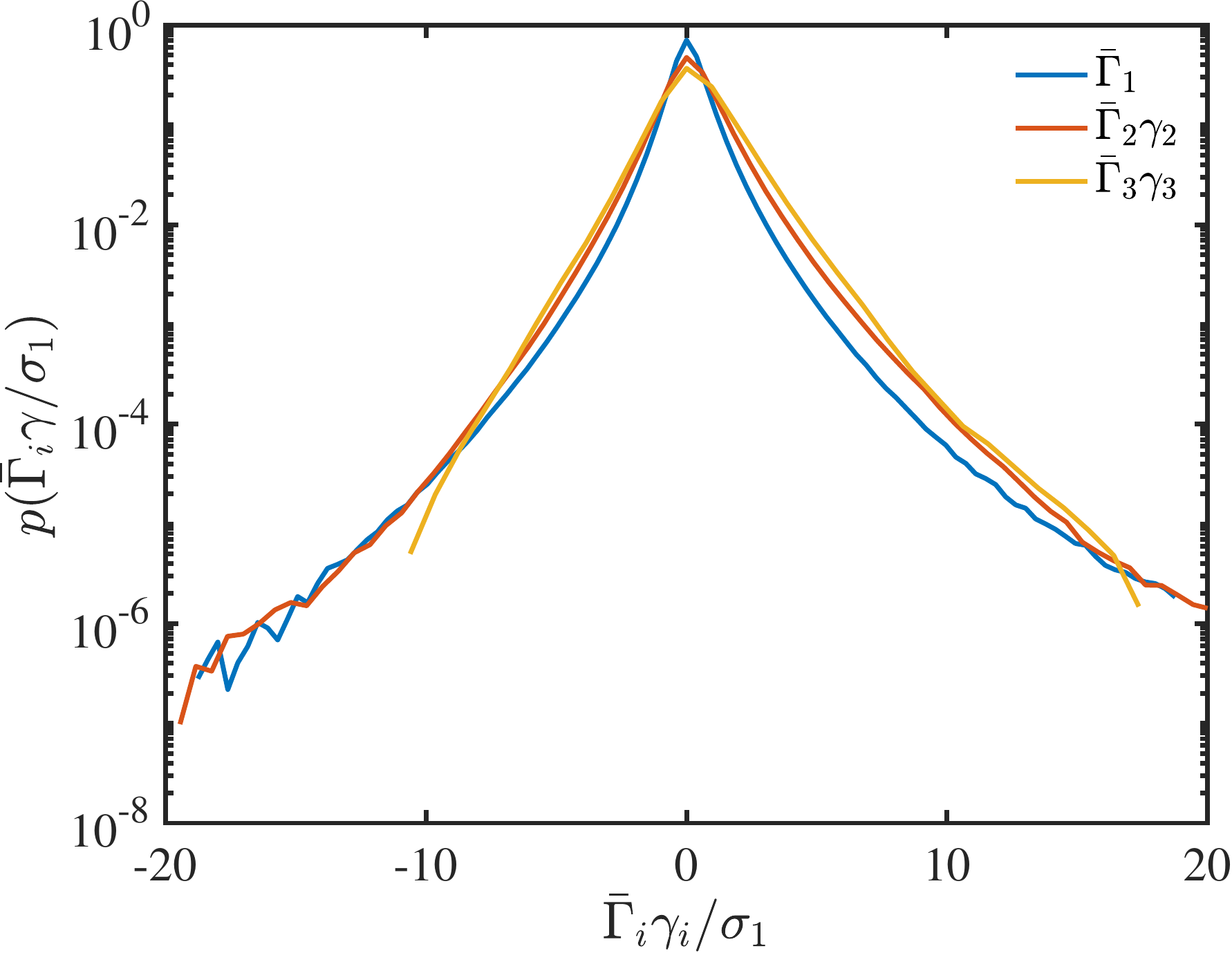}} 
 \end{center}
\caption{ (a) Probability mass distributions of the interscale energy
  transfer and viscous dissipation, $\bar{\Gamma}_1$,
  $\bar{\Gamma}_2$, and $\bar{\Gamma}_3$, at filter cut-offs
  $\bar{\Delta}_1=L/32$, $\bar{\Delta}_2=L/16$, and
  $\bar{\Delta}_3=L/8$, respectively. The results are for DNS using
  sharp Fourier filter and are normalized by the standard deviation of
  $p(\bar{\Gamma}_1)$, denoted by $\sigma_1$. (b) Probability mass
  distributions of the rescaled interscale energy transfer and viscous
  dissipation, $\bar{\Gamma}_1$, $\bar{\Gamma}_2\gamma_2$,
  $\bar{\Gamma}_3\gamma_3$ with $\gamma_2 =
  (\bar{\Delta}_2/\bar{\Delta}_1)^{2/3}$ and $\gamma_3 =
  (\bar{\Delta}_3/\bar{\Delta}_1)^{2/3}$.
  \label{fig:model:test_scaling}}
\end{figure}

Let us introduce the interscale energy transfer and viscous dissipation
at the filter cut-off $\bar{\Delta}$ given by
\begin{equation}
  \overline{\Gamma} = ( \overline{u_i u_j} - \overline{u}_i \overline{u}_j) \overline{S}_{ij} - 2\nu\overline{S}_{ij}\overline{S}_{ij}.
\end{equation}
The modeling assumption proposed here is that the information content
of $p(\bar{\Gamma}_1)$ must be equal to the information content of
$p(\bar{\Gamma}_2 \gamma)$, where $\bar{\Gamma}_1$ and
$\bar{\Gamma}_2$ are $ \overline{\Gamma}$ at two different scales
$\bar{\Delta}_1$ and $\bar{\Delta}_2$, respectively, and $\gamma =
(\bar{\Delta}_1/\bar{\Delta}_2)^{2/3}$ is a scaling factor. This
self-similarity in the information implies that the energy transfer at
two different scales should satisfy
\begin{equation}
 \label{eq:model:SGS_condition}
 p(\bar{\Gamma}_1) \approx
 p(\bar{\Gamma}_2 \gamma )/\gamma.
\end{equation}
The hypothesis in Eq. (\ref{eq:model:SGS_condition}) is corroborated
in figure \ref{fig:model:test_scaling} using DNS data. Figure
\ref{fig:model:test_scaling}(a) shows $p(\bar{\Gamma}_i)$ for three
different values of the filter widths: $\bar{\Delta}_1=L/32$,
$\bar{\Delta}_2=L/16$, and $\bar{\Delta}_3=L/8$.  The scaling
condition from Eq. (\ref{eq:model:SGS_condition}) is tested in figure
\ref{fig:model:test_scaling}(b), which reveals the improved collapse
using the factor $\gamma$. A similar scaling result was observed by
\citet{aoyama2005}.
%

In the case of LES, the interscale energy transfer and dissipation
also depends on the contribution of the SGS model as
\begin{equation}
  \bar{\Gamma} =  ( \overline{u_i u_j} - \overline{u}_i \overline{u}_j) \overline{S}_{ij} -
  2\nu\bar{S}_{ij}\bar{S}_{ij} +  \tau_{ij}^{\mathrm{SGS}}\bar{S}_{ij}.
\end{equation}
Hence, the model proposed aims at minimizing the information lost when
$p(\bar{\Gamma}_1)$ is used to approximate $p(\bar{\Gamma}_2 \gamma)$
in the LES solution with $\bar{\Delta}_1 = 2 \bar{\Delta}$ and
$\bar{\Delta}_2 = 2 \bar{\Delta}_1$. 
The model is formulated using the KL divergence, which
ensures that the average information required for reconstructing
$p(\bar{\Gamma}_2 \gamma )$ is minimum given the information in
$p(\bar{\Gamma}_1)$,
\begin{equation}
  \label{eq:model:SGS_theta}
  \boldsymbol{\theta} = \arg\min_{\boldsymbol{\theta}'}
  \mathrm{KL}\left( \bar{\Gamma}_2  \gamma , \bar{\Gamma}_1 \right),
\end{equation}
where $\boldsymbol{\theta} = (\theta_1,\theta_2)$ from
Eq. (\ref{eq:model:SGS_model}). We will refer to the model as
Information-Preserving SGS model, or as IP-SGS model for
short. \corr{Note that the IP-SGS model only relies on the physical
  assumption that the information content of $p(\bar{\Gamma}_1)$ is
  equal to the information content of $p(\bar{\Gamma}_2 \gamma)$, and
  does not require any DNS data to be trained.}

To validate the model, an LES is conducted using $64^3$ Fourier
modes. The turbulence is driven by a linear forcing with the same $A$
value obtained for the DNS. The kernel selected is the sharp Fourier
filter. The LES entails a severe truncation of the number of degrees
of freedom of the original system: the DNS contains more than 130
million degrees of freedom, whereas in the LES system the number of
degrees of freedom is reduced to only 0.3 million. The IP-SGS model is
implemented as follows: during the LES runtime, statistics are
collected on-the-fly to reconstruct the probability distributions of
$\bar{\Gamma}_2 \gamma$ and $\bar{\Gamma}_1$. Every $100$ time steps,
the model parameters $\theta_1$ and $\theta_2$ are computed from
Eq. (\ref{eq:model:SGS_theta}) using a gradient descent method that
minimizes the KL divergence.

\begin{figure}
  \begin{center}
    \subfloat[]{\includegraphics[width=0.55\textwidth]{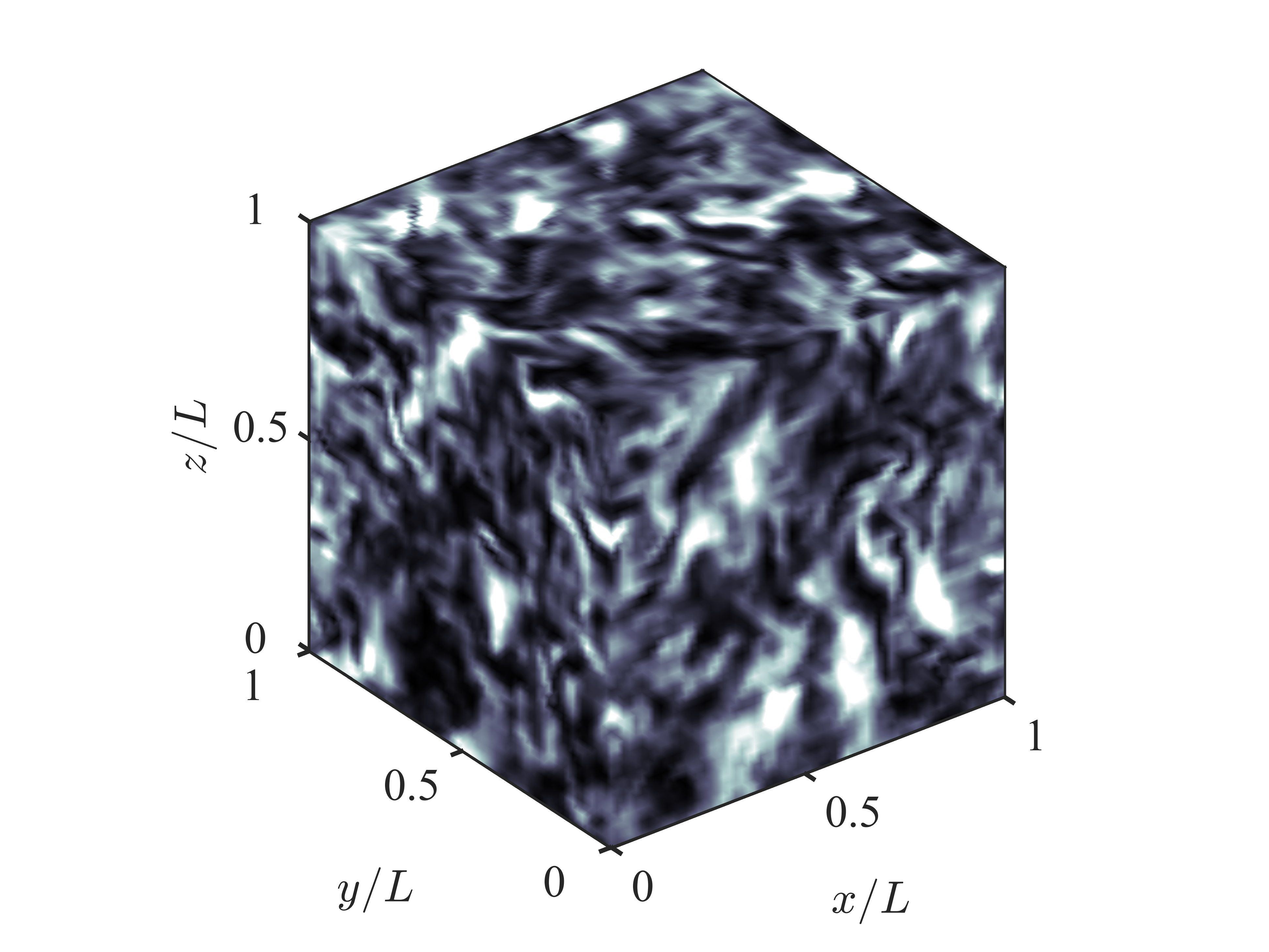}} 
    \subfloat[]{\includegraphics[width=0.44\textwidth]{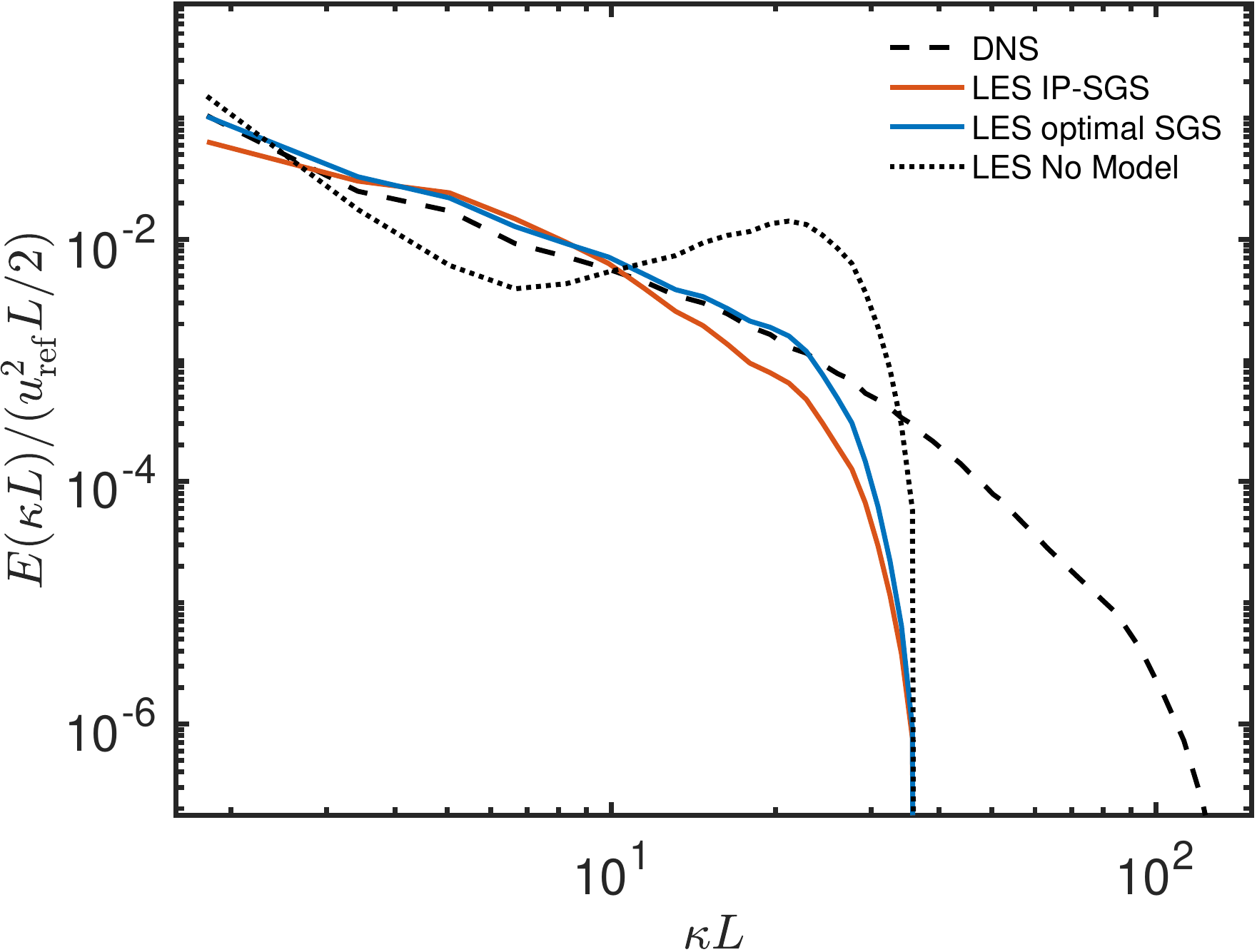}} 
 \end{center}
\caption{ (a) Visualization of the instantaneous turbulent kinetic
  energy $u_i u_i/2$.  The colormap ranges from 0 (dark) to 0.8 of the
  maximum value (light) of turbulent kinetic energy. (b) The kinetic
  energy spectra as a function of the wavenumber for the `exact' DNS
  solution, LES with IP-SGS model, LES with optimal SGS model, and LES
  without SGS model. The energy spectra are normalized by
  $u_\mathrm{ref}^2/2 L$, where $u_\mathrm{ref}^2/2$ is the mean
  kinetic energy of the DNS solution.
  \label{fig:model:results}}
\end{figure}

The performance of the SGS model is evaluated in figure
\ref{fig:model:results} after initial transients in the system.  We
use as figure of merit the kinetic energy spectrum $E(\kappa)$,
where $\kappa$ is the wavenumber. The predictions are compared against
a case without SGS model and the optimal model. The latter is defined
as the SGS model of the form dictated by
Eq. (\ref{eq:model:SGS_model}) with the values of $\theta_1$ and
$\theta_2$ that yield the best prediction of $E(\kappa)$ in the
L$_2$-norm sense. The optimal values of $\theta_1$ and $\theta_2$ were
obtained by a parametric sweep. The results in figure
\ref{fig:model:results} show that the IP-SGS model offers an accuracy
comparable to the optimal model, while providing a non-trivial
improvement with respect to the case without SGS model.  This modeling
exercise demonstrates the viability of the information-theoretic
formulation presented in \S \ref{sec:modeling} as an effective
framework for reduced-order modeling of highly chaotic systems with
large number of degrees of freedom.
%

\section{Control}
\label{sec:control}

%
The first entropic analysis of feedback control for dynamical systems
was proposed by \citet{weidemann1969}, who envisioned the
sensor/actuator device as an information transformation of the data
collected.  The optimal control of time-continuous systems was later
attempted by \citet{saridis1988}, and further extended to discrete
systems by \citet{tsai1992}. \citet{tatikonda2004} investigated the
lower bounds of controllability, observability, and stability of
linear systems under communication constraints between the sensor and
the actuator.  \citet{touchette2004} substantially advanced the
information-theoretic formulation of the control problem by redefining
the concepts of controllability and observability using conditional
entropies.  They posed the problem of control as a task of
entropy-reduction, and proved that the maximum reduction of entropy
achievable in a system after a one-step actuation is bounded by the
mutual information between the control parameters and the current
state of the system.  The analysis by \citet{touchette2004} was
broaden by \citet{delvenne2013} to more general, convex, cost
functions.  In the same vein, \citet{bania2020} recently developed an
information-aided control approach which provides results akin to
those from dynamic programming but at a more affordable computational
cost.  \corr{A review on the thermodynamics of information
  and its connection with feedback control can be found in
  \citet{parrondo2015}}. Despite the merits of the previous works, a
notable caveat is their limited applicability to stochastic control of
problems with a low number of degrees of freedom~\citep[see the
  examples in][]{chen2019}.

In this section, we provide an information-theoretic formulation of
the problem of optimal control for chaotic, high-dimensional dynamical
systems. \corr{New definitions of open/closed-loop control,
  observability, and controllability are introduced in terms of the
  mutual information between different states of the system. The task
  of optimal control is posed as the reduction in uncertainty of the
  controlled state given the information collected by the sensors and
  the action performed by the actuators.  In contrast to the
  traditional formulation of control, which emphasizes the
  differential equations of the dynamical system, our formulation is
  centered in the probability distribution of the states.} The theory
is applied to achieve optimal drag reduction in a wall-bounded
turbulent flow using opposition control at the wall.

\subsection{Formulation}
%
%
Let us denote the state of the system to be controlled at time $t_n$
by $\bQ{}{n}=[\Q{1}{n},...,\Q{N}{n}]$, where $N$ are the total number
of degrees of freedom.  The state of the uncontrolled system is
denoted by $\bQ{u}{n}$ (with subscript $u$) and its dynamics is
completely determined by
\begin{equation}
    \bQ{u}{n+1} = \boldsymbol{f}(\bQ{u}{n}),
\end{equation}
where $\bs{f}$ is the map function introduced in \S
\ref{sec:information_dynamical}.  The system is controlled to a new
state $\bQ{}{n+1}$ (without subscript $u$) by means of a sensor
and an actuator that together constitute the controller.  
The state of the sensor and actuator are denoted by $\mybS$ and
$\mybA$, respectively, and both are considered random variables. 
The properties of the controller are parametrized by the vector $\mybth$ 
(e.g., actuator/sensor locations, actuator intensity and frequency, etc.).  
In general, the full state $\bQ{}{n}$ is inaccessible and only a subset of the 
phase-space is observable by the sensor,
\begin{equation}
    \bS{}{n} = \bs{h}(\bQ{}{n}, \bWs{}{n}; \mybth), \label{eq:sensor_Sn}
\end{equation}
where $\bWs{}{n}$ represents random noise in the
measurements. We will consider that $\bWs{}{n}$ is uncorrelated
with $\bQ{}{n}$. The noise $\bWs{}{n}$ introduces
additional information into the system that can be labeled as spurious,
since it masks the actual information from the state
$\bQ{}{n}$. The actuator gathers the information from the
sensor and acts according to the control law
\begin{equation}
    \bA{}{n} = \bs{g}(\bS{}{n}, \bWa{}{n}; \mybth), \label{eq:controller_An}
\end{equation}
where $\bWa{}{n}$ is an auxiliary random variable which provides additional 
stochasticity (hence, information) to the actuator and is independent of
$\bQ{}{n}$. 
Then, the controlled system is governed by
\begin{equation}
      \bQ{}{n+1} = \bs{f}(\bQ{}{n},\bA{}{n}). \label{eq:controller:new}
\end{equation}

\begin{figure}
    \centering
    \subfloat[\label{fig:control:Qn1}]{\colorlet{c1}{myc1}
\colorlet{c2}{myc3}
\colorlet{c3}{myc4}

\colorlet{c1s}{c1}
\colorlet{c2s}{c2}
\colorlet{c3s}{c3}

\begin{tikzpicture}[scale=1.5,thick,>={Latex[length=.2cm]}]
    
    \begin{footnotesize}

    \pgfmathsetmacro{\rdi}{1}           
    \pgfmathsetmacro{\rdiO}{.9*\rdi}   
    \pgfmathsetmacro{\rdiOm}{.7*\rdiO}  
    \pgfmathsetmacro{\dis}{.7}          

    \coordinate (O1) at (-\dis,0);  
    \coordinate (O2) at (\dis,0);

    \fill[c1s,opacity=.3] (O1) circle (\rdi);
    \fill[c2s,opacity=.3] (O2) circle (\rdi);
    \fill[c3s,opacity=.3] ellipse [x radius=\rdiO, y radius=\rdiOm];

    \draw[c1] (O1) circle (\rdi);
    \draw[c2] (O2) circle (\rdi);
    \draw[c3,name path=ellip] ellipse [x radius=\rdiO, y radius=\rdiOm];
  
    \path (0,-\rdi) --+ (0,-1em);

    \draw[thin,<-] (O1)+(110:\rdi) --++ (118:1.5*\rdi) node [fill=white] {$H(\bQ{}{n})$};
    \draw[thin,<-] (O2)+(70:\rdi) --++ (62:1.5*\rdi) node [fill=white] {$H(\bA{}{n})$};

    \path [name path=line1] (0,0) -- (-30:1.5*\rdi);
    \path [name intersections={of=ellip and line1,by=E}];
    \draw[<-,thin,shorten >=1em] (E) --++ (-35:1.1*\rdi)  
        node[anchor=center]{$H(\bQ{}{n+1})$};

    \end{footnotesize}

\end{tikzpicture}}~\hfill
    \subfloat[\label{fig:control:act}]{\colorlet{c1}{myc1}
\colorlet{c2}{myc2}
\colorlet{c3}{myc3}

\colorlet{c1s}{c1}
\colorlet{c2s}{c2}
\colorlet{c3s}{c3}

\begin{tikzpicture}[scale=1.5,thick,>={Latex[length=.2cm]}]
    
    \begin{footnotesize}

    \pgfmathsetmacro{\rdi}{1}           
    \pgfmathsetmacro{\rdiO}{1.1*\rdi}   
    \pgfmathsetmacro{\rdiOm}{.7*\rdiO}  
    \pgfmathsetmacro{\dis}{.4}          

    \coordinate (O1) at (-\dis,0);  
    \coordinate (O2) at (\dis,0);

    \pgfmathsetmacro{\hei}{sqrt(\rdi^2-\dis^2)}

    \fill[c1s,opacity=.3] (O1) circle (\rdi);
    \fill[c2s,opacity=.3] (O2) circle (\rdi);
    \coordinate (rotP) at (\rdiOm,\hei);
    \begin{scope}[shift={(rotP)}]
        \fill[c3s,rotate around={-30:(-\rdiOm,0)},opacity=.3] 
        ellipse [x radius=\rdiOm, y radius=\rdiO];
    \end{scope}

    \draw[c1] (O1) circle (\rdi);
    \draw[c2] (O2) circle (\rdi);

    \begin{scope}[shift={(rotP)}]
        \draw[c3,rotate around={-30:(-\rdiOm,0)},name path=ellip] 
        ellipse [x radius=\rdiOm, y radius=\rdiO];

        \path [name path=line1] (0,0) -- (0:\rdi);
        \path [name intersections={of=ellip and line1,by=E}];
        \draw[<-,thin] (E) --++ (-10:.65*\rdi)  
            node[anchor=center,fill=white]{$H(\bA{}{n})$};
    \end{scope}

    \draw[thin,<-] (O1)+(110:\rdi) --++ (118:1.5*\rdi) node [fill=white] {$H(\bQ{}{n})$};
    \draw[thin,<-] (O2)+(-45:\rdi) --++(-35:1.6*\rdi) node [fill=white] {$H(\bS{}{n})$};

    \draw[thin,<-] (50:.3*\rdi) -- (0:2.05*\rdi) node [fill=white] 
        {$I(\bA{}{n};\bQ{}{n})$};

    \draw[thin,<-,shorten >=1.5em] (220:.4*\rdi) -- (210:1.72*\rdi) node 
        {$I(\bS{}{n};\bQ{}{n})$};

    \path (0,-\rdi) --+ (0,-1em);

    \end{footnotesize}

\end{tikzpicture}}
    \caption{Schematic of the entropies of the system in
      Eq.~\eqref{eq:controller:new}.  (a) Entropy of the new controlled
      state, ${\protect\bQ{}{n+1}}$, which is bounded by the entropies
      of ${\protect\bQ{}{n}}$ and ${\protect\bA{}{n}}$, as inferred from
      Eq.~\eqref{eq:controller:new}.  (b) Relationship among the
      entropies of the current state, the sensor, and the actuator.}
\end{figure}
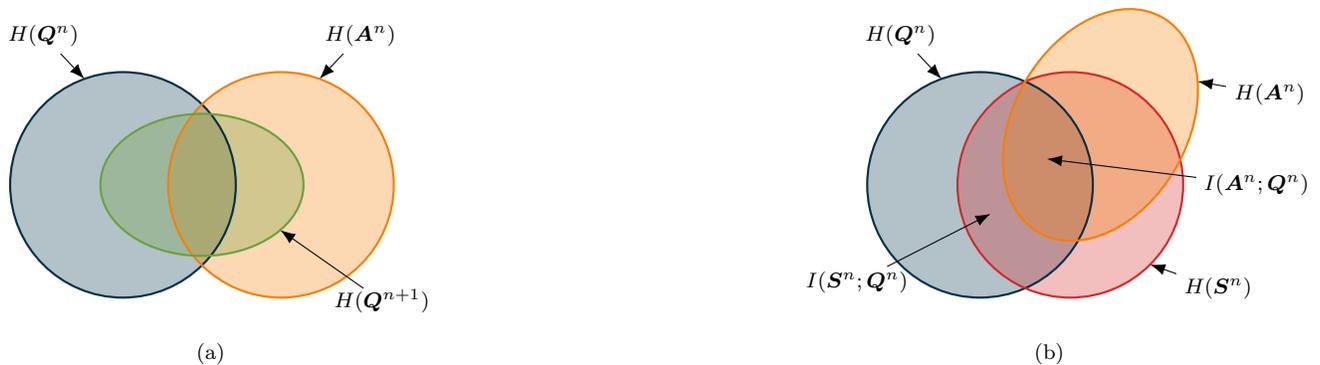

In the case of high-dimensional, chaotic systems, control of the full
state constitutes an impractical task and hence is not the main
concern here. Instead, our goal is to control a few degrees of freedom
which are the most impactful on reducing or enhancing a quantity of
interest. The sub-state to be controlled is denoted by $\bJ{}{n}$, and
is also assumed to be a random variable derived from $\bQ{}{n}$,
\begin{equation}
  \label{eq:control:J}
  \bJ{}{n} = \bs{l}(\bQ{}{n}).
\end{equation}
In many situations, we are interested in the controlled state once the
system has reached the statistically steady state.
In those cases, the time step $n+1$ represents the state of the system
after initial transients.
As an example (in line with the application in \S
\ref{sec:control:example}), we can consider the reduction of drag in
an airfoil by blowing and suction of air over its surface.
%
In this case, $\bQ{}{n}$ is a (high-dimensional) discrete
representation of the velocity and pressure fields in all the domain
surrounding the airfoil, $\bS{}{n}$ are pressure probes at the airfoil
surface, $\bA{}{n}$ is a flow jet at the wall which modifies the air
velocity around the airfoil, and $\bJ{}{n+1}$ represents the
controlled (low-dimensional) drag state after transients.  The goal of
the control law is then to alter the probability distribution of the
drag to reduce \emph{i)} its mean value and \emph{ii)} its standard
deviation to mitigate extreme drag events.

\subsubsection{Open-loop and closed-loop control}

The information in the controlled system flows from the system state
$\bQ{}{n}$ to the sensor $\bS{}{n}$, and from the sensor to the
actuator $\bA{}{n}$.  In the general scenario, the sensor shares the
information with the actuator via a communication channel.  The
capacity of the communication channel between $\bS{}{n}$ and
$\bA{}{n}$ is defined as
\begin{equation}
  \label{eq:control:capacity}
  \mathrm{Ca} = \max_{p(\myindexvar{\bs{s}}{}{n})} I(\bS{}{n};\bA{}{n}),
\end{equation}
where the maximum is taken over all possible input distributions
$p(\myindexvar{\bs{s}}{}{n})$. By virtue of the noisy-channel coding
theorem~\citep{shannon1948}, the capacity in
Eq.~\eqref{eq:control:capacity} provides the highest information rate
(i.e., bits per second) that can be achieved with an arbitrarily small
error probability between the sensor and actuator.

The mutual information between the sensor and actuator from
Eq.~\eqref{eq:control:capacity} provides the grounds for the
definition of open-loop and closed-loop controllers.  A control is
said to be open-loop if
\begin{equation}
 \label{eq:control:open}
I(\bA{}{n};\bQ{}{n})=0,
\end{equation}
i.e., there is no shared information between the actuator and the
state of the system at $t_n$.
%
%
%
Note that an open-loop control can still modify the information content of future 
states because $I(\bQ{}{n+1};\bA{}{n})\geq 0$ by means of the auxiliary
random variable $\bWa{}{n}$. 
Conversely, a control is closed-loop if
\begin{equation}
 \label{eq:control:closed}
I(\bA{}{n};\bQ{}{n})>0,
\end{equation}
namely, the current  state and the actuator share an amount of
 information greater than zero.
From Eq.~\eqref{eq:sensor_Sn} and \eqref{eq:controller_An} it can be shown that
\begin{equation}
I(\bA{}{n};\bQ{}{n}) \leq I(\bS{}{n};\bQ{}{n}).
\end{equation}
Thus, a sufficient condition for open-loop control is
$I(\bS{}{n};\bQ{}{n}) = 0$, whereas $I(\bS{}{n};\bQ{}{n}) > 0$ and
$I(\bA{}{n};\bS{}{n}) > 0$ are necessary conditions for closed-loop
control.
A constraint from Eq. (\ref{eq:controller_An}) is that, if there is no
auxiliary random noise ($\bWa{}{n}$), the actuator cannot contain more
information than the sensor, $H(\bA{}{n})\leq H(\bS{}{n})$.
The relationships among $\bQ{}{n}$, $\bS{}{n}$ and $\bA{}{n}$ are
illustrated in Figure~\ref{fig:control:act}.

It is also interesting to establish the information loss across the
controller. Let us assume an sparse sensor with $N_s$ degrees of
freedom such that $N_s \ll N$. The information content of the noise
scales as $H(\bWs{}{n})\sim \log N_s$, while the information of the
system generally follows $H(\bQ{}{n}) \sim \log N$. It is
then reasonable to assume that $H(\bWs{n}{}) \ll H(\bS{}{n})$ will
hold in most practical situations where $N_s \ll N$. If we further
assume that all the information in the actuator is obtained from the
sensor $\bS{}{n}$, then the hierarchy of information loss across the
controlled system is given by
\begin{equation}
  \label{eq:control:info_loss}
 H(\bQ{}{n}) \gg H(\bS{}{n}) \geq H(\bA{}{n}).
\end{equation}
Equation (\ref{eq:control:info_loss}) shows that, when the
measurements are sparse, the state of the system contains more
information than the sensor, which in turn contains similar or less
information than the actuator. The first inequality in
Eq. (\ref{eq:control:info_loss}) might not hold for systems with a few
degrees of freedom or a large number of sensors. In those situations,
the noise could increase the (spurious) information content of the
sensor to yield $H(\bQ{}{n}) < H(\bS{}{n})$. Nonetheless, here we are
interested in high-dimensional systems controlled using a few
measurements such that Eq. (\ref{eq:control:info_loss}) is likely to
hold.

\subsubsection{Observability and controllability}
\label{subsec:control:obs_con}

Observability (how much we can know about the system) and
controllability (how much we can modify the system) represent two
major pillars of modern control system theory. Here, we formulate the
information-theoretic counterparts of observability and
controllability. Our definitions are motivated by the no-uncertainty
conditions of deterministic systems given by
\begin{subequations}
  \label{eq:control:motivation}
  \begin{gather}
    H(\bJ{}{n} |   \bQ{}{n} ) = 0,   \label{eq:control:motivation1} \\
    H(\bJ{}{n+1} | \bQ{}{n},\bA{}{n}) = 0.   \label{eq:control:motivation2}
  \end{gather}
\end{subequations}
Equation (\ref{eq:control:motivation1}) is a statement about the
observability of $\bJ{}{n}$: there is no uncertainty in $\bJ{}{n}$
when the full state of the system in known.  However, uncertainties
might arise when the information available is only limited to the
state of the sensor. Similarly, Eq. (\ref{eq:control:motivation2})
relates to the concept of controllability: given the system state
$\bQ{}{n}$ and the actuator action $\bA{}{n}$, there is no uncertainty
in the future state $\bJ{}{n+1}$.

Consistently with the remarks above, observability of the state
$\bJ{}{n}$ with respect to the sensor $\bS{}{n}$ is defined as
\begin{equation}
  \label{eq:control:obs}
  O_J \equiv \frac{I(\bJ{}{n} ; \bS{}{n} )}{H(\bJ{}{n})},
\end{equation}
which represents the uncertainty in $\bJ{}{n}$ given the information
from the sensor normalized by the total information in $\bJ{}{n}$ such
that $0 \leq O_J \leq 1$.  \corr{The normalization by $H(\bJ{}{n})$ in
  Eq. (\ref{eq:control:obs}) is introduced to provide a relative
  measure of the observability with respect to the total information
  in $\bJ{}{n}$.} The observability can also be expressed as a
function of $H(\bJ{}{n} | \bS{}{n})$ by
\begin{equation}
  \label{eq:control:obs_2}
  O_J = 1 - \frac{H(\bJ{}{n} | \bS{}{n} )}{H(\bJ{}{n})},
\end{equation}
which shows that the smaller the value of $H(\bJ{}{n} | \bS{}{n} )$
(i.e., the uncertainty in the state $\bJ{}{n}$), the larger the value
of $O_J$, in line with the intuition of observability argued in
Eq. (\ref{eq:control:motivation1}).
We say that a system with targeted variable $\bJ{}{n}$ is
\emph{perfectly observable} with respect to the sensor $\bS{}{n}$ if
and only if there is no uncertainty in the state of the system
conditioned to knowing the state of the sensor, namely,
\begin{equation}
  \label{eq:control:obs_3}
  H(\bJ{}{n} | \bS{}{n} ) = 0,
\end{equation}
which corresponds to $O_J = 1$. Conversely, $O_J=0$ if none of the
information in $\bJ{}{n}$ is accessible to the sensor.  It can be
shown that in the presence of noise in the sensor, the upper bound for
observability is reduced as
\begin{equation}
  \label{eq:control:obs_4}
  0 \leq O_J \leq 1 -  \frac{I(\bWs{}{n};\bS{}{n})}{H(\bJ{}{n})}.
\end{equation}
Thus, the observability of the controller degrades proportionally to
the noise contamination of the sensor as quantified by
$I(\bWs{}{n};\bS{}{n})$, and perfect observability is unattainable
when $I(\bWs{}{n};\bS{}{n})>0$.

Controllability of the future state $\bJ{}{n+1}$ with respect to the
action $\bA{}{n}$ is defined by
\begin{equation}
    \label{eq:control:con}
    C_J \equiv
    \frac{I(\bJ{}{n+1} ;  \bA{}{n} )}{H(\bJ{}{n+1})}, 
\end{equation}
provided that for all future states $\myindexvar{\bs{j}}{}{n+1}$ and
initial condition $\myindexvar{\bs{q}}{}{n}$ there exist a control
$\myindexvar{\bs{a}}{}{n}$ such that
$p(\myindexvar{\bs{j}}{}{n+1}|\myindexvar{\bs{q}}{}{n},\myindexvar{\bs{a}}{}{n})\neq
0$. Equation (\ref{eq:control:con}) quantifies the uncertainty in the
targeted state $\bJ{}{n+1}$ knowing the present information from the
control. Similar to the observability, $C_J$ is bounded by $0\leq C_J
\leq 1$ and can also be cast as
\begin{equation}
    \label{eq:control:con_1}
    C_J = 1 -
    \frac{H(\bJ{}{n+1} |  \bA{}{n} )}{H(\bJ{}{n+1})}, 
\end{equation}
The smaller the value of $H(\bJ{}{n+1} | \bA{}{n} )$ the larger the
controllability of the system (i.e., less uncertainty in the future
outcome). A system is \emph{perfectly controllable} at state
$\bJ{}{n+1}$ if and only if the uncertainty associated with the latter
upon application of the control action $\bA{}{n}$ is zero, i.e., there
exists a non-empty set of control values such that
\begin{equation}
    \label{eq:control:con_2}
    H(\bJ{}{n+1} | \bA{}{n} ) = 0,
\end{equation}
which corresponds to $C_J = 1$. Note that observability depends on the
targeted state at a given time ($\bJ{}{n}$), whereas controllability
relates the future targeted state ($\bJ{}{n+1}$) with the knowledge of
the system in the state at time $t_n$. Overall, the highest
observability and controllability are attained by maximizing the
mutual information between the target state and the control, which
will be leveraged in \S \ref{subsec:control:optimal} when seeking
optimal control strategies.

It is insightful to interpret the definitions in
Eq. (\ref{eq:control:obs}) and Eq. (\ref{eq:control:con}) as the
answer to the question: how much additional information is needed to
completely determine the state $\bJ{}{n}$ at time $t_n$
(observability) and at future times $t_{n+1}$ (controllability)
considering that the state of the control is known at time $t_n$.  The
aforementioned statement can literally be translated as the number of
bits (for example, the size of a digital file) that are required on
average to obtain perfect observability and controllability of the
system. The amount of missing information (MI) to achieve perfect
observability of the state $\bJ{}{n}$ is given by
\begin{equation}
  \mathrm{MI}_O = (1-O_J)H(\bJ{}{n}).
\end{equation}
Analogously, the amount of missing information
for perfect controllability of the state $\bJ{}{n+1}$ is
\begin{equation}
    \mathrm{MI}_C = (1-C_J)H(\bJ{}{n+1}).
\end{equation}
For example, if the observability of the control is $O_J=0.8$ and the
state $\bJ{}{n}$ has 120 megabytes of information, then the control
requires 24 megabytes of additional information to unambiguously
determine $\bJ{}{n}$. A similar example applies to
controllability. Another perspective on $\mathrm{MI}_O$ and
$\mathrm{MI}_C$ is that they signify the minimum number of yes and no
questions about the state $\bJ{}{n}$ that must be asked on average in
order to attain perfect observability and controllability of the
system, respectively.

\subsubsection{Optimal control}
\label{subsec:control:optimal}
Let us consider the quantity of interest $\bJ{}{n+1}$ with 
mean vector and covariance matrix given by
\begin{subequations}
  \label{eq:control:mu_sigma}
\begin{gather}
\mathbb{E}[\bJ{}{n+1}] = \bmun, \\
\mathrm{var}[\bJ{}{n+1}] = \mathbb{E}[(\bJ{}{n+1}-\mathbb{E}[\bJ{}{n+1}])(\bJ{}{n+1}-\mathbb{E}[\bJ{}{n+1}])^T] = \bsgn,
\end{gather}
\end{subequations}
where $\mathbb{E}[\cdot]$ is the expectation operator and superindex
$T$ means transpose. The system is controlled by the tandem
sensor--actuator ($\bS{}{n}$,$\bA{}{n}$) characterized by the
parametrization $\mybth$.
Let us define the control task as the modification of the moments of
$\bJ{}{n+1}$ in Eq. (\ref{eq:control:mu_sigma}) and denote by
$\bmutar$ and $\bsgtar$ the targeted (i.e., desired) mean and variance
for $\bJ{}{n+1}$, respectively.
The goal of the control is to drive the system to the optimal state
$\mybJ^*$, with mean $\bmuopt$ and variance $\bsgopt$, as close as
possible to $\bmutar$ and $\bsgtar$, respectively, using the
controller with optimal parameters $\mybth^*$.
\corr{The search for the optimal control parameters can be posed as
  the reduction in uncertainty of the controlled state given the
  information collected by the sensors and the action performed by the
  actuators.} The latter is formulated as the minimization of
Kullback-Leibler divergence between $\bJ{}{n+1}$ and an auxiliary
state, $\bJhat{}{}$, constructed by shifting the mean of $\bJ{}{n+1}$
to $\bmutar$ and scaling its variance to $\bsgtar$.
Then, the optimal information-theoretic controller is attained for
\begin{subequations}\label{eq:control:optimal_KL_ori}
  \begin{gather}
  \mybth^* = \arg\min_{\mybth} \ \mathrm{KL}( \bJ{}{n+1} ,\bJhat{}{} ), \\
  \bJhat{}{} = \shiftM\bJ{}{n+1} + \bs{b}, \label{eq:control:rescale} 
\end{gather}
\end{subequations}
with
\begin{align*}
\shiftmi_{ij}^2 &= \frac{\sgtar_{ij}}{\mathbb{E}[(\J{i}{n+1}-\mu_i)(\J{j}{n+1}-\mu_j)]}, \quad
\bs{b} = \bmutar - \shiftM\bmun.
\end{align*}
The role of the matrix $\shiftM$ and the vector $\bs{b}$ is to
transform the probability distribution of $\bJ{}{n+1}$ into the
distribution $\bJhat{}{}$ such that $\mathbb{E}[\bJhat{}{}] = \bmutar$
and $\mathrm{var}[\bJhat{}{}] = \bsgtar$.
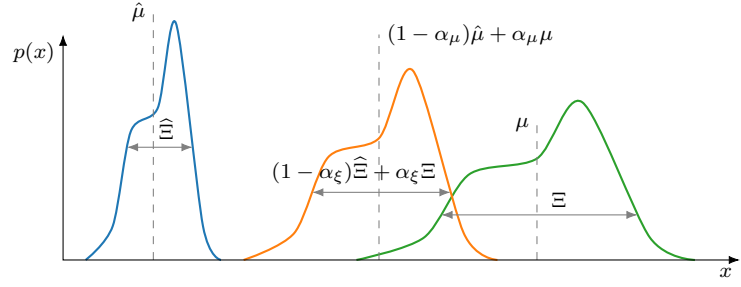
\begin{figure}
    \centering
    \pgfmathsetmacro{\mysca}{3}
\begin{tikzpicture}[scale=\mysca,>={Latex[length=.15cm]},
    sigma/.style={thin,<->,black!50,text=black},
    mean/.style={thin,dashed,black!50,text=black},
    pics/pdf/.style 2 args={code={\draw[scale=\mysca,name path global=#1,thick,color=#2] 
        plot [smooth] coordinates 
        {(0,0) (.3,0.1) (.5,.38) (.8,.45) (1.,.7) (1.3,.1) (1.5,0)};}}]

    \begin{footnotesize}

    \pic[name path=pO]  at (1.3,0) {pdf={pO}{C2}};

    \pic[xscale=.75,yscale=1.2] at (.8,0) {pdf={pI}{C1}};

    \pic[xscale=.4,yscale=1.5] at (0.1,0) {pdf={pT}{C0}};

    \path [name path=rectO] (0,.2) --+ (3,0); 
    \path [name intersections={of = pO and rectO}];
    \draw[sigma] (intersection-1) -- (intersection-2) node[above,pos=.6] {$\sgn$};
    \path [name path=rectI] (0,.3) --+ (3,0); 
    \path [name intersections={of = pI and rectI}];                    
    \draw[sigma] (intersection-1) -- (intersection-2) node[above,pos=.3]
        {$(1-\relfsg)\sgtar + \relfsg\sgn$};
    \path [name path=rectT] (0,.5) --+ (3,0); 
    \path [name intersections={of = pT and rectT}];                   
    \draw[sigma] (intersection-1) -- (intersection-2) node[above,pos=.6] {$\sgtar$};

    \draw[mean] (.4,0) --+ (0,1.1) node[anchor=east,pos=1] {$\mutar$};
    \draw[mean] (1.4,0) --+ (0,1) node[anchor=west,pos=1] 
    {$(1-\relfmu)\mutar + \relfmu\mu$};
    \draw[mean] (2.1,0) --+ (0,.6) node[anchor=east,pos=1] {$\mun$};

    \draw [<->] (0,1) -- (0,0) node[pos=0,anchor=north east] {$p(x)$} -- 
        (3,0) node[anchor=north east] {$x$};
        
    \end{footnotesize} 
\end{tikzpicture}
    \caption{Sketch of the of the probability mass functions illustrating the effect
    of the transformation in Eq.~\eqref{eq:control:rescale} and the relaxation
    factors, $\relfmu$ and $\relfsg$.}
\end{figure}

The optimization method posed in Eq. (\ref{eq:control:optimal_KL_ori})
can be simplified by decomposing the controller parameters space into
three independent sets $\mybth=[\bths \ \bthpa \ \bthaa]$, where
$\bths$ are the sensor parameters (mainly, the sensor locations),
$\bthpa$ are the \emph{passive} actuator parameters (i.e., those that
do not directly modify the state of the system, such as the actuator
locations), and $\bthaa$ are the \emph{active} actuator parameters
(that act on the state of the system, such as the actuator amplitude,
frequency, etc.).
The parameters $\bths$ and $\bthpa$ can be optimized to enhance
observability and controllability, respectively.
It was shown in \S \ref{subsec:control:obs_con} that observability and
controllability improve with the mutual information shared among the
state variable $\bJ{}{n}$ and the controller state.
With this insight, the optimization problem can be simplified as the
following iterative process:
\renewcommand{\labelenumi}{\arabic{enumi}.}
\begin{enumerate}
    \item For a given iteration, $i$, assume $\bthaa$ fixed and solve for $\bths$ and $\bthpa$,
    \begin{subequations}\label{eq:control:os}
    \begin{align}   
        \bths  &\leftarrow \arg\max_{\bths} \ I(\bJ{}{n};\bS{}{n}),\label{eq:control:os1}\\  
      \bthpa &\leftarrow \arg\max_{\bthpa} \ I(\bJ{}{n};\bA{}{n}). \label{eq:control:os2}
    \end{align}
    \end{subequations}
    The initial guess for $\bths$ and $\bthpa$ can be obtained from the 
    uncontrolled system, $\bthaa = \bs{0}$. 
    \item Using the values of $\bths$ and $\bthpa$ from the previous 
    optimization, solve for $\bthaa$ as 
    \begin{align}\label{eq:control:oa}
      \bthaa &\leftarrow \arg\min_{\bthaa} \ 
      \mathrm{KL}(\bJ{}{n+1},\bJhat{}{}), 
    \end{align}
    \item Repeat steps 1 and 2 while $\mathrm{KL}^i < \mathrm{KL}^{i-1}$, being 
    $\mathrm{KL}^i = \mathrm{KL}(\bJ{}{n+1},\bJhat{}{})$ computed at iteration $i$.
\end{enumerate}

To further facilitate the optimization, step 3 can be aided by
introducing relaxation factors such that,
\begin{align}\label{eq:control:rel}
  \mathbb{E}[\bJhat{}{}] =  (1-\relfmu) \bmutar + 
  \relfmu \bmun, \quad
  \mathrm{var}[\bJhat{}{}] = (1 - \relfsg) \bsgtar + \relfsg \bsgn.
\end{align}
within each iteration.
These relaxation factors are within the range $0 < \relfmu,\relfsg \leq 1$ 
during the optimization, taking the limits $\relfmu \rightarrow 0$ and 
$\relfsg \rightarrow 0$ as the iterations advance.
The optimization procedure outlined above is algorithmically
appealing, as it can be performed using gradient ascent/descent
methods such as in reinforcement learning.  One application of
Eqs. \eqref{eq:control:os} and \eqref{eq:control:oa} is discussed in
\S\ref{sec:control:example} for optimal control for drag reduction in
wall turbulence.

\corr{Other variants of the optimization problem can be formulated to
  accommodate different cost functionals. For example, we might be
  interested in completely specifying a targeted probability
  distribution for $\bJ{}{n+1}$ (i.e., control acting over all the
  moments of $\bJ{}{n+1}$).  In that case, the rescaling described in
  Eq.~\eqref{eq:control:rescale} is not needed and the optimization
  problem is posed as
\begin{equation}\label{eq:control:optimal_KL_ori_all_mom}
  \mybth^* = \arg\min_{\mybth} \ \mathrm{KL}( \bJ{}{n+1} ,\bJ{\mathrm{ref}}{} ), \\
\end{equation}
where $\bJ{\mathrm{ref}}{}$ is the prescribed (and thus known)
probability distribution that we aim to attain for the state
$\bJ{}{n+1}$. The optimization in
Eq. (\ref{eq:control:optimal_KL_ori_all_mom}) provides the control
parameters that minimize the error between the probability
distribution of $\bJ{}{n+1}$ and $\bJ{\mathrm{ref}}{}$ in terms of the
L$_1$ norm.}

\subsection{Application: Opposition control for drag reduction in  turbulent channel flows}
\label{sec:control:example}

The enhanced transport of mass, momentum, and heat by turbulent flows
has a significant impact on the design and performance of thermofluid
systems. As such, the need of efficient control strategies to
manipulate turbulent flows remains an ubiquitous task in numerous
engineering applications. Examples of turbulent flow control can be
found in drag reduction for airfoils and pipelines, and enhanced
mixing for combustor chambers and heat exchangers, to name a few
examples. Given its technological importance, flow control remains a
field of active research, and many control strategies have been devised
with different degree of success~\citep{cattafesta2011}.  In the
present section, we consider an application of opposition control for
drag reduction in a turbulent channel flow, where flow is actively
modified at the wall to attenuate turbulence and reduce
drag~\citep{choi1994,hammond1998}. The optimal control strategy is
found using the information-theoretic tools described in \S
\ref{subsec:control:optimal}.

\begin{figure}
    \subfloat[\label{fig:control:channel}]{\newlength{\ddelta}
\setlength{\ddelta}{2.6em}
\begin{tikzpicture}[>={Latex[length=.15cm]},scale=2.5,
    slant/.style={yslant=0.5,xslant=-1}]
    \begin{footnotesize}
    
    \pgfmathsetmacro{\xl}{1.5}
    \pgfmathsetmacro{\zl}{.4*\xl}
    \pgfmathsetmacro{\yl}{.6*\xl}

    \pgfmathsetmacro{\ysla}{.5}
    \pgfmathsetmacro{\xsla}{-1}

    \tikzset{
        walls/.style = {fill=black!40,draw=black!70,fill opacity=.9},
        obsp/.style = {C1,fill=C1,dashed,thick,fill opacity=.1}
    }
    \coordinate (p0) at (.8,.1); 

    \begin{scope}[ yshift=0,yslant=\ysla,xslant=\xsla]
        \draw[walls] (0,0) rectangle (\xl,\zl); 

        \draw[black,->] (p0) --++(0,-.3) node [anchor=south] {$z$};
        \draw[black,->] (p0) --++(.3,0) node [anchor=south] {$x$};

        \draw (0,\zl+.05) --+ (0,.1) coordinate [midway] (c1);
    \end{scope}

    \begin{scope}[yshift=\ddelta*.05,yslant=\ysla,xslant=\xsla]
        \draw[obsp] (0,0) rectangle (\xl,\zl);

        \coordinate (c3) at (0,.2);
    \end{scope}

    \begin{scope}[yshift=\ddelta*.95,yslant=\ysla,xslant=\xsla]
        \draw[obsp] (0,0) rectangle (\xl,\zl);
        \coordinate (c4) at (0,.2);
    \end{scope}

    \begin{scope}[yshift=\ddelta,yslant=\ysla,xslant=\xsla]
        \draw[walls] (0,0) rectangle (\xl,\zl);
        \draw (0,\zl+.05) --+ (0,.1) coordinate [midway] (c2);
    \end{scope}

    \begin{scope}[yshift=\ddelta/2,yslant=\ysla,xslant=\xsla]
        \draw[thick,C0,->] (-.4,.5) --+ (.3,0) node [pos=0,above,anchor=south east] 
        {$U_{\mathrm{bulk}}$};
    \end{scope}

    \draw[black,->] (p0) --++ (0,.4) node [anchor=west] {$y$};
    \draw[black,<->] (c1) -- (c2) node [midway,anchor=east] {$2\delta$};

   
    \coordinate (c5) at (-.8,-.2);
    \draw[->] (c5) -- (c4); 
    \draw[->] (c5) -- (c3) node [fill=white,pos=0] {sensing planes};

    \path (0,-.5em);
    \end{footnotesize}
\end{tikzpicture}}\hfill
    \subfloat[\label{fig:control:vcontrol}]{\begin{tikzpicture}[>={Latex[length=.15cm]}]
    \begin{footnotesize}
    \node[anchor=south west] (fig) at (0,0) {\ig[width=.5\tw]{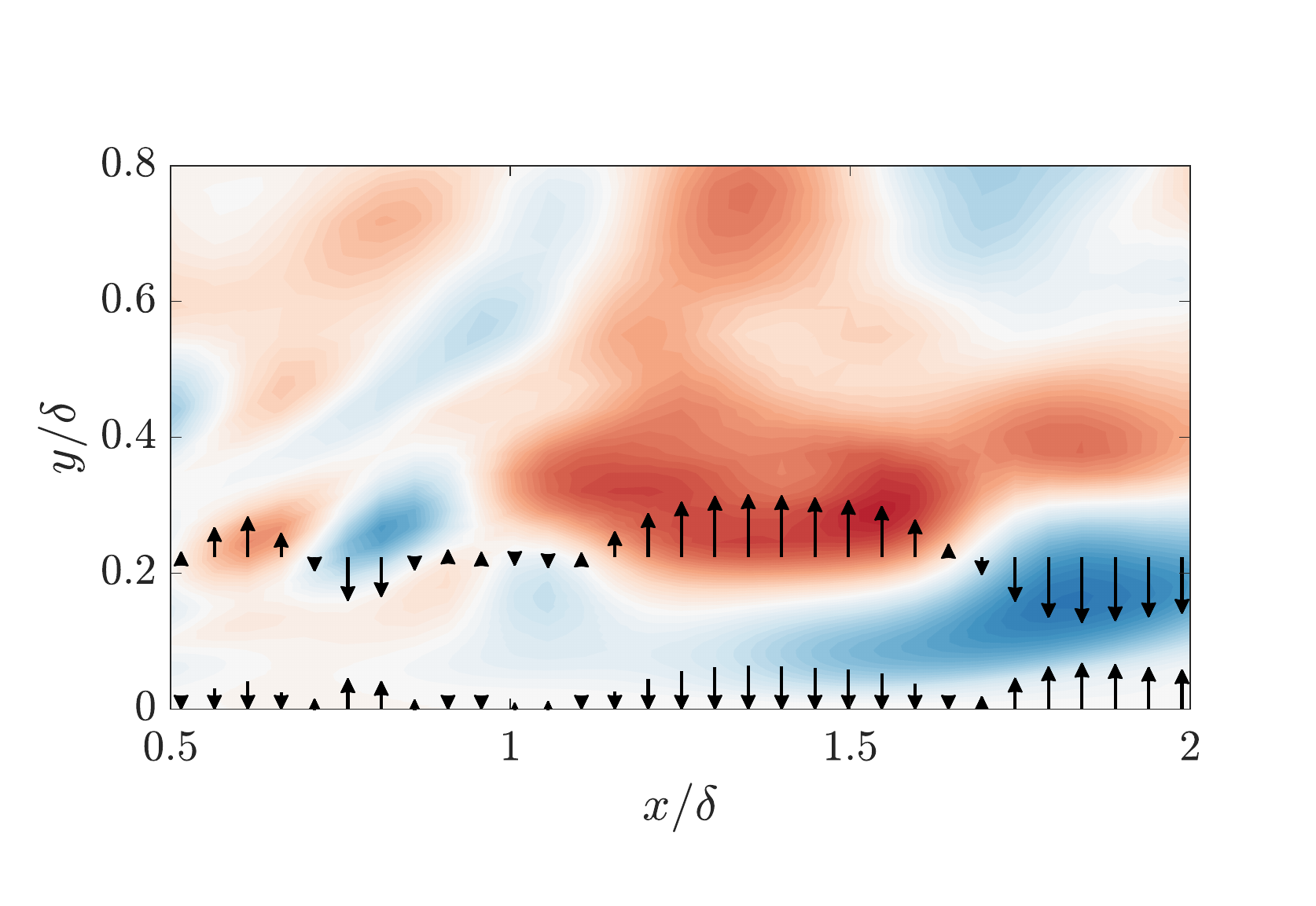}};
    \begin{scope}[x=(fig.south east),y=(fig.north west)]
        \coordinate (c0) at (.142,.39);
        \coordinate (c1) at (.891,.39);
        \draw[very thin,densely dashed] (c0) -- (c1);

        \draw[very thin] ($(c1)+(.02,0)$) --+ (.08,0)
                         (.911,.23) --+ (.08,0);

        \draw[very thin,<->] ($(c1)+(.06,0)$) --++ (0,-.16) 
            node [midway,anchor=west] {$y_s$};

        \node at (.35,.45) {$v(x,y_s,z)$};
    \end{scope}
    \end{footnotesize}
\end{tikzpicture}}
    \caption{(a) Sketch of a channel flow. The mean velocity is in the streamwise ($x$)
    direction. 
    (b) Schematic of the opposition control technique. The contour corresponds to
    the instantaneous vertical velocity on a $z$-plane. Colormap ranges from
    (red) $v^+ = -3.6$ to (blue) $3.6$.}
\end{figure}
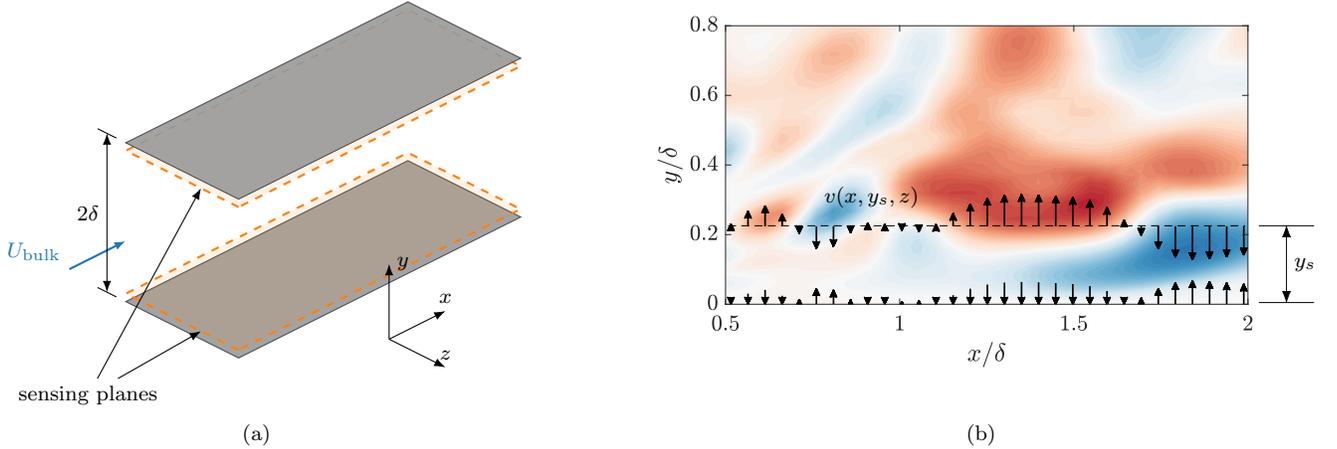

The flow configuration considered is an incompressible turbulent
channel flow~\citep[see section 7.1 in][]{pope2000} comprising the
flow confined between two parallel walls separated by a distance
$2\delta$ as shown in figure \ref{fig:control:channel}.  The
streamwise, wall-normal, and spanwise directions are denoted by $x$,
$y$, and $z$, respectively, and the corresponding velocities are $u$,
$v$, and $w$.  The flow is driven by imposing a constant mass flux in
the streamwise direction which is identical for both the uncontrolled
and controlled cases.  The bottom and top walls are located at
$y=0\delta$ and $y=2\delta$, respectively. The size of the
computational domain is $\pi\delta \times 2\delta \times \pi\delta/2$,
in the streamwise, wall-normal, and spanwise directions,
respectively. The Reynolds number is $\mathrm{Re} = U_{\mathrm{bulk}}
\delta/ \nu \approx 3200$, where $\nu$ is the kinematic viscosity and
$U_{\mathrm{bulk}}$ is the mean streamwise velocity.  The flow is
calculated by direct numerical simulation of the incompressible
Navier-Stokes equations in which all the scales of the flow are
resolved.  The code employed to perform the simulations was presented
and validated in previous studies \citep{lozano2016_brief, bae2018b,
  bae2019, lozano2021}.  In all the simulations, the domain is
discretized into $64\times90\times64$ grid points in the $x$, $y$, and
$z$ directions, respectively, which yields a total number of 368,640
degrees of freedom. The time step of the simulation is fixed to
$\Delta t^{+} \approx 5 \times 10^{-3}$, where $+$ denotes
non-dimensionalization by $\nu$ and the friction velocity $u_{\tau,u}
= \sqrt{ \nu \partial \langle u \rangle/\partial y}|_w$ for the
uncontrolled case. The operator $\langle \cdot \rangle$ signifies
average in $x$, $z$, and time, and the subscript $w$ denotes
quantities evaluated at the wall.

Opposition control is a drag reduction technique based on blowing and
sucking fluid at the wall with a velocity opposed to the velocity
measured at some distance from the wall. In the case of a channel
flow, the measured velocity is located in a wall-parallel plane at a
distance $y_s$ from the wall referred to as the \emph{sensing} plane.
Figure \ref{fig:control:vcontrol} provides an schematic of the problem
setup for opposition control in a turbulent channel flow.  The
optimization problem consists of finding the wall-normal distance of
the \emph{sensing} plane and the blowing/suction velocity of the
actuator.  We consider a prescribed law for the actuator such that the
blowing/suction velocity at the wall is proportional to the measured
velocity in the sensing plane.  The instantaneous wall-normal velocity
at the wall in the controlled case is given by
\begin{equation}
  \label{eq:control:law_oc}
    v(x,0,z) = -\beta v(x,y_{s},z), 
\end{equation}
where $\beta$ is the blowing intensity.  The controller parameter
vector is $\mybth = [\ths,\thaa] = [y_s,\beta]$, and
Eqs.~\eqref{eq:sensor_Sn} and \eqref{eq:controller_An} take the form
of $\myindexvar{S}{}{k} = v(x,y_{s},z)$ and $\myindexvar{A}{}{k} =
-\beta \myindexvar{S}{}{k}$, respectively, at times $t_k=t_{n}$ and
$t_{k}=t_{n+1}$.  A control law equivalent to
Eq. (\ref{eq:control:law_oc}) is applied at the top wall.

The quantity to be controlled is the mean wall shear stress (i.e., the
drag) in the statistically steady state of the controlled channel
flow, denoted by $\J{}{n+1} = \tau_w$. The analysis is conducted
considering two states: the state with no actuation, $\J{}{n}$, and
the final statistically steady state, $\J{}{n+1}$, after actuation has
been applied for a period of time equal to $\Delta t^+ = 60$. The
targeted mean and standard deviation of $\J{}{n+1}$ in the controlled
state are set to $\mutar^+ = 0$ and $\sgtar^{1/2} \approx 0.1
\langle\tau_{w,u}\rangle$, where $\langle\tau_{w,u}\rangle$ is the
mean wall-shear stress of the uncontrolled case.  The auxiliary
probability distribution $\Jhat{}{n+1}$ is defined as in
Eq.~\eqref{eq:control:rescale} using $\mutar$ and $\sgtar$.

Several methods are available for the optimization problem posed in
\S\ref{subsec:control:optimal}.  In our case, a simply gradient
descent algorithm was sufficient to find the optimal state.  The
optimum parameters in Eq.~\eqref{eq:control:os1} and
Eq.~\eqref{eq:control:oa} are respectively computed iteratively as
\begin{align*}
    \ths \leftarrow \ths + \gamma \nabla I(\myindexvar{S}{}{n+1};\myindexvar{J}{}{n+1}), \\
    \thaa \leftarrow \thaa - \gamma \nabla \mathrm{KL}(\myindexvar{J}{}{n+1};\Jhat{}{}),
\end{align*}
where here $n+1$ represents successive controlled states, $\gamma$ is
the step size computed as in \citet{barzilai1988}, and the gradient is
numerically computed using forward finite differences.  The parametric
space is bounded by $\theta_s \equiv y_s^+ \in [0, 180]$ (for the
bottom wall) and $\theta_a \equiv \beta \in [-0.1, 1.1]$. Values of
$\beta$ larger than $1$ were found to be unstable, consistently with
the findings in \citet{chung2011}.  The iteration process is started
by finding the optimal sensor in the uncontrolled state and setting
$\beta = 0$.  The best location of the sensing plane before actuation
is found to be $y_s^+ \approx 9.65$.  The iterative process is then
continued following the steps in \S\ref{subsec:control:optimal}.  To
aid the optimization, $\relfmu$ and $\relfsg$ in
Eq.~\eqref{eq:control:rel} are initially set to $0.6$ and gradually
decreased within each iteration.  The optimal control is found at
$y_s^{+*} \approx 13.9$ and $\beta^* = 1$.
Remarkably, our optimal control coincides with the global optimum reported by \citet{chung2011},
who performed a parametric study varying $y_s$ and $\beta$  in a flow setup identical to the one presented here.

In terms of drag, the optimal control provides $\approx 26\%$
reduction with respect to the uncontrolled state.  A similar value was
reported by other authors \citep{choi1994,hammond1998,chung2011}, 
which is expected since our control parameters are similar.  However, it is important 
to remark that the methodology adopted here radically differs from the approach followed
in the aforementioned references:
instead of conducting simulations using a trial and error approach to
find the optimal parameters, here we rely on the information-theoretic
principles introduced in \S\ref{sec:control}.

For a more quantitative perspective, Figure~\ref{fig:control:tau}
displays the probability mass distribution of the wall shear-stress
for the uncontrolled state and optimally controlled state.  It can be
readily seen that both the mean and the standard deviation are smaller
for the case with the optimal control.  Similarly,
figure~\ref{fig:control:IKL} (left) shows that the Kullback-Leibler
divergence is lower for the controlled case with optimal parameters
than for the uncontrolled case.
Figure~\ref{fig:control:IKL} (right) depicts the mutual information
between the wall shear stress and the wall-normal velocity
at the \emph{sensing} plane for the optimal sensor location,
$y_s^{*,k}$ for the uncontrolled state ($k = n$) and controlled state
($k = n+1$).  Interestingly, the mutual information is larger for the
actuated state, arguably because the wall shear stress in the actuated
case is more correlated with the imposed velocity at the
wall.

\definecolor{fig14c0}{RGB}{125,187,205}
\definecolor{fig14c1}{RGB}{240,148,94}
\begin{figure}
    \subfloat[\label{fig:control:tau}]{\ig[width=.5\tw]{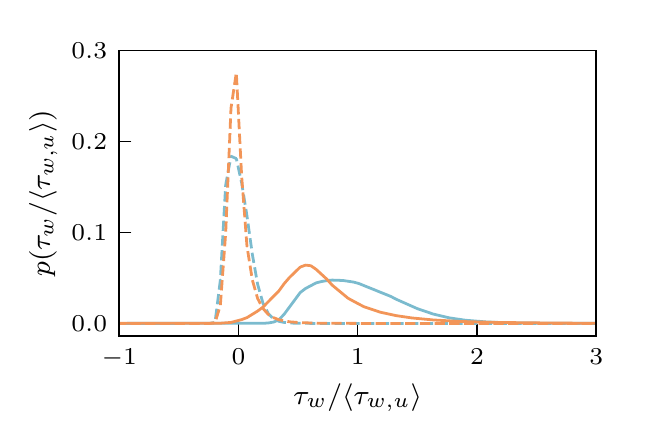}} \hfill
    \subfloat[\label{fig:control:IKL}]{\ig[width=.5\tw]{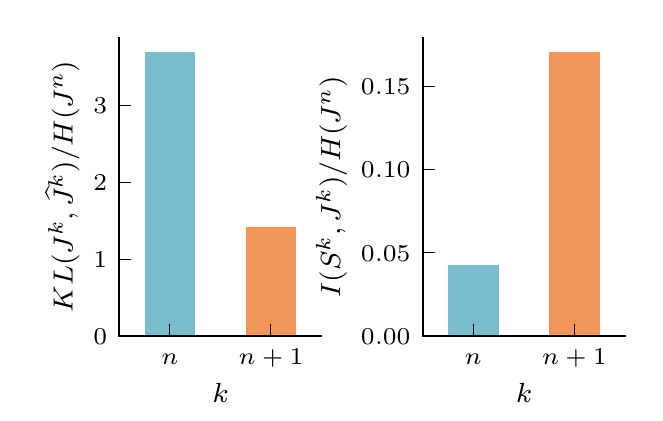}}
    \caption{(a) Probability mass distributions of the wall shear
      stress for the \myls{-}{fig14c0} uncontrolled state ($k = n$)
      and the \myls{-}{fig14c1} controlled state ($k = n+1$).  The
      line styles are \myls{-}{gray} for actual wall shear stress
      distribution ($J^k$) and \myls{--}{gray} for the auxiliary state
      ($\protect\Jhat{}{k}$).
(b) KL divergence between the final state and the auxiliary state
      (left); and Mutual information between the sensor location and
      the state (right), normalized with the entropy of the uncontrolled state.}
\end{figure}

\begin{figure}
    \subfloat[]{\ig[width=.5\tw,trim=.5cm 1.5cm 0cm 2cm,clip]{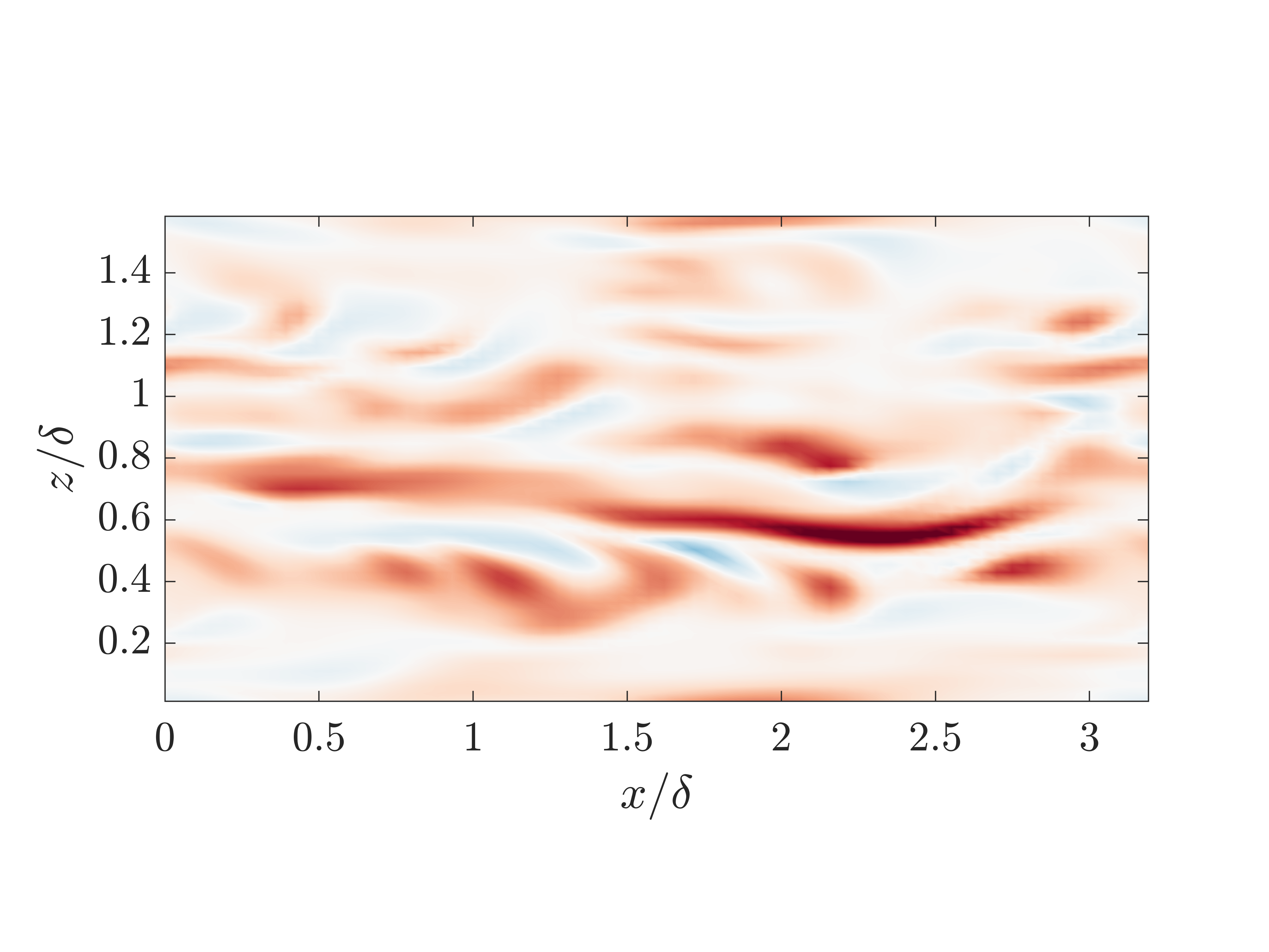}} \hfill
    \subfloat[]{\ig[width=.5\tw,trim=.5cm 1.5cm 0cm 2cm,clip]{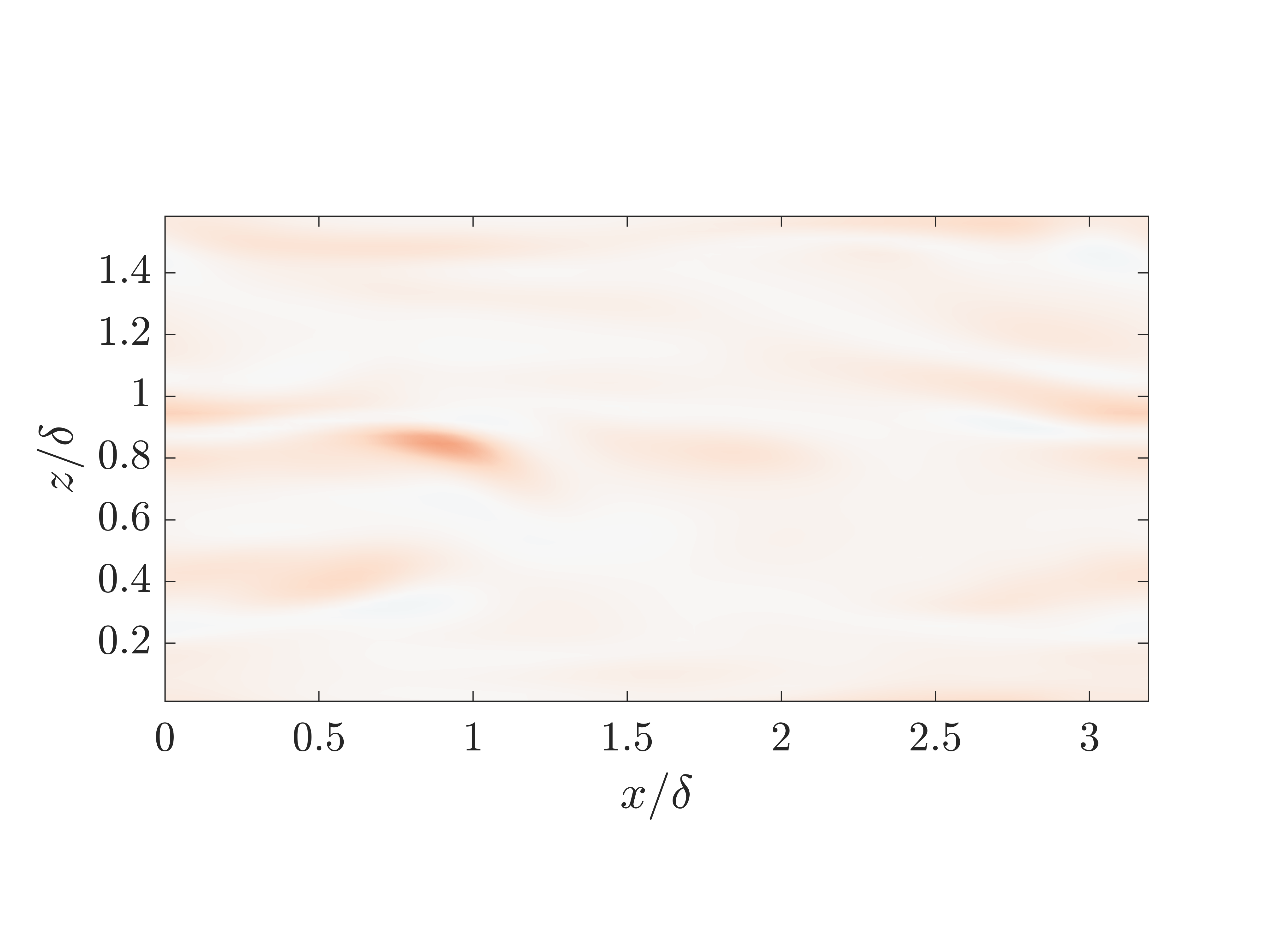}} \hfill
    \caption{Instantaneous tangential Reynolds stress $u'v'$ at $y^+
      \approx 10$ for (a) uncontrolled state and (b) optimally
      controlled state.  The colormap ranges from $u'v'/u_{\tau,u}^2 =
      -6.8$ (blue) to $6.8$ (red).
    \label{fig:control:rey}}
\end{figure}

Finally, to provide additional insight into the effect of the actuation on the
flow, figure~\ref{fig:control:rey} shows the tangential Reynolds
stresses $u'v'$ (where prime denotes fluctuations about the mean
value) at $y_s^+ \approx 10$ for the uncontrolled and  controlled
states.  It can be appreciated that the intensity of the tangential
Reynolds stresses is lower for the case with optimal control,
indicating a suppression of vortical structures near the wall.
%

%
%



\section{Conclusions}

The problems of causality, reduced-order modeling, and control for
chaotic, high-dimensional dynamical systems have been formulated
within the framework of information theory. In the proposed formalism,
the state of the dynamical system is considered a random variable in
which its information (i.e., Shannon entropy) quantifies the average
number of bits to univocally determinate its value.  A key quantity
for the formalization of the theory is the conditional information,
which measures the uncertainty in the state of the system given the
partial knowledge of other states. In contrast to the
equation-centered formulation of dynamical systems, where individual
trajectories are the object of analysis, information theory offers a
more natural approach to the investigation of chaotic systems from the
viewpoint of the probability distributions of the states.

We have argued that statistical asymmetries in the information flux
within the states of the system can be leveraged to measure causality
among variables. As such, the information-theoretic causality from one
variable to another is quantified as the information flux from the
former to the latter. Our definition of causality is motivated by the
information required to attain total knowledge of a future state and
can be interpreted as how much the past information of the system
improves our knowledge of the future state. The formulation of
causality proposed here is grounded on the zero conditional-entropy
condition for deterministic systems and generalizes \corr{to multiple
  variables the definition of causality by \citet{schreiber2000}. The
  quantification of causality proposed also accounts for the
  information flux due to the joint effect of variables, which was
  absent in previous formulations. We have also introduced the
  information leak as the amount of information unaccounted for by the
  observable variables.}

Reduced-order modeling of chaotic systems has been posed as a problem
of conservation of information: modeled systems contain a smaller
number of degrees of freedom than the original system, which in turn
entails a loss of information. Thus, the primary goal of modeling is
to preserve the maximum amount of useful information from the original
system. \corr{We have derived the conditions for maximum
  information-preserving models and shown that accurate models must
  maximize the mutual information between the model state and the true
  state, and minimize the Kullback-Leibler divergence between their
  probabilities}. The mutual information assists the model to
reproduce the dynamics of the original system, while the
Kullback-Leibler divergence enables the accurate prediction of the
statistical quantities of interest.

Lastly, control theory has been cast in information-theoretic terms by
envisioning the controller as a device aimed at reducing the
uncertainty in the future state of the system to be controlled given
the information collected by the sensors and the action performed by
the actuators. We have reformulated the concepts of controllability
and observability using mutual information between the present and
future states. The definitions of open- and closed-loop control have
also been introduced based on the information shared between the
actuator and the system state. The optimization problem was posed as
the minimization of the Kullback-Leibler divergence between the
probability distribution of the controlled state and a targeted state
derived from the latter.

We have applied our information-theoretic framework to advance three
outstanding problems in the causality, modeling, and control of
turbulent flows.  Information-theoretic causal inference was used to
measure the information flux of the turbulent energy cascade in
isotropic turbulence. The principle of maximum conservation of
information was leveraged to devise a subgrid-scale model for
large-eddy simulation of isotropic turbulence.  Finally,
information-theoretic control was utilized to achieve optimal drag
reduction in wall-bounded turbulence using opposition control at the
wall.
%
Overall, information theory offers an elegant formalization of the
problems of causality, modeling, and control for chaotic,
high-dimensional systems, aiding physical interpretation and easing
the tasks of modeling and control all within one unified framework.

\section{Acknowledgments}

This work was supported by the National Science Foundation under Grant
No.  032707-00001. G.~A. was partially supported by STTR with Cascade
Technologies, Inc. and the Naval Air Systems Command.  The authors
acknowledge the MIT SuperCloud and Lincoln Laboratory Supercomputing
Center for providing HPC resources that have contributed to the
research results reported within this paper.

\bibliographystyle{apsrev4-2}
\bibliography{references}

\end{document}